\documentclass[a4paper,times,3p]{elsarticle}
\usepackage[utf8]{inputenc}
\usepackage[T1]{fontenc}

\usepackage{microtype}

\usepackage{booktabs}

\usepackage{amsmath,amsfonts,amsthm,enumerate,graphicx,epstopdf,url,etoolbox,algorithm,algorithmic,listings}
\usepackage{caption}
\usepackage{diagbox}
\usepackage{graphbox}
\usepackage{enumitem}
\usepackage{dirtytalk}
\usepackage{tikz}
\usepackage[export]{adjustbox}
\usepackage{subfig}
\usetikzlibrary{shapes, arrows}
\usepackage{tcolorbox}

\definecolor{cmykcyan}{cmyk}{1,0,0,0}
\definecolor{cmykred}{cmyk}{0,1,1,0}
\definecolor{cmykblack}{cmyk}{0,0,0,1}

\newtheorem*{remark}{Remark}

\usepackage[font=footnotesize,labelfont=bf]{caption}
\usepackage{multirow}

\usepackage{pxfonts}

\usepackage{import}

\graphicspath{./graphics/}

\usepackage[textsize=scriptsize]{todonotes}
\usepackage{setspace}

\makeatletter
\def\ps@pprintTitle{%
 \let\@oddhead\@empty
 \let\@evenhead\@empty
 \def\@oddfoot{}%
 \let\@evenfoot\@oddfoot}
\makeatother

\journal{Computer Methods in Applied Mechanics and Engineering}

\begin{document}

\begin{frontmatter}

\title{Goal-Oriented Adaptive THB-Spline Schemes for PDE-Based Planar Parameterization}

\author[add1]{Jochen Hinz\corref{cor1}}
\ead{j.p.hinz@tudelft.nl}
\author[add3]{Michael Abdelmalik\corref{cor1}}
\ead{michael.amalik@austin.utexas.edu}
\author[add2]{Matthias M\"oller\corref{}}
\ead{m.moller@tudelft.nl}

\cortext[cor1]{Corresponding author}

\address[add1]{Chair of Numerical Modelling and Simulation, Ecole Polytechnique F\'ed\'erale de Lausanne, 1018 Lausanne, Switzerland.}
\address[add3]{Delft Institute of Applied Mathematics, Delft University of Technology, 2628 XE Delft, Netherlands.}
\address[add2]{Oden Institute, University of Texas at Austin, 78712 Austin, Texas, US.}

\begin{abstract}
This paper presents a PDE-based planar parameterization framework with support for \textit{Truncated Hierarchical B-Splines} (THB-splines) which approximates an inversely harmonic geometry parameterization given no more than a boundary correspondence between physical and parametric domains. We accomplish this by requiring that the mapping satisfy the equations of \textit{Elliptic Grid Generation} (EGG) and present an adaptive numerical scheme that can guarantee analysis-suitability by following the structure of the underlying PDE-problem. This is accomplished by sufficiently accurately approximating the exact (folding-free) PDE-solution. For this, we adopt the a posteriori refinement strategy of \textit{Dual Weighted Residual} (DWR) and combine it with goal-oriented cost functions to warrant bijectivity as well as parameterization quality. Hereby, the combination of goal-oriented a posteriori refinement strategies and THB-enabled local refinement avoids over-refinement, in particular in geometries with complex boundaries. \\
To control the parametric properties of the outcome, we introduce the concept of domain optimization. Hereby, the properties of the domain into which the mapping maps inversely harmonically, are optimized in order to fine-tune the parametric properties of the recomputed geometry parameterization.
\end{abstract}

\end{frontmatter}

\section{Introduction}
\label{sect:Adaptive_Introduction}
Isogeometric analysis (IGA), first introduced by Hughes et al. in \cite{cottrell2009isogeometric}, is a numerical technique that aims to bridge the gap between computer aided design (CAD) and (isoparametric) finite element analysis (FEA). This is accomplished by building the geometry mapping from the same spline basis that is used to approximately solve PDE-problems posed over the geometry. As such, spline-based parameterization techniques have received an increased amount of interest in the mathematical community in recent years. Since the CAD pipeline typically provides no more than a spline-based description of the boundary contours of the target geometry, the purpose of all parameterization algorithms is to generate a bijective (folding-free) geometry parameterization from the boundary CAD data. Analogous to mesh quality in classical FEA, the parametric quality of the surface parameterization has a profound impact on the numerical accuracy of the isogeometric analysis \cite{xu2010optimal}. Therefore, besides bijectivity, proficient parameterization algorithms aim at generating parameterizations of high numerical quality. \\
One of the most important applications of IGA lies in shape optimization problems. Since the geometry changes at every shape optimization iteration, algorithms that are differentiable with respect to the design variables (i.e., the boundary control points) have a further advantage since they allow for employing gradient-based shape optimization algorithms which tend to converge in fewer iterations than their zeroth-order counterparts. Another advantage of differentiability is efficiency: as the inner control points are a smooth function of the boundary control points, there is no need for full remeshing after each iteration since cheaper mesh update strategies can be employed. This is also true for settings in which the boundary contours change as a smooth function of time. \\
Traditionally, parameterizations for IGA-applications are built from tensor-product spline spaces. Unfortunately, structured spline technologies do not allow for local refinement as knot insertion in one parametric direction automatically refines a whole row / column of the underlying spline space. For the geometry description, this may result in a very dense spline basis whenever many degrees of freedom (DOFs) are required to properly resolve the boundary contours. As a result, the total number of unknowns (the inner control points) may become infeasibly large, leading to a severe slow-down of the meshing process and / or the isogeometric analysis. \\
To address above efficiency concerns, this paper introduces a PDE-based planar parameterization framework that uses THB-splines \cite{giannelli2012thb}, an unstructured spline technology which allows for local refinement, potentially reducing the required total number of DOFs. A major challenge of unstructured spline technologies is deciding where local refinement is required and where a lower resolution suffices. For this, we employ the principles of \textit{dual weighted residual} (DWR) \cite{rannacher2004adaptive}, an a posteriori refinement technique for PDE problems based on duality considerations. Furthermore, we augment the problem formulation with a mechanism that allows for changing the parametric properties of the PDE solution in order to fine-tune the parametric properties of the mapping operator.

\subsection{Notation}
In this work, we denote vectors in boldface. The $i$-th entry of a vector $\mathbf{x}$ is denoted by $\mathbf{x}_i$ or simply $x_i$ and similarly for the $ij$-th entry of matrices. We make extensive use of vector derivatives. Here, we interchangeably use the denotation
\begin{align}
    \partial_{\mathbf{t}} \mathbf{x} \equiv \frac{\partial \mathbf{x}}{\partial \mathbf{t}}, \quad \text{with} \quad \left[ \frac{\partial \mathbf{x}}{\partial \mathbf{t}} \right]_{ij} = \frac{\partial x_i}{\partial t_j}
\end{align}
for the partial derivative. \\
Furthermore, we frequently work with spline vector spaces $\mathcal{V}_h$. Here, $\left[ \mathcal{V}_h \right]$ refers to the canonical (THB-) spline basis of $\mathcal{V}_h$, which we assume to be clear from context. By default, we employ the abuse of notation
\begin{align}
    \left( \mathcal{V}_h \right)^n = \underbrace{ \mathcal{V}_h \times \cdots \times \mathcal{V}_h }_{n \text{ terms}}.
\end{align}
For better readability, we avoid the parenthesis when no confusion is possible, i.e., $(\mathcal{V}_h)^2 = \mathcal{V}_h^2$.

\subsection{Problem Statement}
Let $\Omega$ denote the target geometry and $\hat{\Omega} = (0, 1)^2$ the parametric domain. In general, we assume that $\Omega$ is topologically equivalent to $\hat{\Omega}$. By $\mathbf{x}: \hat{\Omega} \rightarrow \Omega$, we denote the mapping operator whose components are built from the linear span of the (THB-)spline basis $\left[ \mathcal{V}_h \right] = \{w_1, \ldots, w_N\}$. 
The mapping operator $\mathbf{x}: \hat{\Omega} \rightarrow \Omega$ is of the form:
\begin{align}
\label{eq:Adaptive_mapping}
    \mathbf{x}(\xi, \eta) = \sum_{i \in \mathcal{I}_{I}} \mathbf{c}_i w_i(\xi, \eta) + \sum_{j \in \mathcal{I}_{B}} \mathbf{c}_j w_j(\xi, \eta),
\end{align}
where $\mathcal{I}_{I}$ and $\mathcal{I}_{B}$ refer to the index-sets corresponding to vanishing and nonvanishing basis functions on $\partial \hat{\Omega}$, respectively and $\mathbf{c}_k \in \mathbb{R}^2, \enskip \forall k \in \mathcal{I}_{I} \cup \mathcal{I}_{B}$. Here, $\mathcal{I}_{I}$ corresponds to the subspace $\mathcal{V}_h^\circ = \mathcal{V}_h \cap H^1_0(\hat{\Omega})$. Note that $\mathcal{I}_{I}$ and $\mathcal{I}_{B}$ are mutually disjoint and
\begin{align*}
    \mathcal{I}_{I} \cup \mathcal{I}_{B} = \{1, \ldots, N\}.
\end{align*}
In general, we assume that the $\mathbf{c}_j$ in (\ref{eq:Adaptive_mapping}) are chosen such that $\mathbf{x}\vert_{\hat{\Omega}}$ is a Jordan curve that parameterizes $\partial \Omega$. \\
Then, the purpose of any parameterization algorithm is to choose the $\mathbf{c}_i$ in (\ref{eq:Adaptive_mapping}) such that
\begin{enumerate}
    \item $\mathbf{x}: \hat{\Omega} \rightarrow \Omega$ is bijective,
    \item $\mathbf{x}$ is a parameterization of high numerical quality,
\end{enumerate}
while the $\mathbf{c}_j$ are typically held fixed. The somewhat loosely defined notion of \textit{numerical quality} from point (2.) is a major difficulty in parameterization problems: as the assessment the numerical quality of a parameterization is problem-depended and generally only possible after the it has been completed, a priori quality criteria are inherently heuristic. The commonly applied heuristics will be discussed in Section \ref{sect:Adaptive_related_work}.

\subsection{Related Work}
\label{sect:Adaptive_related_work}
Existing parameterization techniques can be divided into three broad categories:
\begin{enumerate}
    \item Algebraic (direct) methods;
    \item methods based on (constrained and unconstrained) quality cost function optimization;
    \item PDE-based methods.
\end{enumerate}
Algebraic methods (1.) generate a mapping from the solution of a linear system of equations or the evaluation of a closed-form expression. The most-widely used algebraic method is based on the Coon's patch approach \cite{farin1999discrete}. Given the four (known) boundary curves $\mathbf{x}(\xi, 0)$, $\mathbf{x}(1, \eta)$, $\mathbf{x}(\xi, 1)$ and $\mathbf{x}(0, \eta)$, the mapping is constructed by projecting the components of
\begin{align}
\label{eq:Adaptive_Coons}
    \mathbf{x}_{\text{Coons}} = & \enskip (1 - \xi) \mathbf{x}(0, \eta) + \xi \mathbf{x}(1, \eta) \nonumber \\
                              + & \enskip (1 - \eta) \mathbf{x}(\xi, 0) + \eta \mathbf{x}(\xi, 1) \nonumber \\
                              - & \enskip \begin{bmatrix} 1 - \xi & \xi \end{bmatrix} \begin{bmatrix} \mathbf{x}(0, 0) & \mathbf{x}(0, 1) \\ \mathbf{x}(1, 0) & \mathbf{x}(1, 1) \end{bmatrix} \begin{bmatrix} 1 - \eta \\ \eta \end{bmatrix}
\end{align}
onto the spline space $\mathcal{V}_h$. Whenever $\left[ \mathcal{V}_h \right]$ is a tensor-product spline basis, the inner control points can also be computationally inexpensively computed from an explicit formula, see \cite{farin1999discrete}, while in an unstructured setting equation (\ref{eq:Adaptive_Coons}) can be used. \\
Another class of algebraic methods results from minimizing a convex, quadratic cost function $Q(\mathbf{x})$ over the inner control points $\mathbf{c}_i, i \in \mathcal{I}_I$. As before, the boundary control points follow from the boundary contours and are held fixed. $Q(\mathbf{x})$ is typically given by a positively-weighted sum of several cost functions. As such, it takes the form:
\begin{align}
\label{eq:Adaptive_weighted_sum}
    Q(\mathbf{x}) = \sum_i \underbrace{\lambda_i}_{\geq 0} Q_i(\mathbf{x}),
\end{align}
while the minimization problem becomes:
\begin{align}
\label{eq:Adaptive_minimization_quadratic}
    \int_{\hat{\Omega}} Q(\mathbf{x}) \mathrm{d}S \rightarrow \min_{\mathbf{x} \in \mathcal{V}_h^2}, \quad \text{s.t.} \quad \mathbf{x} \vert_{\partial \hat{\Omega}} & = \partial \Omega.
\end{align}
Possible choices for the $Q_i(\mathbf{x})$ in (\ref{eq:Adaptive_weighted_sum}) are \cite{falini2015planar}:
\begin{align}
\label{eq:Adaptive_length}
    Q_{\text{length}}(\mathbf{x}) = \left \| \mathbf{x}_\xi \right \|^2 + \left \| \mathbf{x}_\eta \right \|^2 \quad \text{and} \quad Q_{\text{uniformity}}(\mathbf{x}) = \left \| \mathbf{x}_{\xi \xi} \right \|^2 + 2 \left \| \mathbf{x}_{\xi \eta} \right \|^2 + \left \| \mathbf{x}_{\eta \eta} \right \|^2,
\end{align}
where the latter requires $\mathcal{V}_h \subset C^1(\hat{\Omega})$. The minimization of (\ref{eq:Adaptive_minimization_quadratic}) with the aforementioned quadratic cost functions converges after one iteration of a Newton-type optimization algorithm and can hence be considered of type (1.) as well as type (2.). For an overview of type (1.) approaches, we refer to \cite{gravesen2012planar}. \\
Another convex but quartic cost function is the \textit{Liao}-functional \cite{steinberg1993fundamentals}
\begin{align}
    Q_{\text{Liao}} = g_{11}^2 + 2 g_{12}^2 + g_{22}^2,
\end{align}
where the $g_{ij}$ denote the entries of the metric tensor of the mapping, with
\begin{align}
\label{eq:Adaptive_metric_tensor}
    g_{ij} = \mathbf{x}_{\boldsymbol{\xi}_i} \cdot \mathbf{x}_{\boldsymbol{\xi}_j} \quad \text{and} \quad \boldsymbol{\xi} = (\xi_1, \xi_2)^T \equiv (\xi, \eta)^T.
\end{align}
The minimization of above cost functions is computationally efficient, thanks to convexity, however, the resulting mappings are often folded, i.e., they do not satisfy:
\begin{align}
\label{eq:Adaptive_detJ_positive}
    \det J > 0, \quad \forall (\xi, \eta)^T \in \hat{\Omega}, \quad \text{where} \quad J(\mathbf{x}) = \partial_{\boldsymbol{\xi}} \mathbf{x}
\end{align}
denotes the Jacobian of $\mathbf{x}$. \\
The minimization of nonconvex quality functionals is computationally more demanding but tends to yield better results when convex optimization leads to a folded mapping \cite{steinberg1993fundamentals}. Typical nonconvex quality functionals are:
\begin{itemize}
    \item The \textit{area} functional
    \begin{align}
        \label{eq:Adaptive_Area_costfunc}
        Q_{\text{area}} = \left( \det J \right)^2,
    \end{align}
    which aims to minimize the variance of $\det J$ over $\hat{\Omega}$;
    \item the \textit{orthogonality} functional
    \begin{align}
        Q_{\text{Orthogonality}} = g_{12}^2 \quad \text{or} \quad  Q_{\text{AreaOrthogonality}} = g_{11} g_{22},
    \end{align}
    which is aimed at orthogonalizing the parameter lines;
    \item the \textit{eccentricity} functional
    \begin{align}
        Q_{\text{eccen}} = \left( \frac{ \mathbf{x}_{\xi} \cdot \mathbf{x}_{\xi \xi} }{g_{11}} \right)^2 + \left( \frac{ \mathbf{x}_{\eta} \cdot \mathbf{x}_{\eta \eta} }{g_{22}} \right)^2,
    \end{align}
    which penalizes fast accelerations along the parameter lines.
\end{itemize}
Unfortunately, minimization of the above functionals, in many cases, leads to folding, too. To the best of our knowledge, there are two main ways to prevent the grid from folding:
\begin{enumerate}[label=({\alph*})]
    \item Penalization;
    \item constrained minimization.
\end{enumerate}
Option (a) attempts to prevent grid folding through the modification of existing cost functions with a penalty term, such as
\begin{itemize}
    \item the \textit{Modified Liao} functional
    \begin{align}
        Q_{\text{ML}} = \left( \frac{g_{11} + g_{22}}{\det J} \right)^2.
    \end{align}
\end{itemize}
Adding the Jacobian determinant in the denominator serves the purpose of mitigating the tendency to fold, since the cost functional possesses an infinite barrier close to the boundary of the feasible region. \\
The most widely-used penalty cost functional is the so-called
\begin{itemize}
    \item \textit{Winslow} functional
    \begin{align}
        \label{eq:Adaptive_Winslow}
        Q_{\text{W}} & = \frac{g_{11} + g_{22}}{\det J}.
    \end{align}
\end{itemize}
With $\mathbf{x} \equiv (x, y)^T$, the \textit{Winslow} functional (\ref{eq:Adaptive_Winslow}) follows from performing a pullback of the problem
\begin{align}
\label{eq:Adaptive_inverse_laplace}
    \frac{1}{2} \int_{\Omega} \left \| \boldsymbol{\xi}_x \right \|^2 + \left \| \boldsymbol{\xi}_y \right \|^2 \mathrm{d} \mathbf{x} \rightarrow \min_{\mathbf{x}}, \quad \text{s.t.} \quad \boldsymbol{\xi} \vert_{\partial \Omega} = \partial \hat{\Omega}
\end{align}
into $\hat{\Omega}$. For details, we refer to \cite{gravesen2012planar}. An approach based on the \textit{Winslow} functional can be regarded as the \textit{mapping inverse} counterpart of an approach based on the length functional (\ref{eq:Adaptive_length}). \\ 
In option (b), the minimization is carried out with an added constraint that constitutes a sufficient condition for (\ref{eq:Adaptive_detJ_positive}). For tensor-product B-spline bases, in \cite{xu2011parameterization}, Xu et al. propose a linear convex sufficient condition $L(\mathbf{x}) > 0$ for bijectivity. It is added as a constraint to the minimization problem. If convex cost functions are utilized, this leads to a linear programming problem, which can be computationally inexpensively solved using convex optimization routines. Unfortunately, the set 
\begin{align*}
\left \{ \mathbf{x} \in \mathcal{V}_h^2 \enskip \vert \enskip \mathbf{x}\vert_{\partial \hat{\Omega}} = \partial \Omega \text{ and } L(\mathbf{x}) > 0 \right \}
\end{align*}
may be empty or the constraint may be very restrictive, limiting its applicability to relatively simple shapes. \\
In an effort to allow for more complicated shapes, \cite{xu2011parameterization} and \cite{gravesen2012planar} propose nonlinear nonconvex sufficient conditions for bijectivity. Since the Jacobian determinant $\det J$ is a piecewise-polynomial function of higher polynomial degree itself, it can be projected onto a spline basis that contains it. If all the weights are positive under the expansion, this constitutes a sufficient condition for bijectivity. The nonlinear sufficient condition $N(\mathbf{x}) > 0$ is added as a constraint and the optimization is carried out with a blackbox nonlinear optimization routine (typically, IPOPT \cite{biegler2009large}) that comes with all the drawbacks of nonconvex optimization such as the danger of getting stuck in local minima. A further disadvantage is the need for an initial guess that satisfies the constraints, for which another nonconvex optimization problem has to be solved first. \\ 
While the extension of (penalized or unpenalized) cost function minimization to THB-splines is straightforward, this is not the case for constrained methods, since the constraints are designed for structured splines only. To the best of our knowledge, the only comprehensive overview of planar parameterization techniques for THB-splines can be found in \cite{falini2015planar}, where the application of most of the mentioned (unpenalized) cost functions is studied in a THB-setting. As the optimization is carried out without constraints, folding occurs in the majority of test cases. The paper concludes that the only method potentially capable of dealing with arbitrarily-complex shapes is based on computing $\mathbf{x}$ by approximating the inverse of a map $\mathbf{h}$ which is comprised of a pair of harmonic functions in the target domain $\Omega$, i.e.,
\begin{align}
\label{eq:Adaptive_inverse_mapping_approach}
    \Delta \mathbf{h} = \mathbf{0} \quad \text{in } \Omega, \quad \text{s.t.} \quad \mathbf{h} \vert_{\partial \Omega} = \partial \hat{\Omega}.
\end{align}
The authors propose a two-step approach: First a large number of tuples $\left(\mathbf{h}(\mathbf{x}_j), \mathbf{x}_j \right)$, with $\mathbf{x}_j \in \Omega$ is computed using an isogeometric boundary element method \cite{rjasanow2007fast, sauter2010boundary}, after which the pairs are utilized to approximate $\mathbf{h}^{-1}$ through a least-squares minimization problem with regularization terms. \\
Seeking the mapping as the solution of an inverse-Laplace problem is equivalent to minimizing the \textit{Winslow} functional (\ref{eq:Adaptive_Winslow}), which follows straightforwardly from deriving the Euler-Lagrange equations of the minimization problem (\ref{eq:Adaptive_inverse_laplace}). As $\mathbf{h}$ is a pair of harmonic functions with convex target domain, it follows from the Rad\'o-Kneser-Choquet theorem that $\mathbf{h}$ is a diffeomorphism in the interior of $\Omega$ \cite{gravesen2012planar}, justifying an approximation of its inverse for the purpose of computing a domain parameterization. \\
A major advantage of the two-step approach from \cite{falini2015planar} over a direct minimization of (\ref{eq:Adaptive_Winslow}) is that the latter requires a folding-free initial domain parameterization to avoid division by zero. In the vast majority of cases, however, such a bijection is not available, limiting the method's scope to improving the parametric properties of an already bijective mapping. \\
An advantage of minimizing (\ref{eq:Adaptive_Winslow}), however, is that if the global minimum over $\mathcal{V}_h^2$ has been found, it is clearly bijective, while bijectivity may be lost in the indirect approach, due to numerical inaccuracies. \\
The observation that the impractical minimization of the \textit{Winslow} functional (\ref{eq:Adaptive_Winslow}) is equivalent to solving an inverse Laplace problem has lead to the development of (3.) PDE-based parameterization methods. To acquire a PDE-problem posed over $\hat{\Omega}$, we perform a pullback:
\begin{align}
\label{eq:Adaptive_inverse_mapping_approach_PDE}
    \Delta_\mathbf{x} \boldsymbol{\xi} = \mathbf{0} \quad \text{in } \hat{\Omega}, \quad \text{s.t.} \quad  \mathbf{x} \vert_{\partial \hat{\Omega}} & = \partial \Omega,
\end{align}
where $\Delta_\mathbf{x}$ denotes the Laplace-Beltrami \cite{kreyszig1991differential} operator with respect to $\mathbf{x}$. Problem (\ref{eq:Adaptive_inverse_mapping_approach_PDE}) suffers from the same shortcoming of the \textit{Winslow}-approach: the appearance of a Jacobian determinant in the denominator. However, we may scale the equation by multiplying with any nonsingular $2 \times 2$ tensor $T$. Choosing $T = ( \det \mathbf{x}_{\boldsymbol{\xi}} )^2 \mathbf{x}_{\boldsymbol{\xi}}$ (which is nonsingular in a neighbourhood of $\mathbf{x}$, thanks to the theoretically predicted bijectivity of the PDE-solution), the Jacobian determinant can be removed from (\ref{eq:Adaptive_inverse_mapping_approach_PDE}), leading to the following quasi-linear second-order PDE problem \cite{hinz2018spline}:
\begin{align}
\label{eq:Adaptive_inverse_mapping_approach_PDE_scaled}
    A(\mathbf{x}) \colon H(\mathbf{x}_i) & = 0 \quad \text{in } \hat{\Omega}, \quad \text{for } i \in \{1, 2\} \quad \text{s.t.} \quad \mathbf{x} \vert_{\partial \hat{\Omega}} = \partial \Omega,
\end{align}
where
\begin{align}
\label{eq:Adaptive_A}
H(u)_{ij} = \frac{\partial^2 u}{\partial \boldsymbol{\xi}_i \partial \boldsymbol{\xi}_j} \quad \text{and} \quad A(\mathbf{x}) = \frac{1}{g_{11} + g_{22} + \epsilon} \begin{pmatrix} g_{22} & - g_{12} \\
                                                                         - g_{12} & g_{11} \end{pmatrix},
\end{align}
with the $g_{ij}$ as in (\ref{eq:Adaptive_metric_tensor}) and $\epsilon$ a small positive constant that serves numerical stability (typically, $\epsilon \simeq 10^{-4}$). Furthermore, $A \colon B$ denotes the Frobenius inner product.
\begin{remark}
The purpose of dividing by $g_{11} + g_{22} + \epsilon$ in (\ref{eq:Adaptive_A}) is achieving scaling invariance.
\end{remark}
\noindent In \cite{hinz2018elliptic}, equation (\ref{eq:Adaptive_inverse_mapping_approach_PDE_scaled}) is discretized with a Galerkin approach. The resulting equations are tackled with a Newton-based iterative approach, which is initialized with an algebraic initial guess. \\
The advantages and disadvantages of solving (\ref{eq:Adaptive_inverse_mapping_approach_PDE_scaled}) over a direct minimization of (\ref{eq:Adaptive_Winslow}) are the same as in the indirect approach from \cite{falini2015planar}. Hence, folding resulting from insufficient numerical accuracy can be resolved by recomputing the mapping from a refined spline space. \\
In this paper, we will present several schemes for approximately solving (\ref{eq:Adaptive_inverse_mapping_approach_PDE_scaled}) with THB-spline bases. A major challenge in a THB-setting is deciding where a high resolution is needed. Since the approach is PDE-based, we adopt the a posteriori refinement strategy of \textit{dual weighted residual}, which is the topic of Section \ref{subsect:Adaptive_DWR}.
\section{Solution Strategies}
\label{sect:Adaptive_Solution_Strategies}
In this section we present several solution strategies to approximately solve (\ref{eq:Adaptive_inverse_mapping_approach_PDE_scaled}). As the resulting equations are nonlinear, we base our solution strategy on iterative approaches. Initial guesses are always constructed using an algebraic method (see Section \ref{sect:Adaptive_Introduction}). \\
Let 
\begin{align}
\label{eq:THB_basis_main}
\mathcal{U}^{\mathbf{f}} = \left \{ \mathbf{v} \in \mathcal{V}^2 \enskip \vert \enskip \mathbf{v} = \mathbf{f} \text{ on } \partial \hat{\Omega} \right \}
\end{align}
and let $\mathcal{U}^{\mathbf{f}}_h$ be the set resulting from replacing $\mathcal{V}$ by the finite-dimensional $\mathcal{V}_h \subset \mathcal{V}$ in (\ref{eq:THB_basis_main}). We have
\begin{align}
\label{eq:THB_basis_main_discrete}
\mathcal{U}^{\mathbf{f}}_h = \left \{ \mathbf{v} \in \mathcal{V}_h^2 \enskip \vert \enskip \mathbf{v} = \mathbf{f} \text{ on } \partial \hat{\Omega} \right \}.
\end{align}
\begin{remark}
For $\mathcal{U}^{\mathbf{f}}_h$ in (\ref{eq:THB_basis_main_discrete}) to be nonempty, we have to assume that $\mathbf{f}$ restricted to $\partial \hat{\Omega}$ is contained in $\mathcal{V}_h^2$, which may necessitate a projection of the Dirichlet data onto the finite-dimensional (THB-)spline space $\mathcal{V}_h^2$.
\end{remark}
\noindent Let $\mathbf{x}_D$ be such that $\mathbf{x}_D \vert_{\partial \hat{\Omega}}$ parameterizes $\partial \Omega$. For convenience we assume that $\mathbf{x}_D \in \mathcal{V}_h^2 \setminus \mathcal{U}^{\mathbf{0}}_h$. In an IGA-setting, (\ref{eq:Adaptive_inverse_mapping_approach_PDE_scaled}) suggests a discretization of the form:
\begin{align}
\label{eq:Adaptive_inverse_mapping_approach_PDE_scaled_discretized}
    & \text{find } \mathbf{x}_h \in \mathcal{U}^{\mathbf{x}_D}_h \quad \text{s.t.} \quad F(\mathbf{x}_h, \boldsymbol{\sigma}_h)= 0 \quad \forall \boldsymbol{\sigma}_h \in \mathcal{U}^{\mathbf{0}}_h,
\end{align}
with
\begin{align}
\label{eq:Adaptive_F_operator}
    F(\mathbf{x}, \boldsymbol{\sigma}) = \sum_{i = 1}^2 \int \limits_{\hat{\Omega}} \boldsymbol{\tau}_i(\boldsymbol{\sigma}, \mathbf{x}) A(\mathbf{x}) \colon H(\mathbf{x}_i) \mathrm{d} S,
\end{align}
for some $\boldsymbol{\tau}: \mathcal{U}^{\mathbf{0}} \times \mathcal{V}^2 \rightarrow L_2(\hat{\Omega}, \mathbb{R}^2)$. Unless stated otherwise, in the following, we assume $\boldsymbol{\tau}(\boldsymbol{\sigma}, \mathbf{x}) = \boldsymbol{\sigma}$. \\
As second order derivatives appear in (\ref{eq:Adaptive_F_operator}), in (\ref{eq:THB_basis_main}) we take $\mathcal{V}_h \subset H^2(\hat{\Omega})$.
\subsection{Newton Approach}
\label{subsect:Adaptive_Newton}
In the following, we briefly recapitulate the approach from \cite{hinz2018elliptic}, which is designated for tensor-product NURBS bases but can also be applied in a THB-setting. By
\begin{align}
    B^\prime(u, \ldots, z) \equiv \left. \frac{\partial B(u + \epsilon z, \ldots)}{\partial \epsilon} \right \vert_{\epsilon = 0},
\end{align}
we denote the Gateaux derivative of any differentiable form $B(\cdot, \ldots)$ with respect to its first argument. Given $\mathbf{x}^k \in \mathcal{U}^{\mathbf{x}_D}_h$, we compute the Newton increment from
\begin{align}
\label{eq:Adaptive_main_Newton}
    \text{find } \delta \mathbf{x}^k \in \mathcal{U}_h^{\mathbf{0}} \quad \text{s.t.} \quad F^\prime(\mathbf{x}^k, \boldsymbol{\sigma}_h, \delta \mathbf{x}^k) = - F(\mathbf{x}^k, \boldsymbol{\sigma}_h), \quad \forall \boldsymbol{\sigma}_h \in \mathcal{U}_h^{\mathbf{0}}.
\end{align}
Upon completion, we update $\mathbf{x}^{k + 1} = \mathbf{x}^k + \kappa \delta \mathbf{x}^k$ for some $\kappa \in (0, 1]$, whose optimal value is estimated using a line search routine. Above steps are repeated until the residual norm is deemed sufficiently small. \\
Optionally, derivative evaluations of the form $F^\prime(\mathbf{x}^k, \boldsymbol{\sigma}_h, \mathbf{v})$ may be approximated using finite differences:
\begin{align}
\label{eq:Adaptive_F_Gateaux_approximation}
    F^\prime(\mathbf{x}^k, \boldsymbol{\sigma}_h, \mathbf{v}) \simeq \frac{F(\mathbf{x}^k + \epsilon \mathbf{v}, \boldsymbol{\sigma}_h) - F(\mathbf{x}^k, \boldsymbol{\sigma}_h)}{\epsilon},
\end{align}
for $\epsilon$ small. Solving (\ref{eq:Adaptive_main_Newton}) using a suitable Krylov-subspace method only requires computing derivative evaluations $F^\prime(\mathbf{x}^k, \boldsymbol{\sigma}_h, \mathbf{v})$, which may be approximated using (\ref{eq:Adaptive_F_Gateaux_approximation}), leading to a Newton-Krylov algorithm that avoids the expensive assembly of the Jacobian matrix in (\ref{eq:Adaptive_main_Newton}). The optimal choice of $\epsilon$ in (\ref{eq:Adaptive_F_Gateaux_approximation}) is discussed in \cite{knoll2004jacobian}.

\subsection{Pseudo-Transient Continuation}
\label{subsect:Adaptive_Pseudotime}
In this technique, we seek the steady-state solution of the problem
\begin{align}
    \text{find } \mathbf{x}_h(\boldsymbol{\xi}, t) \in \mathcal{U}^{\mathbf{x}_D}_h, \quad \text{s.t.} \quad \left \langle \partial_t \mathbf{x}_h, \boldsymbol{\sigma}_h \right \rangle = -F(\mathbf{x}_h, \boldsymbol{\sigma}_h), \quad \forall \boldsymbol{\sigma}_h \in \mathcal{U}_h^{\mathbf{0}},
\end{align}
with
\begin{align}
    \left \langle \partial_t \mathbf{x}_h, \boldsymbol{\sigma}_h \right \rangle = \int \limits_{\hat{\Omega}} \boldsymbol{\sigma}_h \cdot \partial_t \mathbf{x}_h \mathrm{d} S.
\end{align}
Here, we only consider the choice $\boldsymbol{\tau}(\boldsymbol{\sigma}, \mathbf{x}) = \boldsymbol{\sigma}$. We discretize in time using backward Euler. Introducing $\delta \mathbf{x}^k = \mathbf{x}^{k+1} - \mathbf{x}^k$, with $F(\mathbf{x}^{k+1}, \boldsymbol{\sigma}_h) \simeq F(\mathbf{x}^{k}, \boldsymbol{\sigma}_h) + F^\prime(\mathbf{x}^k, \boldsymbol{\sigma}_h, \delta \mathbf{x}^k)$, we compute the temporal increment from
\begin{align}
    \text{find } \delta \mathbf{x}^k \in \mathcal{U}_h^{\mathbf{0}}, \quad \text{s.t.} \quad \left \langle \tfrac{\delta \mathbf{x}^k}{\delta t^k}, \boldsymbol{\sigma}_h \right \rangle + F^\prime(\mathbf{x}^k, \boldsymbol{\sigma}_h, \delta \mathbf{x}^k) = - F(\mathbf{x}^k, \boldsymbol{\sigma}_h), \quad \forall \boldsymbol{\sigma}_h \in \mathcal{U}_h^{\mathbf{0}},
\end{align}
where $\delta t^k$ denotes the time-step during the $k$-th iteration. As proposed in \cite{kelley1998convergence}, we base the time-step selection on the following recursive formula
\begin{align}
    \delta t^k = \delta t^{k - 1} \frac{ \left \| \mathbf{F}(\mathbf{x}^{k-1}) \right \|_2 }{ \left \| \mathbf{F}(\mathbf{x}^{k}) \right \|_2}, \quad \text{with} \quad \left \| \mathbf{F}(\mathbf{x}) \right \|^2_2 = \sum \limits_{\boldsymbol{\sigma}_h \in \left[ \mathcal{U}_h^{\mathbf{0}} \right]} F(\mathbf{x}, \boldsymbol{\sigma}_h)^2.
\end{align}
The iteration is terminated once $\| \mathbf{x}^{k} - \mathbf{x}^{k-1} \|$ is sufficiently small (in a suitable norm).

\subsection{Picard Iteration}
\label{subsect:Adaptive_Picard}
In the following, we present a Picard-based iterative scheme that is loosely based on the default approach from the rich literature of classical meshing techniques \cite{thompson1998handbook}. As opposed to Sections \ref{subsect:Adaptive_Newton} and \ref{subsect:Adaptive_Pseudotime}, we base the scheme on a linearize-then-discretize approach, rather than the converse. Note that for given $\mathbf{x} = (x, y)^T$, we have
\begin{align}
    A(\mathbf{x}) = C^T(\mathbf{x}) C(\mathbf{x}), \quad \text{with} \quad C(\mathbf{x}) = \frac{1}{\sqrt{g_{11} + g_{22} +\epsilon}} \begin{pmatrix} \frac{\partial y}{\partial \eta} & - \frac{\partial y}{\partial \xi} \\
    - \frac{\partial x}{\partial \eta} & \frac{\partial x}{\partial \xi} \end{pmatrix}.
\end{align}
As such, $A(\mathbf{x})$ is symmetric positive semi-definite (SPSD) for all $\mathbf{x}$ and symmetric positive definite (SPD) for $\mathbf{x}: \hat{\Omega} \rightarrow \Omega$ bijective. Let us introduce the operator $\mathbf{K}: C^2(\hat{\Omega}, \mathbb{R}^2) \times C^2(\hat{\Omega}, \mathbb{R}^2) \times \mathbb{R}^+ \rightarrow C^0(\hat{\Omega}, \mathbb{R}^2)$ with components
\begin{align}
\label{eq:THB_Picard_stabilized_operator}
    \mathbf{K}_i(\mathbf{x}, \mathbf{y}, \mu) = A_\mu(\mathbf{y}) \colon H(\mathbf{x}_i) - \mu \Delta_{\boldsymbol{\xi}} \mathbf{y}_i, \quad \text{where} \quad A_\mu(\mathbf{y}) = A(\mathbf{y}) + \mu I^{2 \times 2}.
\end{align}
Note that for $\mu > 0$, $A_\mu(\mathbf{x})$ is SPD and that for all choices of $\mu$, $\mathbf{K}_i(\mathbf{x}, \mathbf{x}, \mu) = A(\mathbf{x}) \colon H(\mathbf{x}_i)$. For given $\mu > 0$, we seek $\mathbf{x}$ as the limit $k \rightarrow \infty$ of the recursive sequence
\begin{align}
\label{eq:Adaptive_Picard_not_discretized}
    \text{find } \mathbf{x}^{k+1} \quad \text{s.t.} \quad \mathbf{K}(\mathbf{x}^{k+1}, \mathbf{x}^k, \mu) = \mathbf{0}, \quad \text{and} \quad \mathbf{x}^{k+1} = \mathbf{x}_D \quad \text{on} \quad \partial \hat{\Omega}.
\end{align}
To discretize (\ref{eq:Adaptive_Picard_not_discretized}), let us introduce the semi-linear form $G_\tau: \mathcal{V}^2 \times \mathcal{V}^2 \times \mathbb{R}^+ \times \mathcal{U}^{\mathbf{0}} \rightarrow \mathbb{R}$ with
\begin{align}
\label{eq:Adaptive_Picard_G_tau}
    G_\tau(\mathbf{x}, \mathbf{y}, \mu, \boldsymbol{\sigma}) & = \sum_{i=1}^2 \int \limits_{\hat{\Omega}} \boldsymbol{\tau}_i(\boldsymbol{\sigma}, \mathbf{y}) \left( A_\mu (\mathbf{y}) \colon H(\mathbf{x}_i)- \mu \Delta_{\boldsymbol{\xi}} \mathbf{y}_i \right) \mathrm{d} S.
\end{align}
Given $\mathbf{x}^k$, we compute $\mathbf{x}^{k+1} \in \mathcal{U}^{\mathbf{x}_D}$ as the solution of
\begin{align}
\label{eq:Adaptive_main_Picard}
    \text{find } \mathbf{x}^{k+1} \in \mathcal{U}^{\mathbf{x}_D} \quad \text{s.t } \quad G_\tau(\mathbf{x}^{k+1}, \mathbf{x}^k, \mu, \boldsymbol{\sigma}_h) = 0, \quad \forall \boldsymbol{\sigma}_h \in \mathcal{U}^{\mathbf{0}},
\end{align}
where, as before, $\mathcal{U}^{\mathbf{x}_D} = \{ \mathbf{v} \in \mathcal{V}^2 \enskip \vert \enskip \mathbf{v} = \mathbf{x}_D \text{ on } \partial \hat{\Omega} \}$. \\
The discretization of (\ref{eq:Adaptive_main_Picard}) follows straightforwardly from replacing $\mathcal{V}$ by the finite dimensional $\mathcal{V}_h \subset \mathcal{V}$. Equation (\ref{eq:Adaptive_main_Picard}) leads to a decoupled (block-diagonal) system of elliptic equations in nonvariational (or non-divergence) form \cite{lakkis2011finite}. Inspired by \cite{gallistl2017variational}, here we consider the choices
\begin{align}
\label{eq:THB_tau_choices}
    \boldsymbol{\tau}^\text{Id}(\boldsymbol{\sigma}, \mathbf{y}) = \boldsymbol{\sigma}, \quad \boldsymbol{\tau}^\text{div}(\boldsymbol{\sigma}, \mathbf{y}) = \gamma(\mathbf{y}) \Delta_{\boldsymbol{\xi}} \boldsymbol{\sigma} \quad \text{and} \quad \boldsymbol{\tau}^\text{ls}_i(\boldsymbol{\sigma}, \mathbf{y}) = A_\mu(\mathbf{y}) \colon H(\boldsymbol{\sigma}_i),
\end{align}
where
\begin{align}
    \gamma(\mathbf{y}) = \frac{\operatorname{trace}(A_\mu(\mathbf{y}))}{A_\mu \colon A_\mu (\mathbf{y})}.
\end{align}
A Picard scheme results from iterating on (\ref{eq:Adaptive_main_Picard}) until $ \| \mathbf{x}^{k+1} - \mathbf{x}^{k} \|$ is negligibly small.
\begin{remark}
Adding artificial diffusion in (\ref{eq:THB_Picard_stabilized_operator}) stabilizes the linearized discrete equation from (\ref{eq:Adaptive_main_Picard}). In the absence of stabilization (i.e., $\mu = 0$), (\ref{eq:Adaptive_main_Picard}) can be ill-posed in rare cases, depending on the previous iterate $\mathbf{x}^k$. This is also true for a Newton-based approach. Whenever an invalid iterate is encountered in the Newton approach, we fall back on the techniques from this section. \\
For $\mu > 0$, well-posedness of (\ref{eq:Adaptive_main_Picard}) with the choices from (\ref{eq:THB_tau_choices}) is discussed in \cite{gallistl2017variational} and \cite{blechschmidt2019error}. Stabilizing a Newton-based approach constitutes a topic for future research.
\end{remark}

\begin{figure}[h!]
  \centering
  \subfloat[The unrefined domain.]{\includegraphics[width=.45\linewidth]{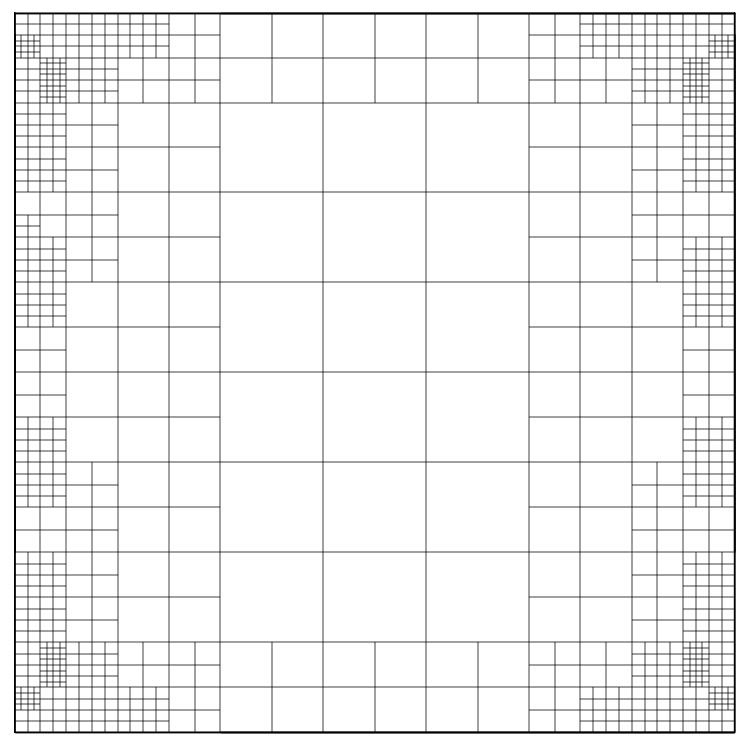}\label{fig:solution_strategies_domain}}
  \qquad
  \subfloat[The uniformly refined domain.]{\includegraphics[width=.45\linewidth]{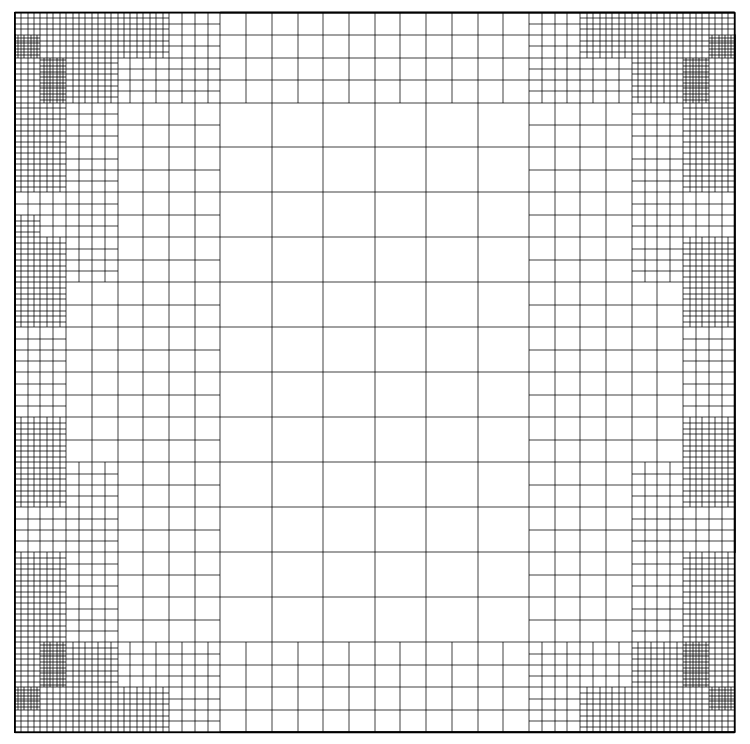}\label{fig:solution_strategies_domain_ref}}
  \caption{The THB-refined parametric domains used in the computations of the parameterizations from Figures \ref{fig:solution_strategies} and \ref{fig:solution_strategies_ref}.}
  \label{fig:solution_strategies_domains}
\end{figure}

\begin{figure}[h!]
  \centering
  \subfloat[Reference mapping acquired from minimizing the Winslow functional over the domain from Figure \ref{fig:solution_strategies_domain}.]{\includegraphics[width=.45\linewidth]{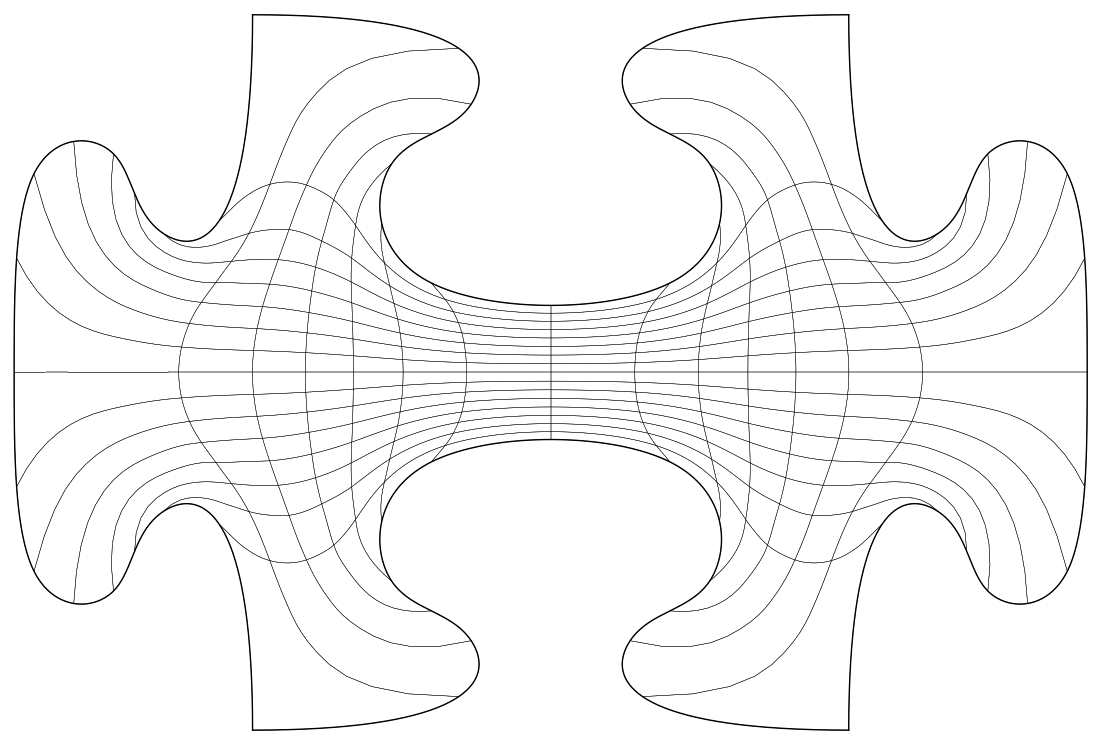}\label{fig:solution_strategies_wins}}
  \centering $\quad$
  \subfloat[The parameterization for the choice $\boldsymbol{\tau} = \boldsymbol{\tau}^\text{Id}$.]{\includegraphics[width=.45\linewidth]{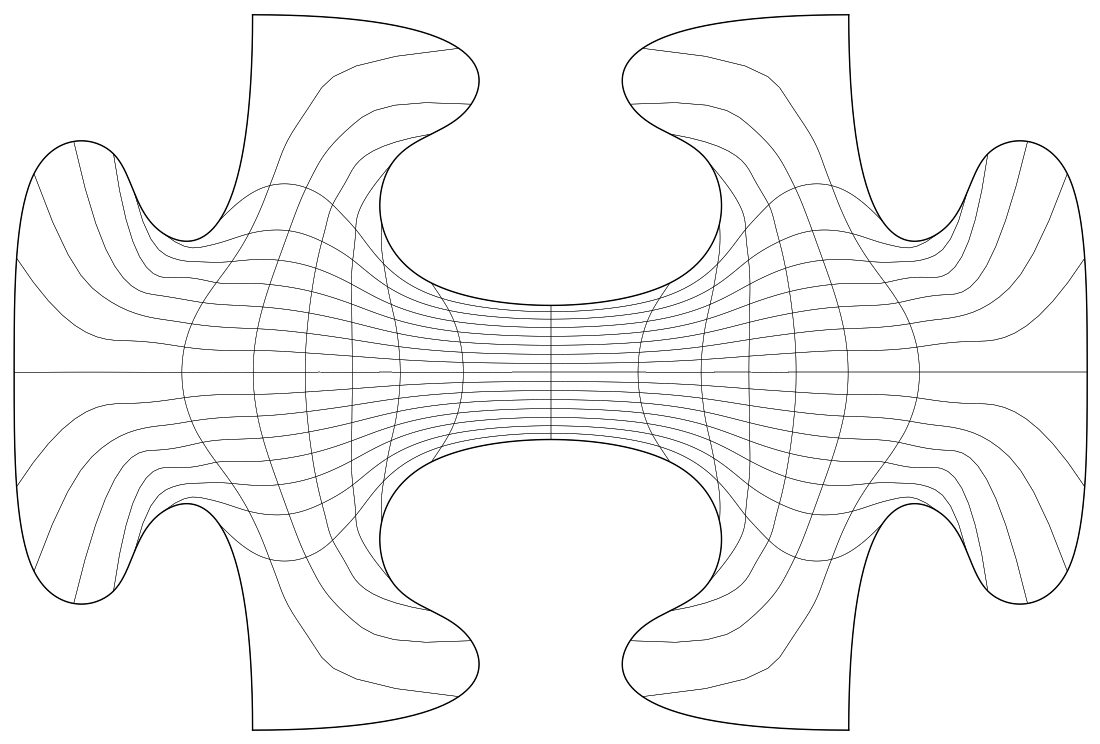}\label{fig:solution_strategies_id}} \\
  \centering
  \subfloat[The parameterization for the choice $\boldsymbol{\tau} = \boldsymbol{\tau}^\text{ls}$.]{\includegraphics[width=.45\linewidth]{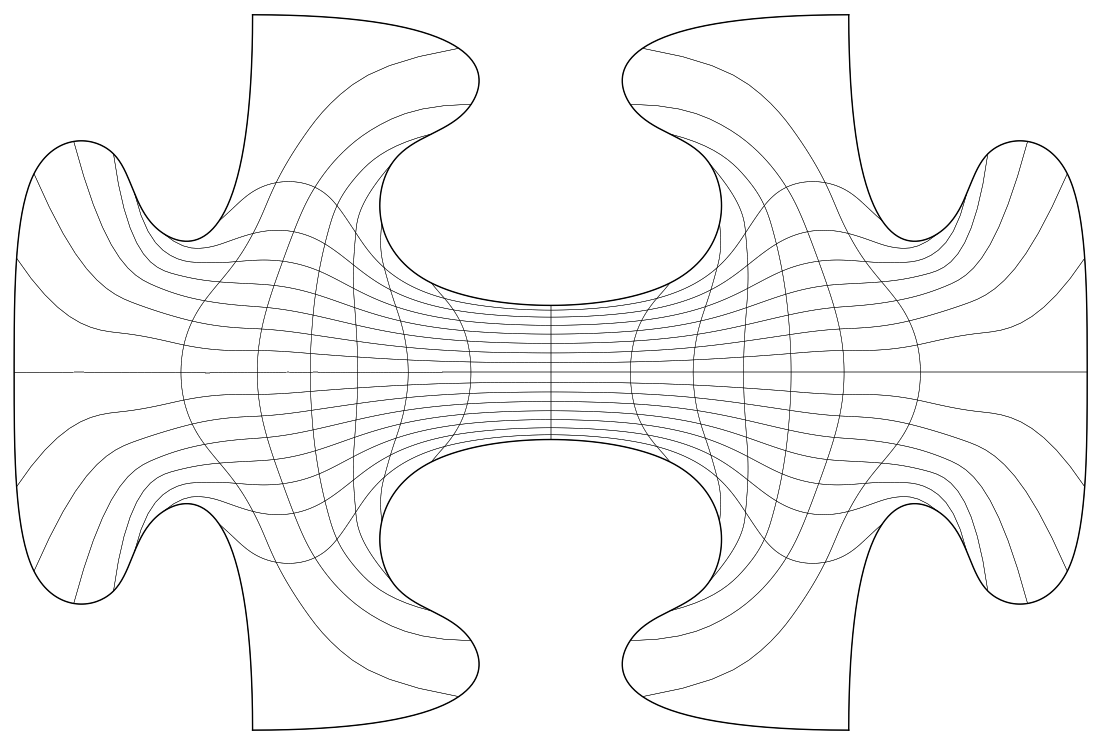}\label{fig:solution_strategies_ls}}
  \centering $\quad$
  \subfloat[The parameterization for the choice $\boldsymbol{\tau} = \boldsymbol{\tau}^\text{div}$.]{\includegraphics[width=.45\linewidth]{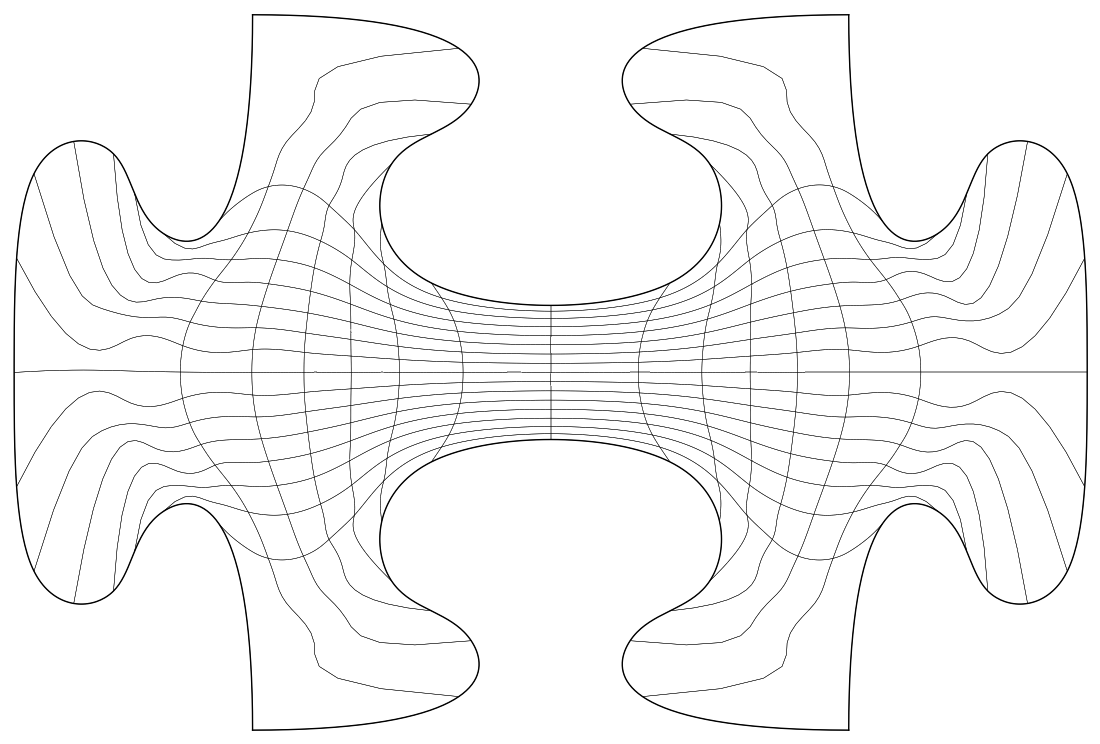}\label{fig:solution_strategies_div}}
  \centering
  \caption{Parameterizations acquired using the various discretization techniques.}
  \label{fig:solution_strategies}
\end{figure}

\begin{figure}[h!]
  \centering
  \subfloat[Winslow.]{\includegraphics[width=.45\linewidth]{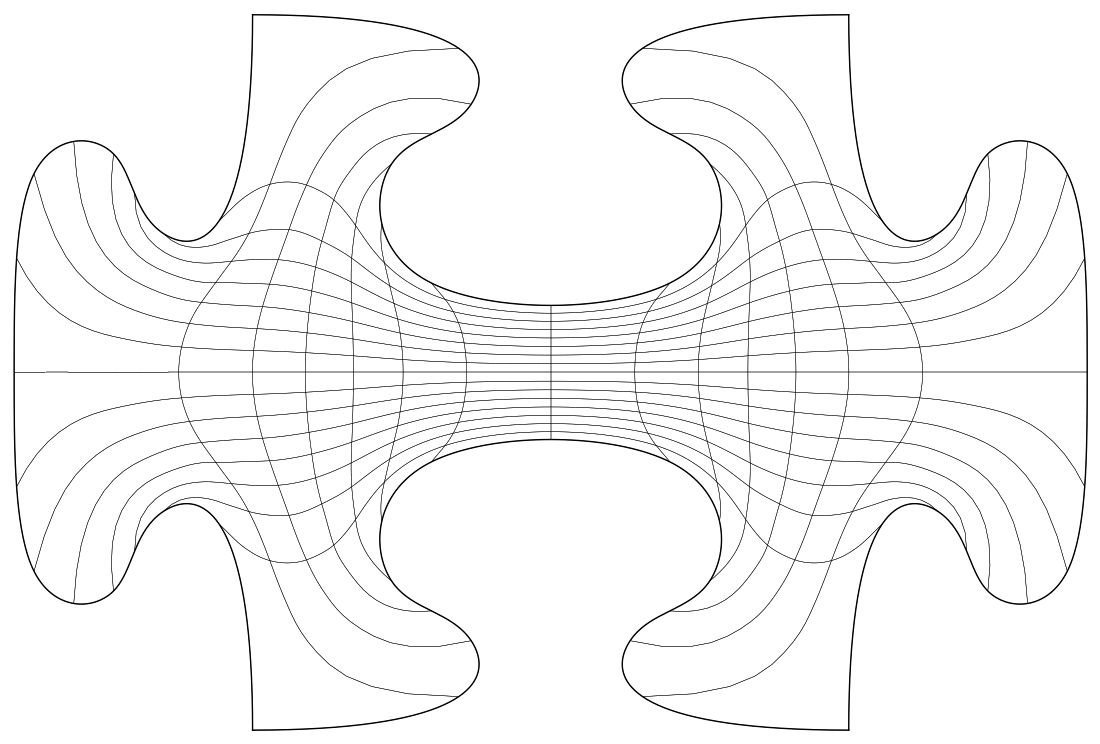}\label{fig:solution_strategies_wins_ref}}
  \centering $\quad$
  \subfloat[The refined parameterization for $\boldsymbol{\tau} = \boldsymbol{\tau}^\text{Id}$.]{\includegraphics[width=.45\linewidth]{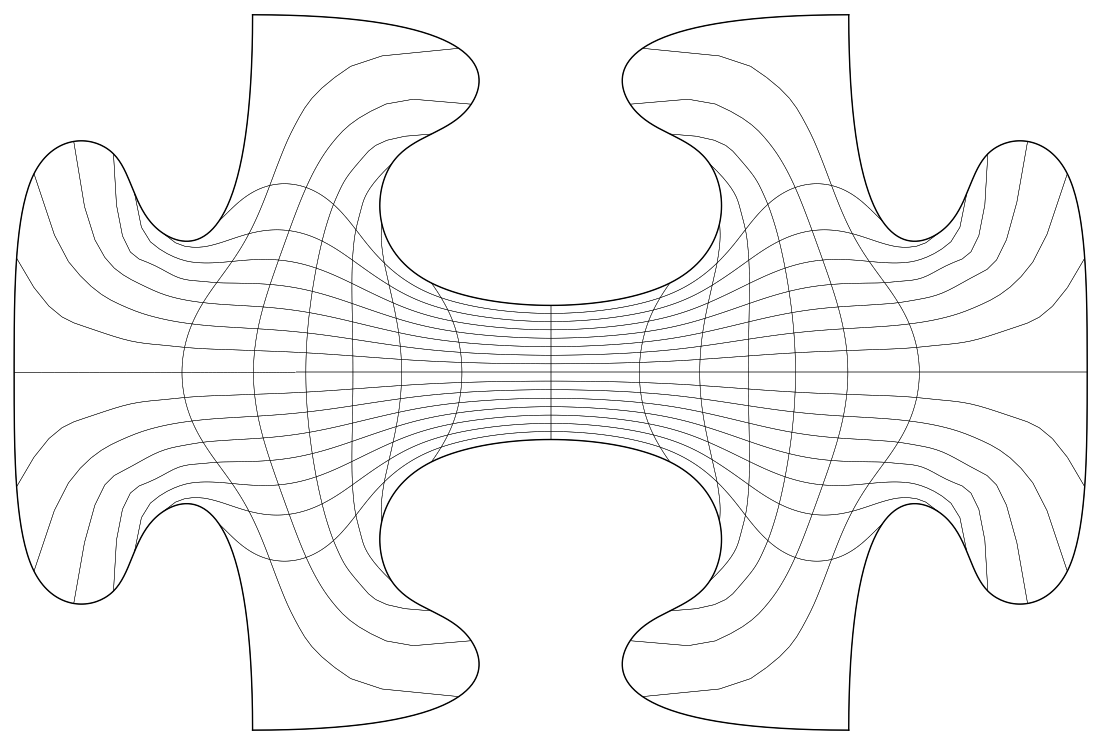}\label{fig:solution_strategies_id_ref}} \\
  \centering
  \subfloat[The refined parameterization for $\boldsymbol{\tau} = \boldsymbol{\tau}^\text{ls}$.]{\includegraphics[width=.45\linewidth]{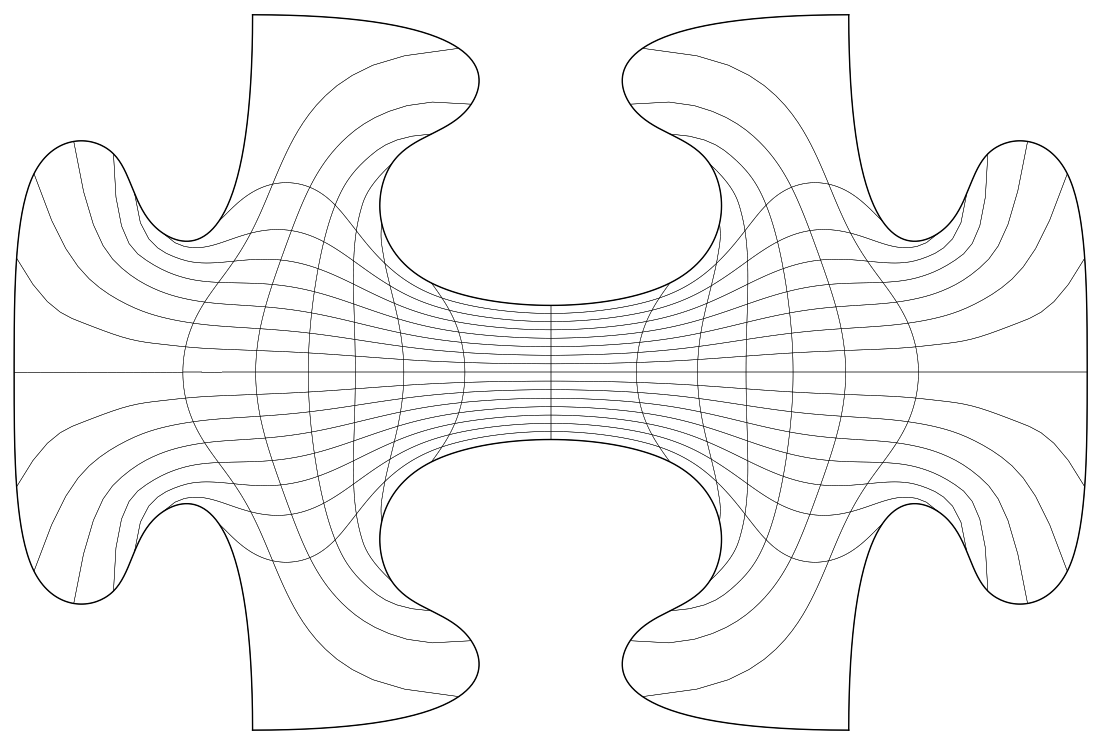}\label{fig:solution_strategies_ls_ref}}
  \centering $\quad$
  \subfloat[The refined parameterization for $\boldsymbol{\tau} = \boldsymbol{\tau}^\text{div}$.]{\includegraphics[width=.45\linewidth]{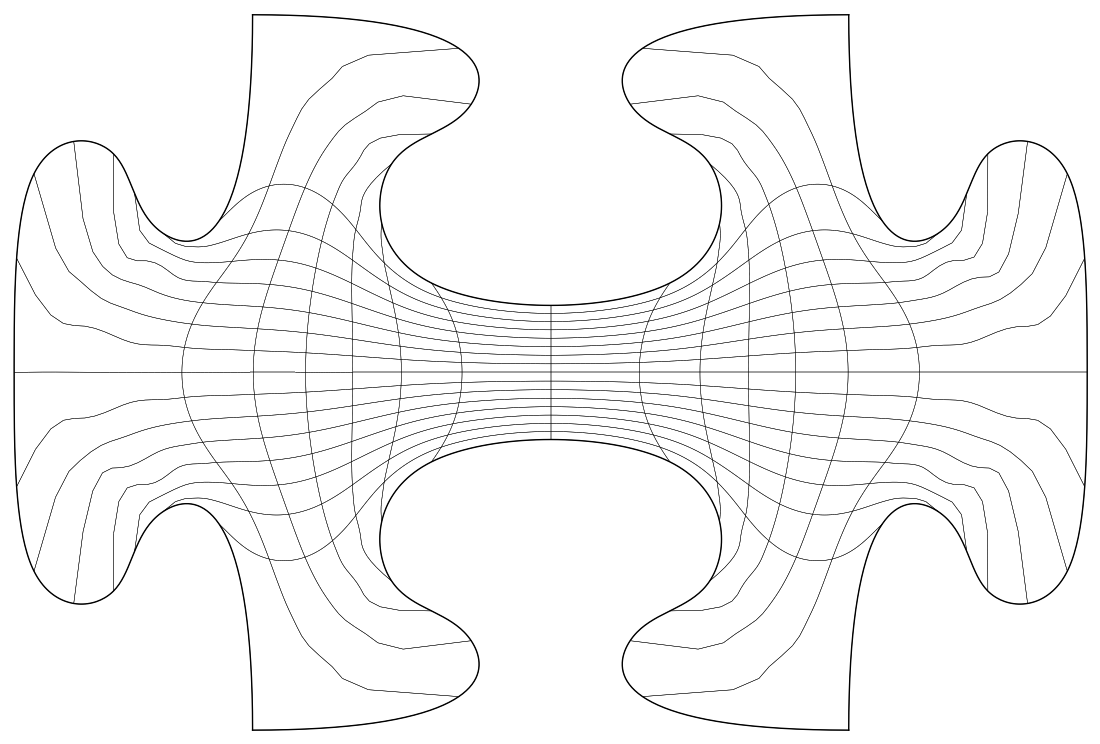}\label{fig:solution_strategies_div_ref}}
  \centering
  \caption{Parameterizations acquired using the various discretization techniques over the uniformly refined domain.}
  \label{fig:solution_strategies_ref}
\end{figure}

\begin{table}[h!]
\centering
\begin{tabular}{l|c|c|c|c}
  \backslashbox{\small refinement}{\small method}  & Direct & $\boldsymbol{\tau} = \boldsymbol{\tau}^\text{Id}$ & $\boldsymbol{\tau} = \boldsymbol{\tau}^\text{ls}$ & $\boldsymbol{\tau} = \boldsymbol{\tau}^\text{div}$ \\ \hline
$h$   & $4.784$ & $4.849$ & $4.913$ & $4.974$ \\ 
$h/2$ &         & $4.787$ & $4.790$ & $4.815$ \\                                    
\end{tabular}
\caption{Evaluation of the Winslow functional with the various parameterizations.}
\label{tab:solution_strategies_Winslow}
\end{table}

\subsection{Direct Approach}
\label{subsect:Adaptive_Direct_Minimization}
Assuming a bijective initial guess $\mathbf{x}^0 \in \mathcal{U}^{\mathbf{x}_D}_h$ is available, we may alternatively compute an approximately inversely harmonic parameterization by a direct minimization of the Winslow functional (\ref{eq:Adaptive_Winslow}). Let
\begin{align}
\label{eq:Adaptive_Winslow_costfunc}
    L_W(\mathbf{x}) = \int_{\hat{\Omega}} Q_\text{W}(\mathbf{x}) \mathrm{d} S
\end{align}
denote the evaluation of the Winslow function (see equation (\ref{eq:Adaptive_Winslow})), whose domain is the set of all bijective $\mathbf{x}$. To conform with the topic of this paper, we compute the minimizer over $\mathcal{U}^{\mathbf{x}_D}_h$ as the solution of the following discretized PDE problem:
\begin{align}
\label{eq:Adaptive_Winslow_PDE_discretized}
    \text{find } \mathbf{x}_h \in \mathcal{U}^{\mathbf{x}_D}_h \quad \text{s.t.} \quad L_W^\prime(\mathbf{x}_h, \boldsymbol{\sigma}_h) = 0, \quad \forall \boldsymbol{\sigma}_h \in \mathcal{U}_h^{\mathbf{0}}.
\end{align}
We solve (\ref{eq:Adaptive_Winslow_PDE_discretized}) with one of the approaches from Sections \ref{subsect:Adaptive_Newton} and \ref{subsect:Adaptive_Pseudotime}. Typically $\mathbf{x}^0$ is the solution of one of the indirect methods presented in Sections \ref{subsect:Adaptive_Newton} to \ref{subsect:Adaptive_Picard}. In practice, we have often encountered convergence failure even when $\mathbf{x}^0$ is bijective. As a rule of thumb, we retry solving (\ref{eq:Adaptive_Winslow_PDE_discretized}) with a refined $\mathbf{x}^0$, resulting from an indirect approach, if converge is not reached after a few iterations.
\begin{remark}
If a measure of quality of the solution results from substituting into (\ref{eq:Adaptive_Winslow_costfunc}), a direct approach yields the best outcome.
\end{remark}

\subsection{Example: Puzzle Piece}
Figure \ref{fig:solution_strategies} shows the various parameterizations of a puzzle piece geometry, resulting from solving the discretized equations with the Newton-approach from Section \ref{subsect:Adaptive_Newton} and the different choices of $\boldsymbol{\tau}: \mathcal{U}^{\mathbf{0}} \times \mathcal{V}^2 \rightarrow L_2(\hat{\Omega}, \mathbb{R}^2)$ from (\ref{eq:THB_tau_choices}). For Newton, stabilization is avoided, i.e., $\mu = 0$. All methods lead to a bijective outcome. However, the figure shows noticeable differences in the parametric properties between the various methods, in particular in the protruded parts and in particular in Figure \ref{fig:solution_strategies_div}. Upon uniform refinement of the underlying element segmentation of $\hat{\Omega}$ (see Figure \ref{fig:solution_strategies_domains}) and recomputation over the associated refined THB-spline basis, the differences become less pronounced, suggesting that all schemes are consistent. The associated parameterizations are depicted in Figure \ref{fig:solution_strategies_ref}. Table \ref{tab:solution_strategies_Winslow} shows the outcomes of substituting the various parameterizations into (\ref{eq:Adaptive_Winslow_costfunc}). Not surprisingly, the choice $\boldsymbol{\tau} = \boldsymbol{\tau}^\text{div}$ fares the worst while the table suggests that $\boldsymbol{\tau} = \boldsymbol{\tau}^{Id}$ is the best choice. Upon refinement, the $\boldsymbol{\tau} = \boldsymbol{\tau}^{Id}$ and $\boldsymbol{\tau} = \boldsymbol{\tau}^{ls}$ parameterizations become virtually indistinguishable from the global minimizer over the coarse space, which is also reflected in table \ref{tab:solution_strategies_Winslow}. \\
As documented in the literature \cite{steinberg1993fundamentals}, all parameterizations suffer from the well-known pathologies of inversely harmonic maps, such as the tendency to yield large elements within protruded parts. Fortunately, in a THB-setting this can be compensated for by performing local refinement in the affected regions. Mitigating the impact of these pathologies will be the topic of Section \ref{sect:Adaptive_Domain_Optimization}.

\section{A Basic Scheme Based on a Posteriori Refinement}
\label{sect:Adaptive_Basic_Scheme}
One of the main challenges of PDE-based parameterization is selecting an appropriate finite-dimensional spline space $\mathcal{V}_h$. For this, we employ the technique of \textit{Dual Weighted Residual}, which will be the topic of Section \ref{subsect:Adaptive_DWR}.
\subsection{Dual Weighted Residual}
\label{subsect:Adaptive_DWR}
\textit{Dual Weighted Residual}, is an a posteriori error estimation technique that is based on duality considerations. Consider a semi-linear differential form $A(u, \phi)$ (which is linear in $\phi$). We consider the problem
\begin{align}
    \label{eq:Adaptive_DWR_main}
    \text{find } u \in \mathcal{V}^\circ \quad \text{s.t.} \quad A(u, \phi) = f(\phi), \quad \forall \phi \in \mathcal{V}^\circ,
\end{align}
for some linear functional $f(\cdot)$ and a suitably-chosen vector space $\mathcal{V}$, with $\mathcal{V}^\circ = \mathcal{V} \cap H^1_0(\hat{\Omega})$. We seek an approximate solution $u_h \in \mathcal{V}_h^\circ$ with $\mathcal{V}_h \subset \mathcal{V}$ by solving a discretized counterpart of (\ref{eq:Adaptive_DWR_main})
\begin{align}
\label{eq:Adaptive_DWR_main_discretized}
    \text{find } u_h \in \mathcal{V}_h^\circ \quad \text{s.t.} \quad A(u_h, \phi_h) = f(\phi_h), \quad \forall \phi_h \in \mathcal{V}_h^\circ.
\end{align}
Let $L(u)$ be such that
\begin{align}
\label{eq:Adaptive_DWR_DeltaL}
    \Delta L(u_h) \equiv L(u) - L(u_h)
\end{align}
is a quantity of interest (which for instance measures the global quality of the approximation). Furthermore, let
\begin{align}
    \rho(u, \psi) = f(\psi) - A(u, \psi)
\end{align}
denote the residual. \\
If $z$ is the solution of
\begin{align}
    \text{find } z \in \mathcal{V}^\circ \quad \text{s.t.} \quad A^\prime(u, \phi, z) = L^\prime(u, \phi), \quad \forall \phi \in \mathcal{V}^\circ,
\end{align}
we have
\begin{align}
    \Delta L(u_h) = \rho(u_h, z - \psi_h) + R_h(e),
\end{align}
for arbitrary $\psi_h \in \mathcal{V}_h^\circ$ and some $R_h$ that is quadratic in $e \equiv u - u_h$ \cite{rannacher2004adaptive}. In practice, we neglect $R_h$ and approximate $z$ by the solution of the discrete adjoint equation
\begin{align}
\label{eq:Adaptive_DWR_discrete_adjoint}
    \text{find } z_h \in \overline{\mathcal{V}}_h^\circ \quad \text{s.t.} \quad A^{\prime}(u_h, \sigma_h, z_h) = L^{\prime}(u_h, \sigma_h), \quad \forall \sigma_h \in \overline{\mathcal{V}}_h^\circ,
\end{align}
for some adjoint (THB-)spline space $\overline{\mathcal{V}}_h \subset \mathcal{V}$. Hence,
\begin{align}
\label{eq:Adaptive_adjoint_equation_discrete}
    \Delta L(u_h) \simeq \rho(u_h, z_h - \psi_h) = \sum_{w_i \in \left[\mathcal{V}_h \right]} \rho(u_h, w_i(z_h - \psi_h)) \equiv \sum_i \mathbf{r}_i(u_h),
\end{align}
thanks to semi-linearity of $A(\cdot, \cdot)$ and the partition of unity property associated with $\left[ \mathcal{V}_h \right]$. \\
The motivation to use an adjoint spline space that differs from $\mathcal{V}_h$ is the fact that substituting any $z_h \in \mathcal{V}_h^\circ$ in (\ref{eq:Adaptive_adjoint_equation_discrete}) results in $\Delta L(u_h) = 0$ due to Galerkin orthogonality, making it a meaningless error estimate. \\
The appeal of using (\ref{eq:Adaptive_adjoint_equation_discrete}) is that a scalar quantity of interest $\Delta L(u_h)$ is transformed into an integral quantity over $\hat{\Omega}$, which in turn is decomposed into the basis function wise contributions $\mathbf{r}_i(u_h)$. The vector $\mathbf{r}(u_h)$ may then be utilized in the selection of basis functions for goal-oriented refinement (see Section \ref{subsect:Adaptive_Refinement_Strategies}).
\begin{remark}
If $u_h$ is a very inaccurate approximation of $u$, the discrete adjoint solution $z_h$ will be inaccurate regardless of the choice of $\overline{\mathcal{V}}_h$. Heuristically, we have rarely encountered this situation in the examples considered in this work. In case refinement is ineffective, the procedure should be restarted with a uniformly refined initial basis.
\end{remark}

\subsection{Applications to PDE-Based Parameterization}
In this section, we apply the methodology from Section \ref{subsect:Adaptive_DWR} to the PDE-based parameterization problem (\ref{eq:Adaptive_inverse_mapping_approach_PDE_scaled_discretized}). Let $\mathbf{x}_D$ be the canonical extension of the Dirichlet data as introduced in (\ref{eq:Adaptive_inverse_mapping_approach_PDE_scaled_discretized}). With $\mathbf{x}_h = \mathbf{x}_D + \mathbf{x}_0$, we may write (\ref{eq:Adaptive_inverse_mapping_approach_PDE_scaled_discretized}) in the equivalent form
\begin{align}
\label{eq:Adaptive_EGG_problem_DWR_form}
    \text{find } \mathbf{x}_0 \in \mathcal{U}_h^{\mathbf{0}} \quad \text{s.t.} \quad F(\mathbf{x}_D + \mathbf{x}_0, \boldsymbol{\sigma}_h) = 0, \quad \forall \boldsymbol{\sigma}_h \in \mathcal{U}_h^{\mathbf{0}}.
\end{align}
In the formalism of (\ref{eq:Adaptive_DWR_main_discretized}), we hence have $A(\mathbf{x}, \boldsymbol{\sigma}) = F(\mathbf{x}_D + \mathbf{x}, \boldsymbol{\sigma})$ and $f(\boldsymbol{\sigma}) = 0$. Alternatively, we may absorb the dependence on $\mathbf{x}_D$ in $f(\cdot)$. As before, the relation between $\mathcal{V}_h$ and $\mathcal{U}_h^{\mathbf{f}}$ follows from (\ref{eq:THB_basis_main_discrete}). \\
We would like to design scalar cost functions ($L(u)$ in (\ref{eq:Adaptive_DWR_DeltaL})) to aid us in refining an a priori chosen basis $[\mathcal{V}_h]$ such that after recomputing the solution over the refined space $\mathcal{V}_h^R \supset \mathcal{V}_h$,
\begin{enumerate}
    \item $\mathbf{x}_h^R$ is bijective;
    \item $\mathbf{x}_h^R$ approximates $\mathbf{x}$ well.
\end{enumerate}
In a discrete setting, we may relax the condition that $\mathbf{x}_h^R$ be bijective by the condition that $\mathbf{x}_h^R$ has a positive Jacobian determinant in all quadrature points $\Xi = \{\boldsymbol{\xi}_1^q, \ldots, \boldsymbol{\xi}_M^q \}$. \\
As such, let $\mathbf{x}_h$ be the solution of (\ref{eq:Adaptive_EGG_problem_DWR_form}) over the space $\mathcal{V}_h$ and let 
\begin{align}
\label{eq:Adaptive_DWR_quadrature_subset}
    \Xi_{\mathbf{-}} = \left \{ \boldsymbol{\xi}_i^q \in \Xi \enskip \vert \enskip \det J(\mathbf{x}_{h}) < 0 \text{ in } \boldsymbol{\xi}_i^q \right \}.
\end{align}
To address (potential) lack of bijectivity, we propose the following goal-oriented cost function:
\begin{align}
\label{eq:Adaptive_bijectivity_costfunc}
    L_{\Xi}(\mathbf{x}) = \sum_{\boldsymbol{\xi}_i^q \in \Xi_{\mathbf{-}}} \det J(\mathbf{x})(\boldsymbol{\xi}_i^q),
\end{align}
such that
\begin{align}
    \Delta L_{\Xi}(\mathbf{x}_h) = \underbrace{L_\Xi(\mathbf{x})}_{\geq 0} - \underbrace{L_\Xi(\mathbf{x}_h)}_{\leq 0} \geq 0,
\end{align}
with equality if and only if $\Xi_{\mathbf{-}} = \emptyset$. Here, the inequality $L_\Xi(\mathbf{x}) \geq 0$ follows from the Rad\'o-Kneser-Choquet theorem (see Section \ref{sect:Adaptive_Introduction}) while $L_\Xi(\mathbf{x}_h) \leq 0$ follows from (\ref{eq:Adaptive_DWR_quadrature_subset}). According to (\ref{eq:Adaptive_adjoint_equation_discrete}), we may approximate
\begin{align}
\label{eq:Adaptive_DWR_basis_function_split}
    \Delta L_\Xi(\mathbf{x}_h) & \simeq -F(\mathbf{x}_h, \mathbf{z}_h - \boldsymbol{\psi}_h)
                                   = \sum_{w_i \in \left[ \mathcal{V}_h \right]} - F \left(\mathbf{x}_h, w_i( \mathbf{z}_h - \boldsymbol{\psi}_h ) \right) \equiv \sum_i \mathbf{r}_i(\mathbf{x}_h).
\end{align}
Typically, we choose $\boldsymbol{\psi}_h$ as the $L_2(\hat{\Omega}, \mathbb{R}^2)$-projection of $\mathbf{z}_h$ onto $\mathcal{U}_h^{\mathbf{0}}$.
\begin{remark}
Even though subtracting a nonzero $\boldsymbol{\psi}_h \in \mathcal{U}_h^{\mathbf{0}}$ does not alter the outcome on the right hand side of (\ref{eq:Adaptive_DWR_basis_function_split}), it does influence its decomposition into the basis function wise contributions $\mathbf{r}_i(\mathbf{x}_h)$. Here, the proposed choice of $\boldsymbol{\psi}_h$ serves to retain the sharpness of the error bound.
\end{remark}
\noindent Using the basis function wise decomposition of the quantity of interest $\Delta L_\Xi(\mathbf{x}_h)$, the procedure selects a subset of the $w_i \in \left[ \mathcal{V}_h \right]$ and marks them for refinement. We propose selection criteria in Section \ref{subsect:Adaptive_Refinement_Strategies}. \\
After refinement of $\mathcal{V}_h$, we recompute the mapping from the enriched basis $\mathcal{V}_h^R$ and if necessary repeat above steps until discrete bijectivity (over $\Xi$) has been achieved.
\begin{remark}
For better performance, we always use the prolonged coarse-grid solution as an initial guess for recomputing the mapping under the refined basis.
\end{remark}
\noindent Upon completion, we may choose to settle for the (possibly inaccurate but with respect to the $\boldsymbol{\xi}_i^q \in \Xi$ analysis-suitable) resulting mapping $\mathbf{x}_h^R$, or we may choose to further improve its accuracy with respect to the exact solution. As the exact solution of the PDE problem is equal to the minimizer of the Winslow function (see Section \ref{sect:Adaptive_related_work})
\begin{align*}
    L_W(\mathbf{x}) = \int_{\hat{\Omega}} \frac{g_{11} + g_{22}}{ \det J} \mathrm{d} S,
\end{align*}
by choosing $- L_W(\mathbf{x})$ as a cost function, we acquire the quantity of interest
\begin{align}
\label{eq:Adaptive_dWinslow_costfunc}
    \Delta L_W(\mathbf{x}_h) = -L_W(\mathbf{x}) + L_W(\mathbf{x}_h) \geq 0,
\end{align}
with equality for $\| \mathbf{x} - \mathbf{x}_h \|_{H^1(\hat{\Omega})} = 0$. As such, (\ref{eq:Adaptive_dWinslow_costfunc}) may serve as a measure for the distance of $\mathbf{x}_h$ to $\mathbf{x}$. As before, we approximate (\ref{eq:Adaptive_dWinslow_costfunc}) by substituting the discrete adjoint solution $\mathbf{z}_h$ in (\ref{eq:Adaptive_adjoint_equation_discrete}) and base refinement criteria on the basis function wise contributions to (\ref{eq:Adaptive_dWinslow_costfunc}). The steps of refinement, recomputation and adjoint estimation may be repeated until the estimate $\vert \Delta L_W(\mathbf{x}_h) \vert \simeq \vert -F(\mathbf{x}_h, \mathbf{z}_h - \boldsymbol{\psi}_h) \vert$ is deemed sufficiently small. \\
The above methodology is compatible with the direct approach from Section \ref{subsect:Adaptive_Direct_Minimization}. A typical workflow consists of computing a bijection $\mathbf{x}_h$ under the cost function (\ref{eq:Adaptive_bijectivity_costfunc}) using the PDE-based approach and continuing to improve parametric quality using (\ref{eq:Adaptive_dWinslow_costfunc}). Furthermore, once a bijective $\mathbf{x}_h$ has been found, it may serve as an initial guess for the direct approach from Section \ref{subsect:Adaptive_Direct_Minimization}.

\subsection{Choice of Adjoint Basis}
\label{subsect:Adaptive_Adjoint_Basis}
Problem (\ref{eq:Adaptive_DWR_discrete_adjoint}) requires choosing a suitable dual spline space $\mathcal{V} \supset \boldsymbol{\overline{\mathcal{V}}}_h \neq \mathcal{V}_h$, which typically results from uniformly refining $\mathcal{V}_h$ (in either $h$ or the polynomial degree $p$), leading to a $\sim 4$-fold increase in the number of DOFs associated with the (linear) discrete adjoint equation. In a THB-setting, we have the luxury of choosing $\overline{\mathcal{V}}_h$ reminiscent of the role of $K$-refinement \cite{hughes2005isogeometric} in a structured spline setting. Let $(p, \alpha)$ be the degree and regularity of $\mathcal{V}_h$ (which we assume to be equal in both directions for convenience) and let $\mathcal{T}$ denote the corresponding decomposition of $\hat{\Omega}$ into elements. We define $\overline{\mathcal{V}}_K(\mathcal{V}_h)$ as the richest (dimensionality-wise) THB space of degree $p+1$ and regularity $\alpha + 1$ that is compatible with the elements in $\mathcal{T}$. Typically, we have $\operatorname{dim} \overline{\mathcal{V}}_K(\mathcal{V}_h) \simeq \operatorname{dim} \mathcal{V}_h$. \\
While taking $\overline{\mathcal{V}}_h = \mathcal{V}_{h/2}$ yields more accurate adjoint solutions $\mathbf{z}_h \in (\overline{\mathcal{V}}_h^\circ)^2$, we have found the choice $\boldsymbol{\overline{\mathcal{V}}}_h = \overline{\mathcal{V}}_K(\mathcal{V}_h)$, to be sufficient for refinement based on both (\ref{eq:Adaptive_DWR_quadrature_subset}) and (\ref{eq:Adaptive_dWinslow_costfunc}). As such, solving the discrete adjoint equation becomes a cheap operation.

\subsection{Refinement Strategies}
\label{subsect:Adaptive_Refinement_Strategies}
The decomposition into basis function wise contributions $\mathbf{r}_i(\mathbf{x}_h)$ introduced in (\ref{eq:Adaptive_adjoint_equation_discrete}) is particularly useful in a THB-setting since elementwise refinement may not change the dimension of the underlying THB spline space. In the following, we present several strategies for using $\mathbf{r}(\mathbf{x}_h)$ to mark basis functions $w_i \in \left[ \mathcal{V}_h \right]$ for refinement. We define the vectors $\mathbf{w}$ and $\mathbf{\tilde{r}}$ with
\begin{align}
\label{eq:Adaptive_weighted_residual}
    \mathbf{w}_i = \int_{\hat{\Omega}} w_i \mathrm{d} S \quad \text{and} \quad \mathbf{\tilde{r}}_i = \frac{\mathbf{r}_i}{\mathbf{w}_i}.
\end{align}
Furthermore, we let $\mathbf{\tilde{r}}_\text{max} = \max_i |\mathbf{\tilde{r}}_i|$ and $\mathcal{I} = \{1, \ldots, \vert \mathcal{V}_h \vert \}$. Inspired by \cite{prudhomme2003computable}, we define
\begin{align}
\label{eq:Adaptive_Absolute_Maximum_Threshold}
    \mathcal{I}^\alpha_\text{max} = \{ i \in \mathcal{I} \enskip \vert \enskip | \mathbf{\tilde{r}}_i | \geq \beta \mathbf{\tilde{r}}_\text{max} \}
\end{align}
as the index-set of absolutely weighted contributions that exceed the value $\beta \mathbf{\tilde{r}}_\text{max}$, for some $\beta \in [0, 1]$. The $i \in \mathcal{I}^\alpha_\text{max}$ then constitute the indices corresponding to basis functions whose supporting elements $\mathcal{E}^k \in \mathcal{T}$, from the $k$-th level in the element hierarchy, are replaced by finer counterparts $\mathcal{E}^{k+1}$ from the $(k+1)$-th level. Note that the function may, due to preceding refinements of other functions, be already partially supported by $\mathcal{E}^{l} \in \mathcal{T}$, with $l \geq k + 1$. In this case only the coarsest supporting elements $\mathcal{E}^k$ are refined. As a result, upon constructing the canonical THB-spline space over the refined $\mathcal{T}$, $w_i \in \left[ \mathcal{V}_h \right]$ is replaced by several functions from the next level in the hierarchy, leading to a local increase of the DOFs. Basis function wise refinement ensures that always at least one function is removed from the basis and replaced by several finer ones. Naturally, hierarchical refinement based on THB-splines is a somewhat more involved process than this manuscript suggests. For more details, we refer to \cite{giannelli2012thb, van2020adaptive}. \\
Since both (\ref{eq:Adaptive_bijectivity_costfunc}) and (\ref{eq:Adaptive_dWinslow_costfunc}) are strictly positive quantities of interest, disregarding negative contributions in (\ref{eq:Adaptive_weighted_residual}) is a plausible strategy, too. Heuristically, this strategy mildly reduces the total number of required DOFs until bijectivity is achieved. However, this comes at the expense of a larger number of the required a posteriori refinements, which are limited to typically no more than $3-4$ using (\ref{eq:Adaptive_Absolute_Maximum_Threshold}).

\subsection{Results}
\label{subsect:Adaptive_Results}
To demonstrate the appeal of local refinement made possible by THB splines, in the following, we present parameterizations for the U.S. state of Indiana, the German province of North Rhine-Westphalia and the country of Austria, all of which have complicated boundaries but relatively simple interior. The initial basis $\left[ \mathcal{V}_h \right]$ results from refining an initial grid comprised of $7 \times 7$ elements by the boundaries until the contours of $\Omega$ are approximated sufficiently well. In all cases, $\mathcal{V}_h$ is a bicubic hierarchical space. We take $\overline{\mathcal{V}}_h = \overline{\mathcal{V}}_K(\mathcal{V}_h)$ (see Section \ref{subsect:Adaptive_Adjoint_Basis}) and base refinement on (\ref{eq:Adaptive_Absolute_Maximum_Threshold}) with $\beta = 0.2$. The numerical scheme has been implemented in the open-source finite element library \textit{Nutils} \cite{gertjan_van_zwieten_2020_4071707}.

\begin{figure}[h!]
  \centering
  \subfloat{\includegraphics[valign=c, width=0.4 \linewidth]{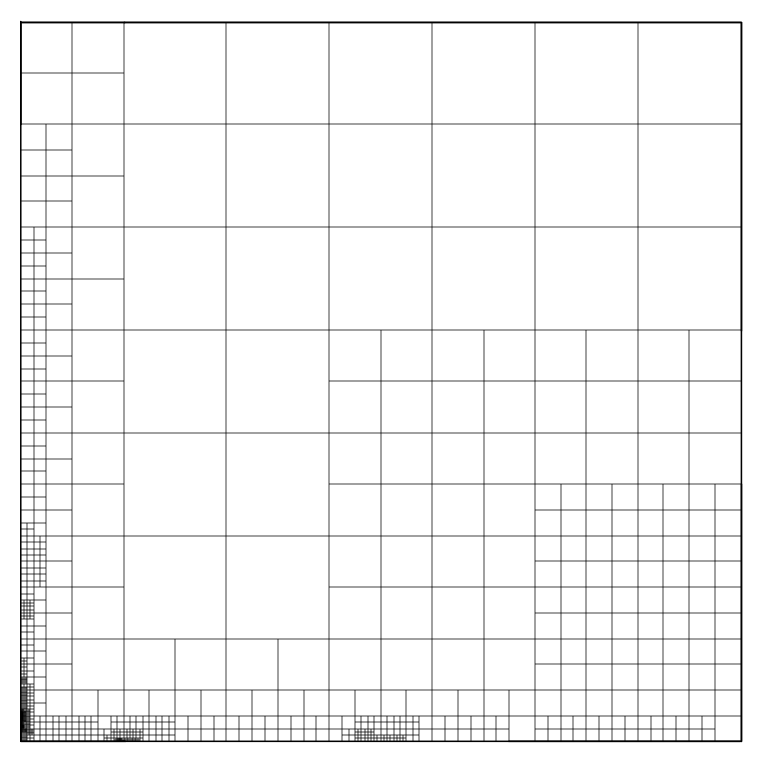}}\label{fig:Indiana_domain} $\quad$
  \centering
  \subfloat{\includegraphics[align=c, width=0.55 \linewidth]{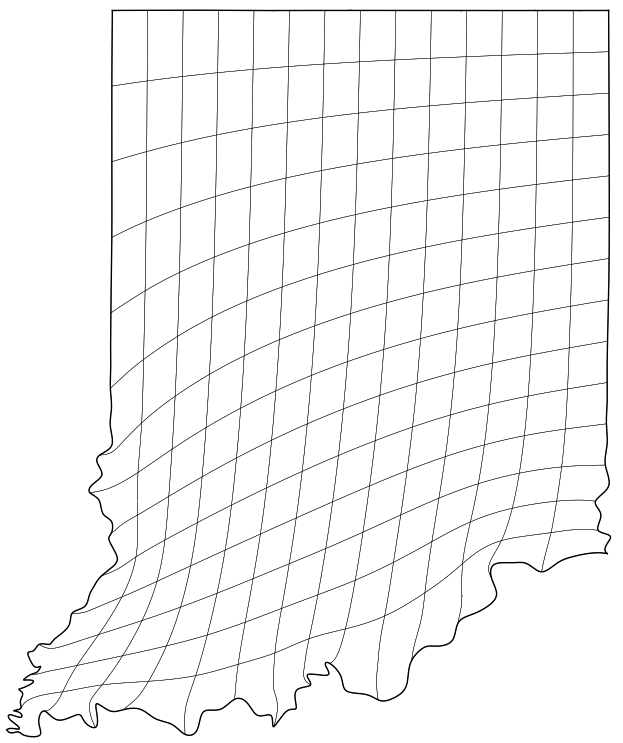}\label{fig:Indiana_final}}
  \caption{The domain with canonical bicubic basis of $2338$ DOFs (left) and the THB-spline parameterization of the U.S. state of Indiana (right).}
  \label{fig:Indiana}
\end{figure}

\begin{figure}[h!]
  \centering
  \subfloat{\includegraphics[align=c, width=0.3 \linewidth]{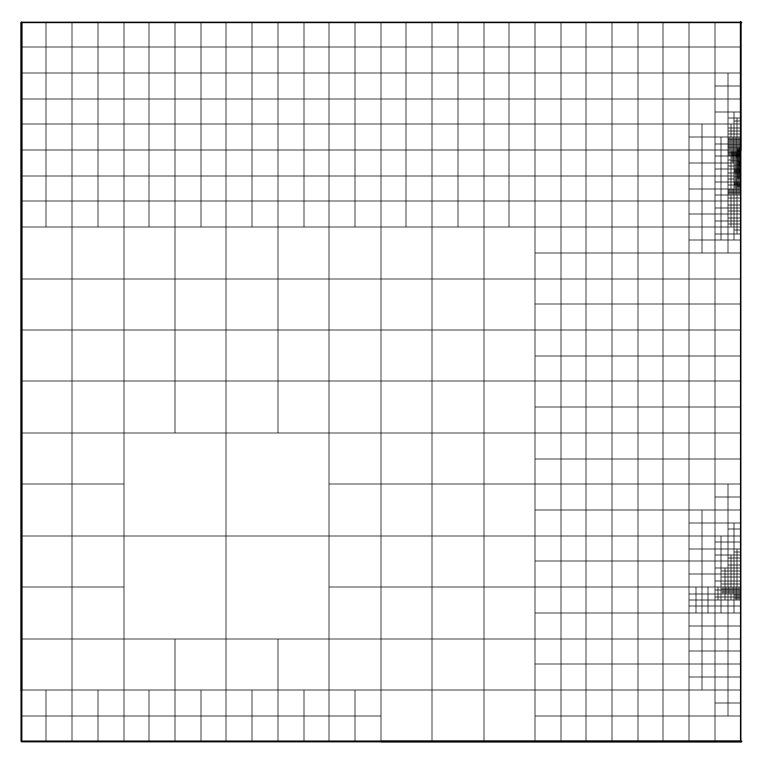}}\label{fig:NRW_domain} $\qquad$
  \centering
  \subfloat{\includegraphics[align=c, width=0.55 \linewidth]{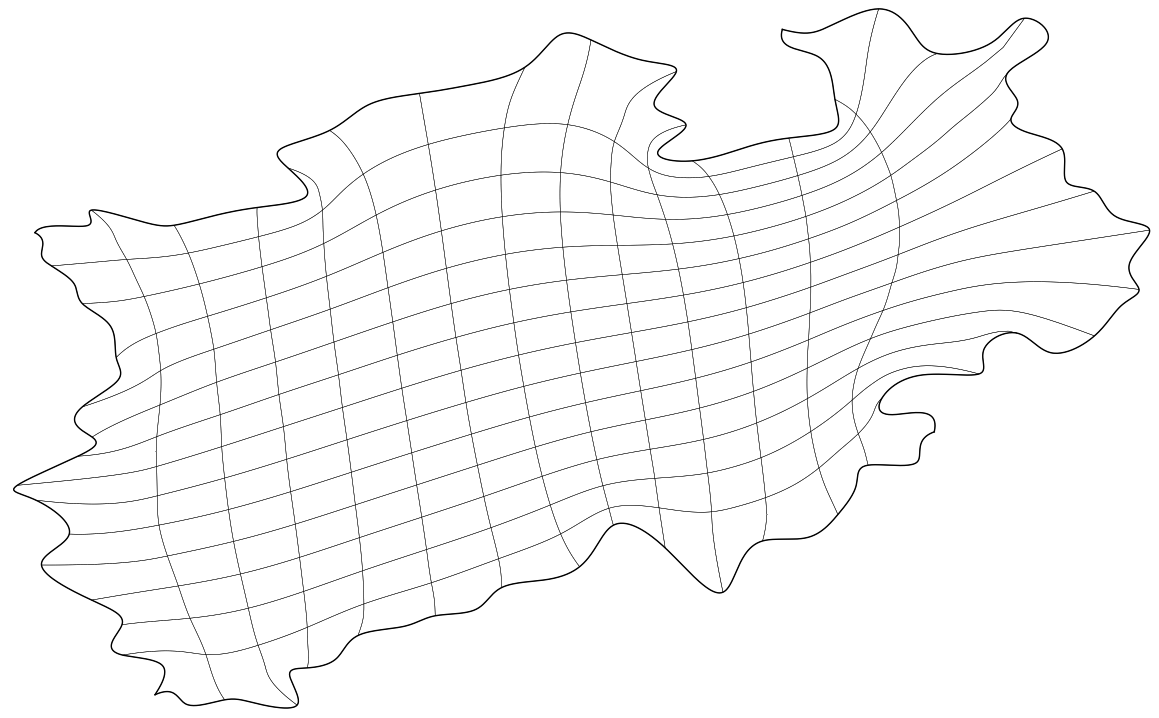}\label{fig:NRW_final}}
  \caption{The domain with bicubic basis of $2676$ DOFs (left) and the THB-spline parameterization of the German province of North Rhine-Westphalia (right).}
  \label{fig:NRW}
\end{figure}

\begin{figure}[h!]
  \centering
  \subfloat{\includegraphics[align=c, width=0.27 \linewidth]{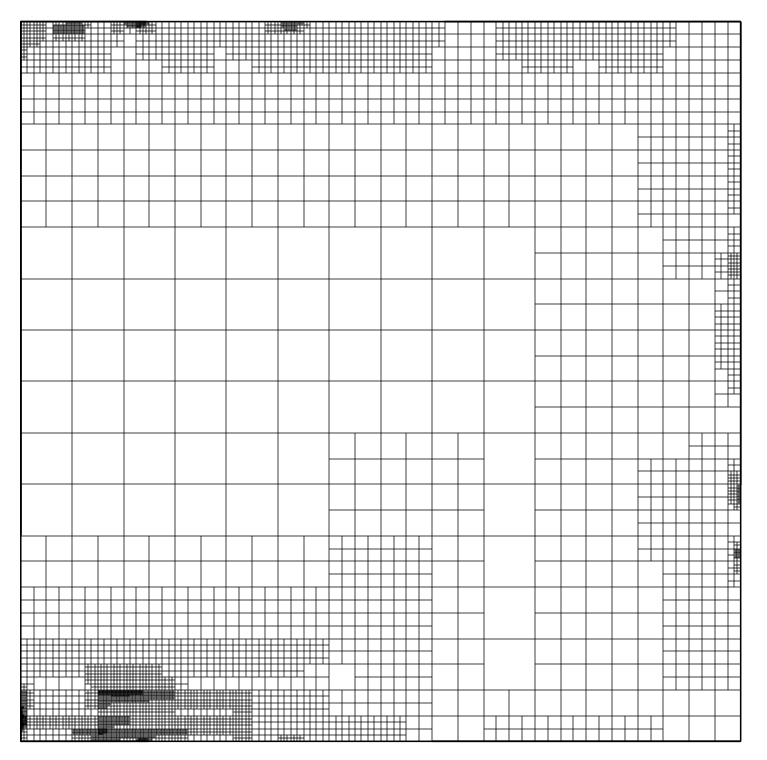}}\label{fig:Austria_domain}
  \centering
  \subfloat{\includegraphics[align=c, width=0.72 \linewidth]{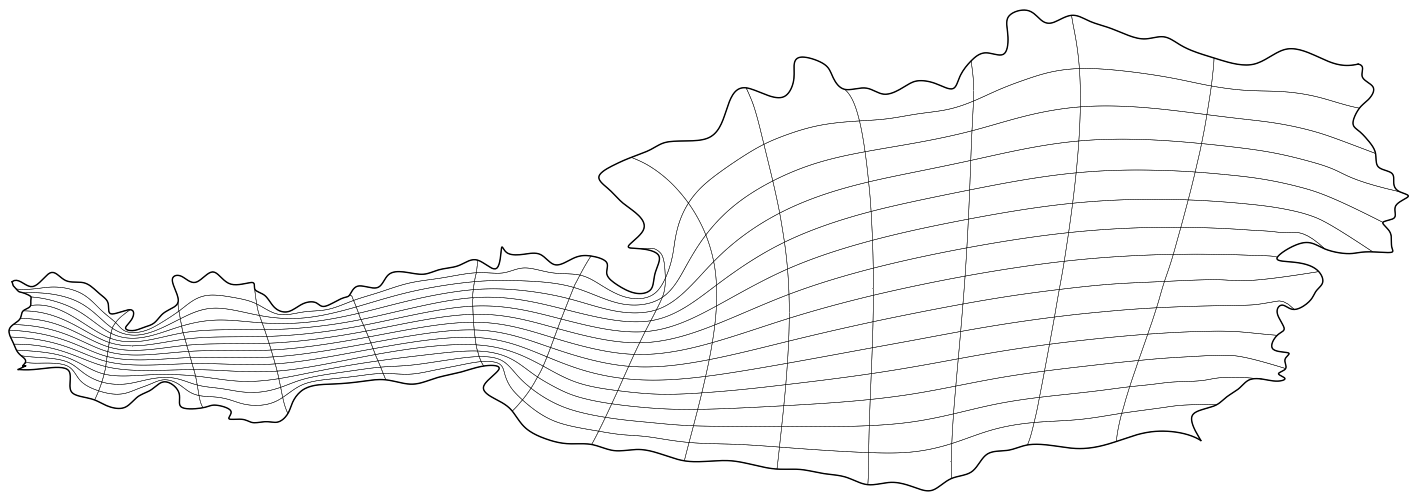}\label{fig:Austria_final}}
  \caption{The domain with bicubic basis comprised of $9640$ DOFs (left) and the THB-spline parameterization of Austria (right).}
  \label{fig:Austria}
\end{figure}

\noindent Figures (\ref{fig:Indiana}) to (\ref{fig:Austria}) clearly demonstrate the DOF savings made possible by local refinement. Not surprisingly, refinement especially affects the protruded and concave areas close to the boundaries. \\
At every refinement level, parameterizations were computed using the Newton-Krylov approach from Section \ref{subsect:Adaptive_Newton}. They were post-processed with the direct approach from Section \ref{subsect:Adaptive_Direct_Minimization} once bijectivity had been achieved. \\
The iterative solver typically converges after $4-5$ nonlinear iterations on the coarsest level plus another $2-3$ iterations per a posteriori refinement. Once bijectivity is achieved, initializing the direct approach from (\ref{eq:Adaptive_Winslow_PDE_discretized}) with the PDE solution typically leads to convergence after fewer than $3$ iterations.
\section{Domain Optimization}
\label{sect:Adaptive_Domain_Optimization}
As demonstrated in Section \ref{subsect:Adaptive_Results}, the approach from Section \ref{sect:Adaptive_Basic_Scheme} can handle challenging geometries. However, it lacks the flexibility of precisely controlling the parametric properties of the outcome, which may lead to undesirable features, such as large elements (see Figures \ref{fig:solution_strategies} and \ref{fig:solution_strategies_ref}). As such, in the following we present a framework that allows for more flexibility, where we pay particular attention to mitigating the aforementioned pathologies associated with inversely harmonic maps. \\
Instead of mapping inversely harmonically into a domain $\hat{\Omega}$ with a Cartesian coordinate system, we now define it through a parameterization $\mathbf{s}: \hat{\Omega} \rightarrow \hat{\Omega}$. For convenience, we assume that the boundary correspondence $\mathbf{s} \vert_{\partial \hat{\Omega}}: \partial \hat{\Omega} \rightarrow \partial \hat{\Omega}$ is the identity. Suppose that $\mathbf{x}^*: \hat{\Omega} \rightarrow \Omega$ solves the equation
\begin{align}
    \Delta_{\mathbf{x}} \boldsymbol{\xi} = \mathbf{0}, \quad \text{s.t.} \quad \mathbf{x} \vert_{\partial \hat{\Omega}} = \mathbf{x}_D(\boldsymbol{\xi}),
\end{align}
for $\mathbf{x}$. Then, if $\mathbf{x}(\boldsymbol{\xi})$ is the solution of
\begin{align}
\label{eq:Adaptive_Strong_Equations_reparameterized}
    \Delta_{\mathbf{x}} \mathbf{s}(\boldsymbol{\xi}) = \mathbf{0}, \quad \text{s.t.} \quad \mathbf{x} \vert_{\partial \hat{\Omega}} = \mathbf{x}_D(\mathbf{s}(\boldsymbol{\xi})),
\end{align}
it clearly satisfies $\mathbf{x} = \mathbf{x}^* \circ \mathbf{s}$, thanks to the fact that $\mathbf{x}_D \circ \mathbf{s} = \mathbf{x}_D$ on $\partial \hat{\Omega}$ (i.e., the boundary condition does not change upon pullback). As such, we may approximate compositions $\mathbf{x}^* \circ \mathbf{s}$ by solving the discretized counterpart of (\ref{eq:Adaptive_Strong_Equations_reparameterized}). \\
Introducing the set of vectors
\begin{align}
    \mathbf{p}^{ij}(\mathbf{s}) = -T^{-1} \frac{\partial^2 \mathbf{s}}{\partial \boldsymbol{\xi}_i \partial \boldsymbol{\xi}_j}, \quad \text{with} \quad T = \partial_{\boldsymbol{\xi}} \mathbf{s} \quad \text{and} \quad (i, j) \in \{1, 2\} \times \{1, 2 \},
\end{align}
it can be shown that with $\mathbf{s} = \mathbf{s}(\boldsymbol{\xi})$, (\ref{eq:Adaptive_Strong_Equations_reparameterized}) can be reformulated as \cite[Chapter 4]{thompson1998handbook}
\begin{align}
A(\mathbf{x}) \colon \left( H(\mathbf{x}_i) + P^1(\mathbf{s}) \frac{\partial \mathbf{x}_i}{\partial \xi} + P^2(\mathbf{s}) \frac{\partial \mathbf{x}_i}{\partial \eta} \right) = 0 \quad i \in \{1, 2\}, \quad \text{s.t.} \quad \mathbf{x} \vert_{\partial \hat{\Omega}} = \mathbf{x}_D \vert_{\partial \hat{\Omega}}.
\end{align}
Here, the matrices $P^1$ and $P^2$ satisfy
\begin{align}
\label{eq:Adaptive_Control_matrices}
    P^k_{ij}(\mathbf{s}) = \mathbf{p}_{k}^{ij}(\mathbf{s}), \quad k \in \{ 1, 2 \}.
\end{align}
Therefore, we introduce
\begin{align}
\label{eq:Adaptive_F_operator_reparameterized}
    F(\mathbf{x}, \boldsymbol{\sigma}, \mathbf{s}) = \sum_{i=1}^2 \int_{\hat{\Omega}} \boldsymbol{\tau}( \boldsymbol{\sigma}, \mathbf{x})_i A(\mathbf{x}) : \underbrace{\left( H(\mathbf{x}_i) + P^1(\mathbf{s}) \frac{\partial \mathbf{x}_i}{\partial \xi} + P^2(\mathbf{s}) \frac{\partial \mathbf{x}_i}{\partial \eta} \right) \det T(\mathbf{s})}_{\tilde{H}(\mathbf{x}_i, \mathbf{s})} \mathrm{d}S,
\end{align}
and for given $\mathbf{s}(\boldsymbol{\xi})$, we solve
\begin{align}
\label{eq:Adaptive_inverse_mapping_approach_PDE_scaled_discretized_reparameterized}
    & \text{find } \mathbf{x}_h \in \mathcal{U}^{\mathbf{x}_D}_h \quad \text{s.t.} \quad F(\mathbf{x}_h, \boldsymbol{\sigma}_h, \mathbf{s})= 0 \quad \forall \boldsymbol{\sigma}_h \in \mathcal{U}^{\mathbf{0}}_h,
\end{align}
in order to approximate $\mathbf{x}^* \circ \mathbf{s}$. Unless stated otherwise, we utilize the Newton approach from Section \ref{subsect:Adaptive_Newton} with $\boldsymbol{\tau}(\boldsymbol{\sigma}, \mathbf{x}) = \boldsymbol{\sigma}$. We can apply the Picard approach from Section \ref{subsect:Adaptive_Picard} by replacing $H(\mathbf{x}_i) \rightarrow \tilde{H}(\mathbf{x}_i, \mathbf{s})$ in equation (\ref{eq:Adaptive_Picard_G_tau}). In the following, we present several strategies for choosing $\mathbf{s}$ to improve the parametric properties of the composite mapping.

\subsection{Exploiting the Maximum Principle}
\label{sect:Adaptive_Eploiting_Maximum_Principle}
Clearly, for well-posedness of (\ref{eq:Adaptive_inverse_mapping_approach_PDE_scaled_discretized_reparameterized}), $\mathbf{s}: \hat{\Omega} \rightarrow \hat{\Omega}$ should not fold. As the control mapping maps into a convex domain, we may exploit the fact that if it is the solution to a second order elliptic problem in divergence form, it is necessarily a bijection \cite{bauman2001univalent}. Thus, let $\mathbf{s} = (\mathbf{s}_1, \mathbf{s}_2)^T$ be such that
\begin{align}
\label{eq:Adaptive_forward_laplace_control_mapping}
    \nabla_{\boldsymbol{\xi}} \cdot \left(D \nabla_{\boldsymbol{\xi}} \mathbf{s}_i \right) = 0 \quad i \in \{1, 2\}, \quad \text{in} \quad \hat{\Omega}, \quad \text{s.t.} \quad \mathbf{s}(\boldsymbol{\xi}) = \boldsymbol{\xi} \text{ on } \partial \hat{\Omega},
\end{align}
where $D: \hat{\Omega} \rightarrow \mathbb{R}^{2 \times 2}$ is an SPD diffusivity tensor. In the following, we assume that an accurate approximation $\mathbf{x}^*_h$ of $\mathbf{x}^*$ has been computed using the methodology from Section \ref{sect:Adaptive_Basic_Scheme}. In order to mitigate the impact of the well-known pathologies of inversely harmonic maps (see Section \ref{sect:Adaptive_Solution_Strategies}), we may select $D$ in (\ref{eq:Adaptive_forward_laplace_control_mapping}) such that the value of
\begin{align}
\label{eq:Adaptive_L_Area}
    L_{\text{Area}}(\mathbf{x}_h) = \int_{\hat{\Omega}} \det J(\mathbf{x}_h)^2 \mathrm{d} S
\end{align}
is expected to decrease with respect to $\mathbf{x}^*$ (see (\ref{eq:Adaptive_Area_costfunc})). Note that
\begin{align}
\label{eq:Adaptive_Poisson_Area_equation}
    \left( \det J (\mathbf{x}^* \circ \mathbf{s}) \right)^2 & \simeq \left( \det \partial_{\mathbf{s}} \mathbf{x}_h^* \right)^2 \det J(\mathbf{s})^2 \nonumber \\
    & = \left( \det \partial_{\mathbf{s}} \mathbf{x}_h^* \right)^2 \left( g_{11} g_{22} - g_{12}^2 \right)_{\boldsymbol{\xi} \rightarrow \mathbf{s}} \nonumber \\
    & \leq \frac{1}{2} \left( \det \partial_{\mathbf{s}} \mathbf{x}_h^* \right)^2 \left( g_{11} + g_{22} \right)^2_{\boldsymbol{\xi} \rightarrow \mathbf{s}},
\end{align}
where the subscript $\boldsymbol{\xi} \rightarrow \mathbf{s}$ indicates that the $g_{ij}$ between brackets refer to the metric induced by $\mathbf{s}(\boldsymbol{\xi})$. Given that $\mathbf{s}(\boldsymbol{\xi}) = \boldsymbol{\xi}$ initially, (\ref{eq:Adaptive_Poisson_Area_equation}) suggests a convex optimization problem of the form
\begin{align}
\label{eq:Adaptive_PoissonArea_minimization}
    L_\text{PoissonArea}(\mathbf{s}, k) \rightarrow \min \limits_{\mathbf{s} \in \mathcal{V}_h^2}, \quad \text{s.t.} \quad \mathbf{s}(\boldsymbol{\xi}) = \boldsymbol{\xi}\text{ on } \partial \hat{\Omega},
\end{align}
where
\begin{align}
\label{eq:Adaptive_PoissonArea_costfunc}
    L_\text{PoissonArea}(\mathbf{s}, k) = \int_{\hat{\Omega}} \left( \det \partial_{\boldsymbol{\xi}} \mathbf{x}_h^* \right)^{k} \left( \| \partial_{\boldsymbol{\xi}} \mathbf{s}_1 \|^2 + \| \partial_{\boldsymbol{\xi}} \mathbf{s}_2 \|^2 \right) \mathrm{d}S,
\end{align}
for recomputing $\mathbf{s}(\boldsymbol{\xi})$. As such, we are solving the discretized equations corresponding to (\ref{eq:Adaptive_forward_laplace_control_mapping}) with
\begin{align}
\label{eq:Adaptive_Diffusivity_Area}
D = \left( \det \partial_{\boldsymbol{\xi}} \mathbf{x}_h^* \right)^{k} I^{2 \times 2}.
\end{align}
Even though the exact solution of (\ref{eq:Adaptive_forward_laplace_control_mapping}) does not fold, the discretized counterpart may fold due to extreme diffusive anisotropy. This can be counteracted by reducing the value of $k$. Alternatively, (\ref{eq:Adaptive_PoissonArea_costfunc}) can be utilized for DWR-based a posteriori refinement to achieve bijectivity and accuracy of $\mathbf{s}: \hat{\Omega} \rightarrow \hat{\Omega}$. \\
Upon completion, we compute $\mathbf{x}_h \in \mathcal{V}_h^2$ using the control mapping $\mathbf{s}: \hat{\Omega} \rightarrow \hat{\Omega}$, with a posteriori refinement if necessary.
\begin{figure}[h!]
  \centering
  \subfloat[]{\includegraphics[width=.45\linewidth]{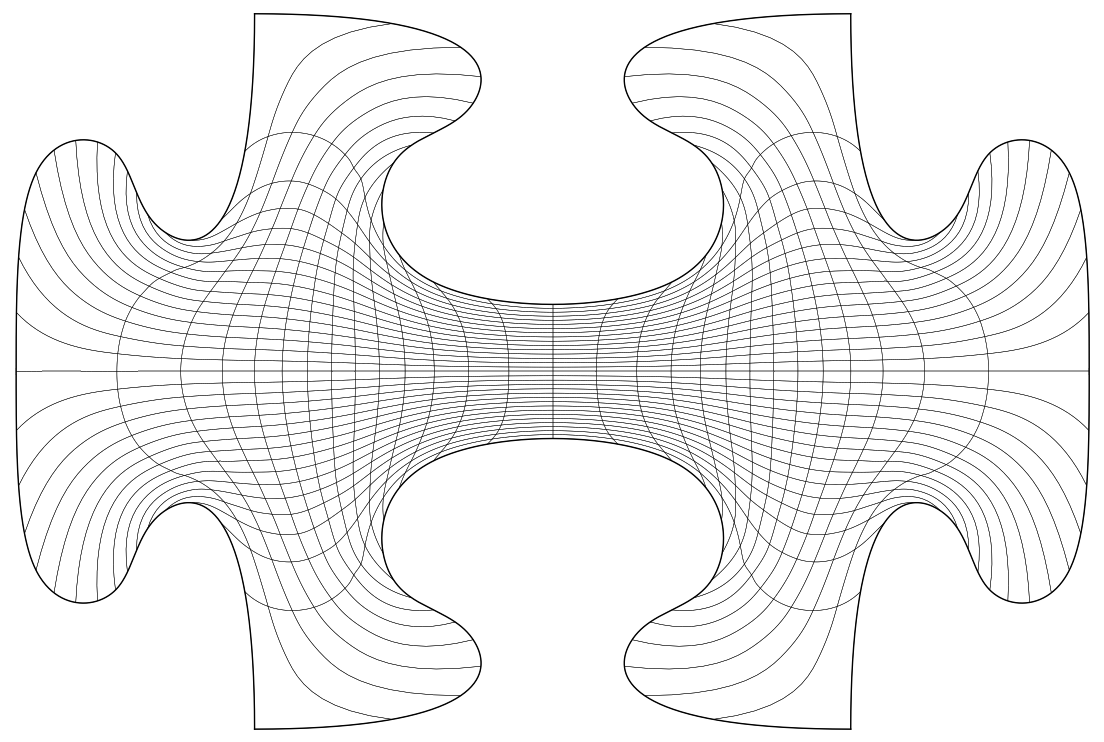}\label{fig:PoissonAreaPuzzle_k_0}} \quad
  \centering
  \subfloat[]{\includegraphics[width=.45\linewidth]{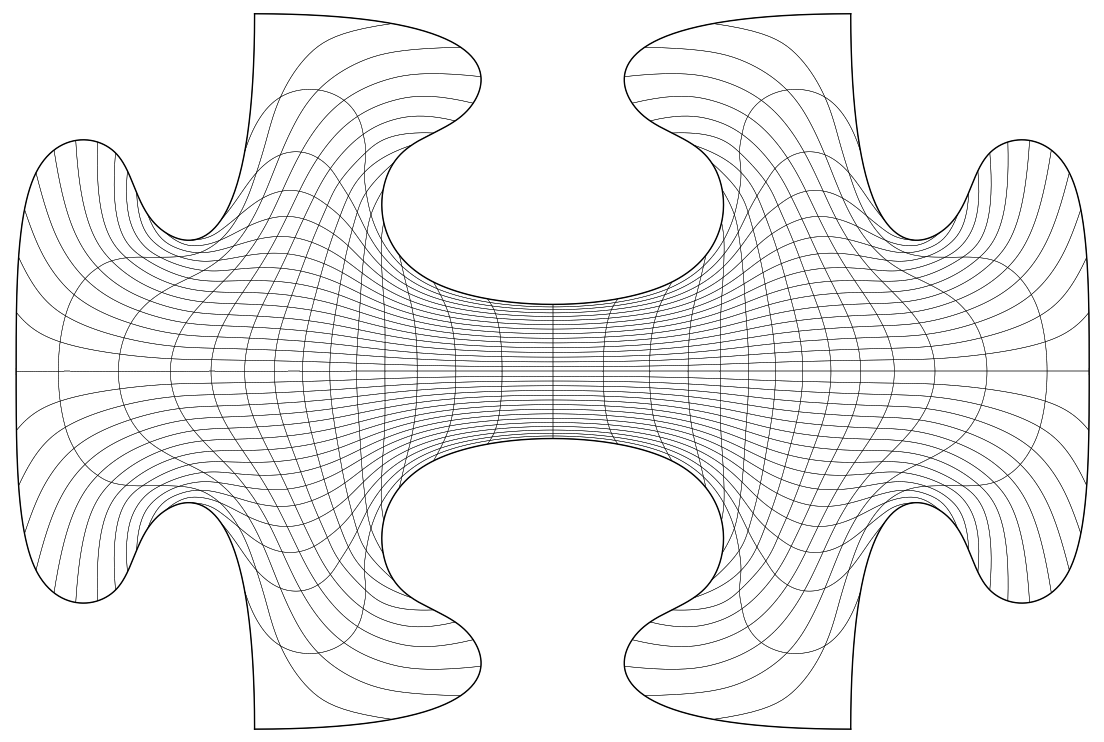}\label{fig:PoissonAreaPuzzle_k_05}} \\
  \centering
  \subfloat[]{\includegraphics[width=.45\linewidth]{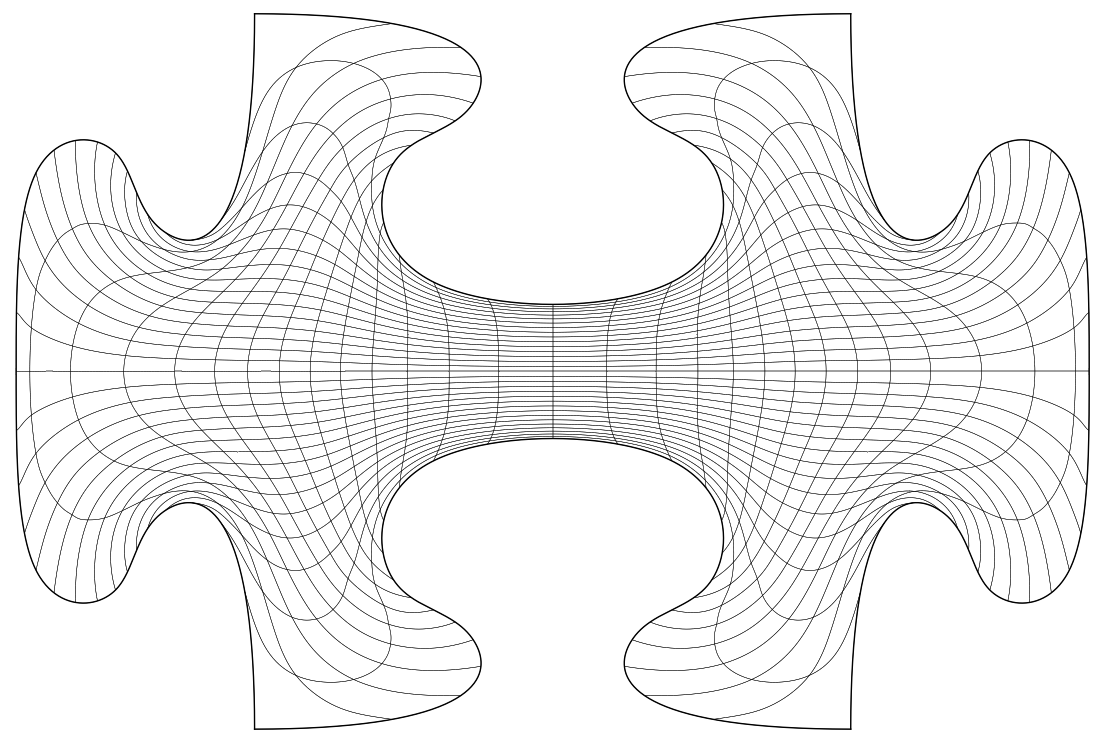}\label{fig:PoissonAreaPuzzle_k_1}} \quad
  \centering
  \subfloat[]{\includegraphics[width=.45\linewidth]{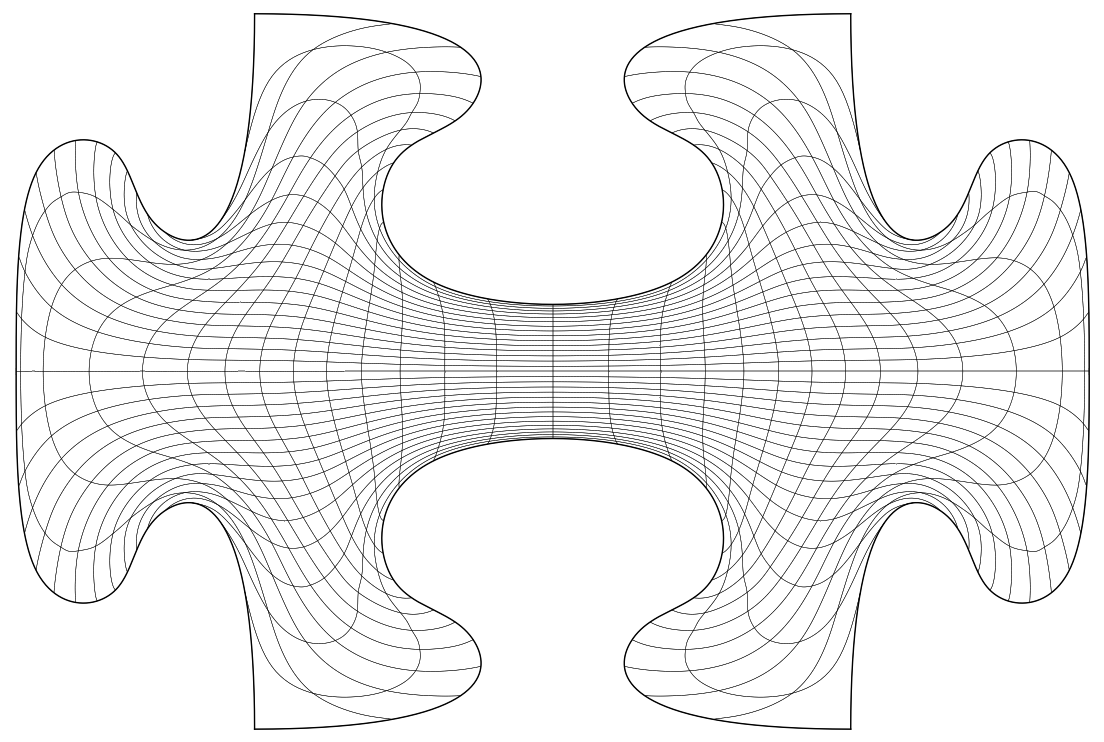}\label{fig:PoissonAreaPuzzle_k_15}}
  \centering
  \caption{Several parameterizations of the puzzle piece with reparameterization based on (\ref{eq:Adaptive_forward_laplace_control_mapping}) and (\ref{eq:Adaptive_Diffusivity_Area}) with the reference parameterization $k=0$ (a), reparameterization with $k=0.5$ (b), $k=1$ (c) and $k=1.5$ (d). }
  \label{fig:PoissonAreaPuzzle}
\end{figure}

\begin{figure}[h!]
  \centering
  \subfloat[]{\includegraphics[align=c, width=.45\linewidth]{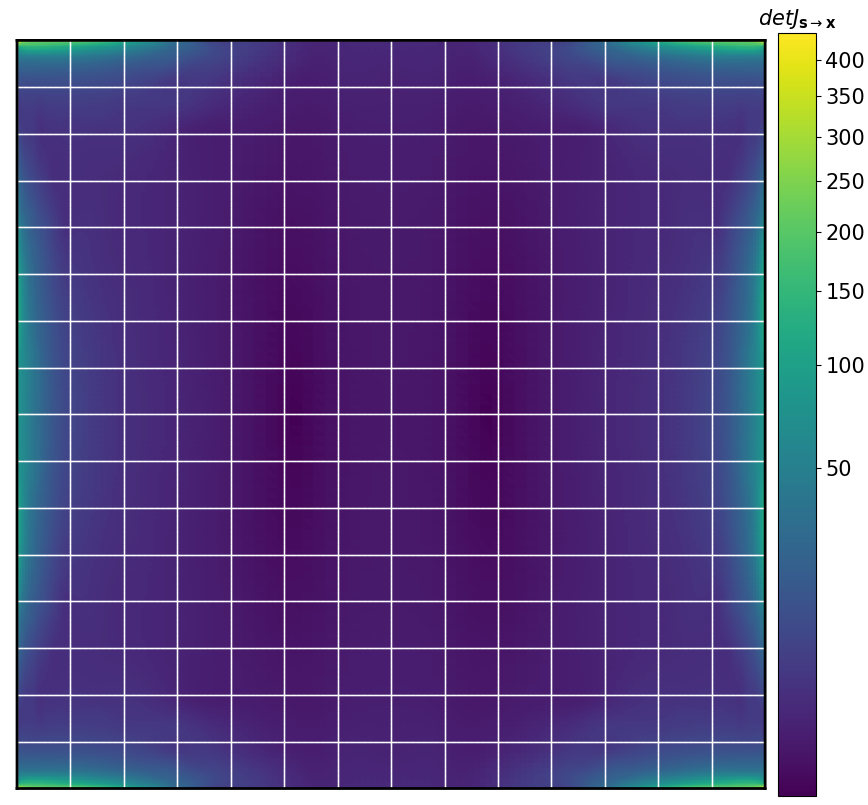}\label{fig:PoissonAreaPuzzle_ReparameterizedDomain_normal}} \quad
  \centering
  \subfloat[]{\includegraphics[align=c, width=.45\linewidth]{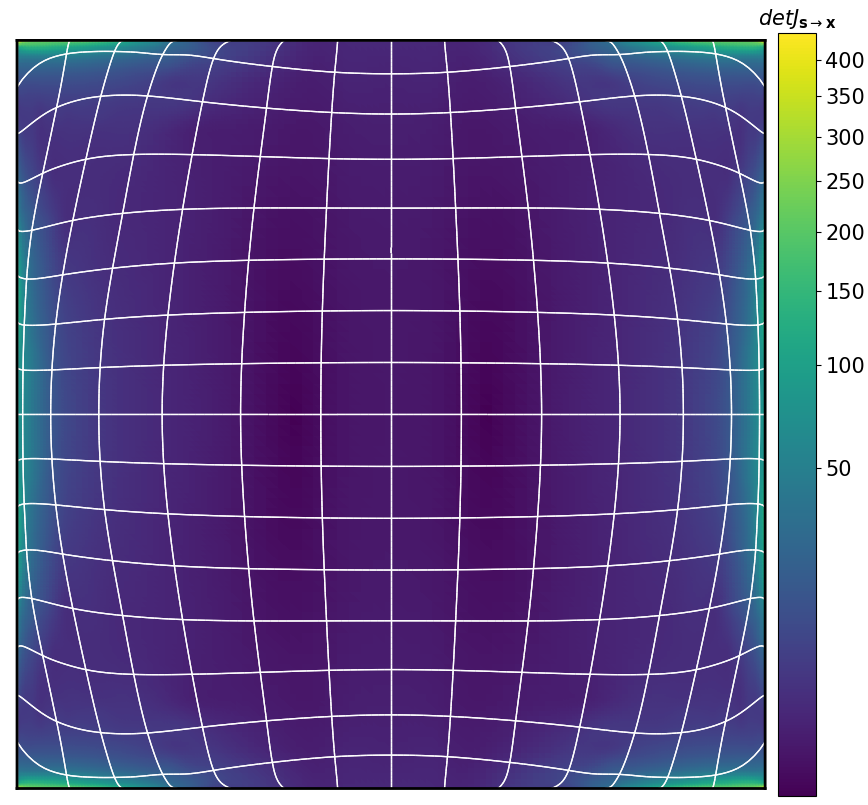}\label{fig:PoissonAreaPuzzle_ReparameterizedDomain_ReparameterizedDomain}}
  \caption{Plots showing the reference domain (a) and the reparameterized domain based on (\ref{eq:Adaptive_forward_laplace_control_mapping}) and (\ref{eq:Adaptive_Diffusivity_Area}) with $k=1.5$ (b). The figure clearly shows that the elements are contracted wherever $\det \partial_{\boldsymbol{\xi}} \mathbf{x}_h^*$ is large.}
  \label{fig:PoissonAreaPuzzle_ReparameterizedDomain}
\end{figure}

\begin{table}[h!]
\centering
\begin{tabular}{c|c|c|c|c}
$k$   & $0$ & $0.5$ & $1$ & $1.5$ \\ \hline
$L_\text{Area}(\mathbf{x}_h) \times 10^{-2}$ & $3.291$ & $2.077$ & $1.439$ & $1.299$
\end{tabular}
\captionsetup{justification=centering}
\caption{Evaluation of $L_\text{Area}(\mathbf{x}_h)$ for various values of k.}
\label{tab:L_PoissonArea_Puzzle_Piece}
\end{table}
\noindent Figure \ref{fig:PoissonAreaPuzzle} shows puzzle piece geometry parameterizations for various values of $k$, while Table \ref{tab:L_PoissonArea_Puzzle_Piece} contains the outcomes of substituting into (\ref{eq:Adaptive_L_Area}). Both clearly demonstrate that the methodology has the desired effect, with more drastic outcomes for larger values of $k$. Figure \ref{fig:PoissonAreaPuzzle_ReparameterizedDomain} shows the isolines of $\mathbf{s}(\boldsymbol{\xi})$ before and after reparameterization with $k=1.5$. All parameterizations were computed with the reference basis corresponding to Figure \ref{fig:PoissonAreaPuzzle_k_0}. No a posteriori refinements were necessary.

\begin{figure}[h!]
  \centering
  \subfloat[]{\includegraphics[align=c, width=.30\linewidth]{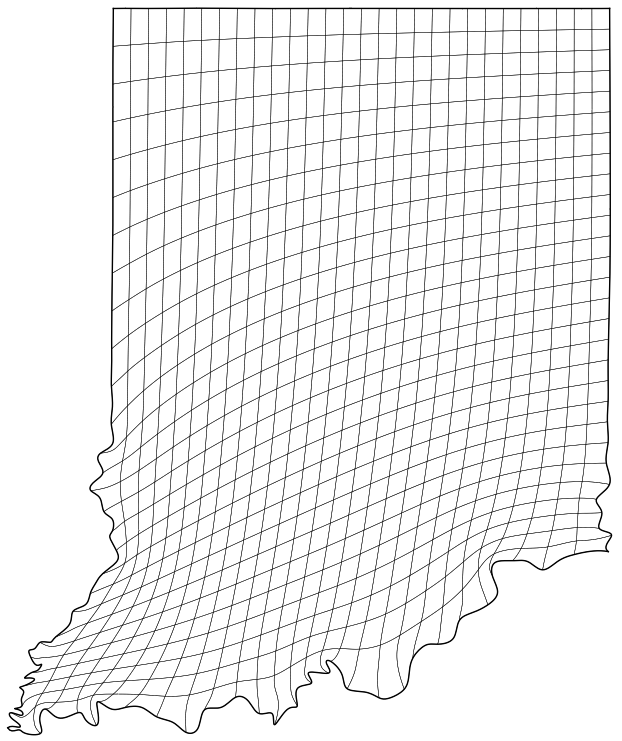}\label{fig:PoissonAreaIndiana_k_0}}
  \centering
  \subfloat[]{\includegraphics[align=c, width=.30\linewidth]{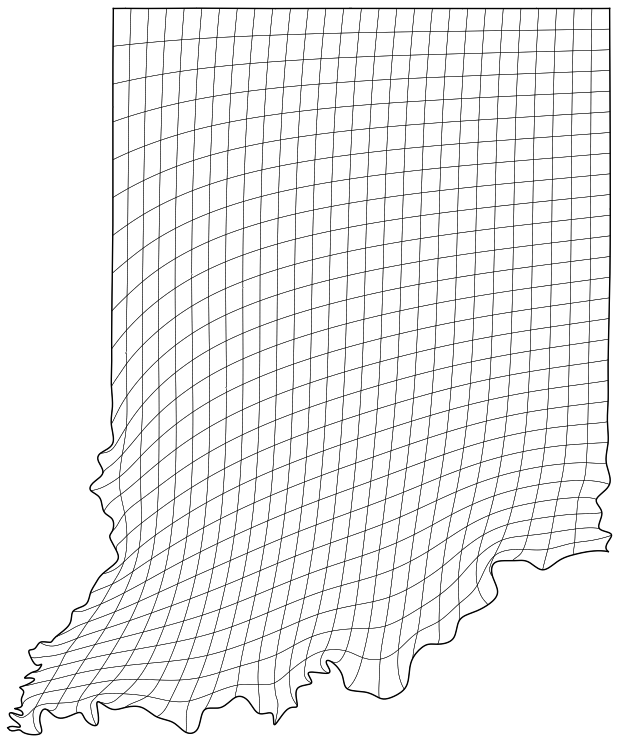}\label{fig:PoissonAreaIndiana_k_1}} \qquad
  \centering
  \subfloat[]{\includegraphics[align=c, width=.30\linewidth]{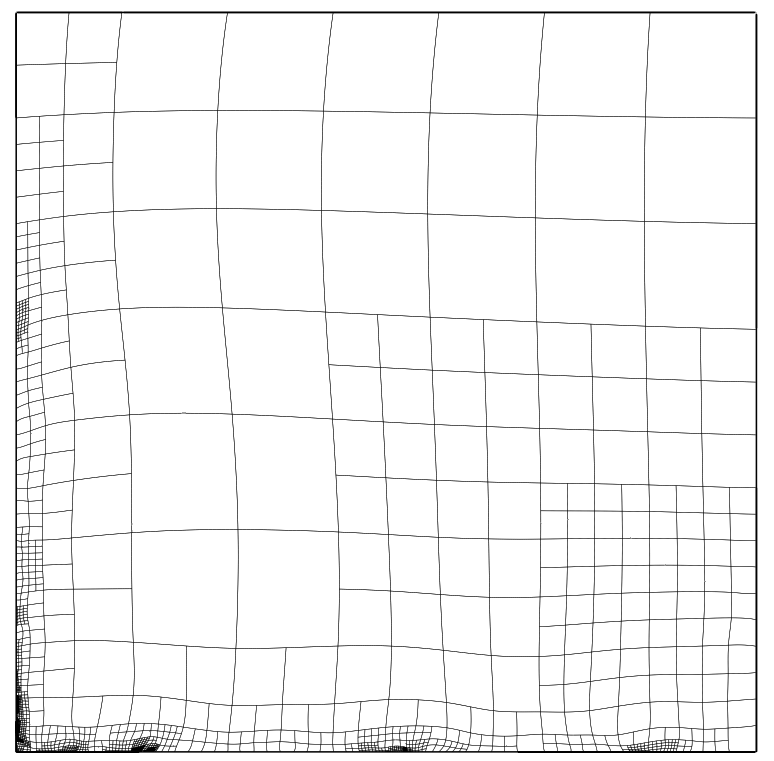}\label{fig:PoissonAreaIndiana_k_1_domain}}
  \caption{Parameterization of the U.S. state of Indiana with $k=0$ (a), $k=1$ (b) and the corresponding reparameterized domain (c). }
  \label{fig:PoissonAreaIndiana}
\end{figure}
Figure \ref{fig:PoissonAreaIndiana} shows parameterizations of the U.S. state of Indiana for $k=0$ and $k=1$. Contrary to Table \ref{tab:L_PoissonArea_Puzzle_Piece}, with
\begin{align*}
    L_\text{Area}(\mathbf{x}_h) = 1.049 \times 10^2 \quad \text{for} \quad k=0 \quad \text{and} \quad L_\text{Area}(\mathbf{x}_h) = 1.008 \times 10^2 \quad \text{for} \quad k=1,
\end{align*}
the effect is very mild. Restricting the integrals to $\eta < 1/7$, however, the difference becomes more pronounced with
\begin{align*}
    L_\text{Area}(\mathbf{x}_h) = 17.167 \quad \text{and} \quad L_\text{Area}(\mathbf{x}_h) = 14.118,
\end{align*}
respectively. Unsurprisingly from the shape of the geometry, the difference is most striking close to the lower boundary, which can also be seen in the figure. A posteriori refinement was necessary in Figure \ref{fig:PoissonAreaIndiana_k_1}. \\
Heuristically, reparameterization based on the maximum principle proves to be one of the most robust and effective choices for a wide range of geometries while being computationally efficient. This is thanks to the fact that it addresses the known pathologies of inversely harmonic maps, while also yielding smooth solutions, which preserves smoothness of the composite mapping.

\subsection{Constrained Domain Optimization}
The concept of reparameterizing the domain in order to alter the parametric properties of the recomputed geometry parameterization can be further extended in a way more reminiscent of the well-known cost function minimization approach (see Section \ref{sect:Adaptive_Introduction}). Given an accurate approximation $\mathbf{x}^*_h$ of $\mathbf{x}^*$ (see Section \ref{sect:Adaptive_Eploiting_Maximum_Principle}), we define the metric $G_{\mathbf{s} \rightarrow \mathbf{x}} = \partial_{\mathbf{s}} \mathbf{x}^T \partial_{\mathbf{s}} \mathbf{x}$, which is initially given by 
\begin{align*}
    G_{\mathbf{s} \rightarrow \mathbf{x}} = \partial_{\boldsymbol{\xi}} \mathbf{x}^{*T} \partial_{\boldsymbol{\xi}} \mathbf{x}^{*} \simeq \partial_{\boldsymbol{\xi}} \mathbf{x}^{*T}_h \partial_{\boldsymbol{\xi}} \mathbf{x}^{*}_h.
\end{align*}
Hence, in order to optimize $\mathbf{x}_h(\boldsymbol{\xi})$, we optimize $\mathbf{s}(\boldsymbol{\xi})$ in the metric induced by $G_{\mathbf{s} \rightarrow \mathbf{x}}$. With
\begin{align}
    g_{ij}^\mathbf{s} = \left[\partial_{\boldsymbol{\xi}} \mathbf{s}^T G_{\mathbf{s} \rightarrow \mathbf{x}} \partial_{\boldsymbol{\xi}} \mathbf{s} \right]_{ij} \quad \text{and} \quad J_{ij}^\mathbf{s} = [\partial_{\boldsymbol{\xi}} \mathbf{x}_h^* \partial_{\boldsymbol{\xi}} \mathbf{s}]_{ij},
\end{align}
we define domain optimization cost functions $Q_i^\mathbf{s}(\mathbf{s})$ by replacing $g_{ij} \rightarrow g_{ij}^{\mathbf{s}}$ and $J_{ij} \rightarrow J_{ij}^\mathbf{s}$ in the $Q_i$ introduced in equation (\ref{eq:Adaptive_weighted_sum}) (see Section \ref{sect:Adaptive_Introduction}). We may nevertheless choose to add terms of the form $Q_j(\mathbf{s})$, which should then be regarded as regularization terms. Let $\mathcal{U}_h^\square = \{ \mathbf{v} \in \mathcal{V}_h^2 \enskip \vert \enskip \mathbf{v} = \boldsymbol{\xi} \text{ on } \partial \hat{\Omega} \}$. A domain optimization problem takes the form
\begin{align}
\label{eq:Adaptive_Domain_Optimization_Problem}
    \int_{\hat{\Omega}} Q(\mathbf{s}) \mathrm{d}S \rightarrow \min \limits_{\mathbf{s} \in \mathcal{U}_h^\square}, \quad \text{s.t.} \quad \mathbf{C}(\mathbf{s}) \geq \mathbf{0},
\end{align}
with
\begin{align}
\label{eq:Adaptive_Domain_Q}
    Q(\mathbf{s}) = \sum_i \lambda_i^{\mathbf{s}} Q_i^{\mathbf{s}}(\mathbf{s}) + \sum_j \lambda_j Q_j(\mathbf{s}).
\end{align}
Here, the constraint $\mathbf{C}(\mathbf{s}) \geq 0$ ensures that the minimizer of (\ref{eq:Adaptive_Domain_Optimization_Problem}) does not fold. In the following, we list all choices of $\mathbf{C}(\mathbf{s})$ that come to mind. \\
Given the element segmentation $\mathcal{T}$ of $\hat{\Omega}$, by $\mathcal{V}^{p, \alpha}(\mathcal{T})$ we denote the canonical THB-space with order $p$ and regularity $\alpha$ that is compatible with $\mathcal{T}$. Note that $\alpha \leq p - 1$. Clearly, if $\mathcal{V}_h$ has order $p$ and regularity $\alpha \leq p-1$, this implies that 
\begin{align*}
    \det \partial_{\boldsymbol{\xi}} \mathbf{s} \in \mathcal{V}^{2p -1, \alpha-1}(\mathcal{T}).
\end{align*}
As such, we also have 
\begin{align*}
    \det \partial_{\boldsymbol{\xi}} \mathbf{s} \in \mathcal{V}^{2p-1, -1}(\mathcal{T}).
\end{align*}
Hence, we can base the constraint on B\'ezier extraction, in which we require that all weights of projecting $\det \partial_{\boldsymbol{\xi}} \mathbf{s}$ onto $\mathcal{V}^{2p-1, -1}(\mathcal{T})$ be positive. Let $\mathbf{\hat{d}}$ be the corresponding vector of weights. We have
\begin{align}
\label{eq:Adaptive_Bezier_Constraint}
    \mathbf{\hat{d}}(\mathbf{s}) = \hat{M}^{-1} \mathbf{\hat{f}}(\mathbf{s}) > \mathbf{0}, \quad \text{where} \quad \mathbf{\hat{f}}_i(\mathbf{s}) = \int \limits_{\hat{\Omega}} \hat{\phi}_i \det \partial_{\boldsymbol{\xi}} \mathbf{s} \mathrm{d}S,
\end{align}
with
\begin{align}
\hat{\phi}_i \in \left[ \mathcal{V}^{2p-1, -1}(\mathcal{T}) \right] \quad \text{and} \quad \hat{M}_{i,j} = \int \limits_{\hat{\Omega}} \hat{\phi}_i \hat{\phi}_j \mathrm{d}S.
\end{align}
Note that $\hat{M}$ is block-diagonal with $\vert \mathcal{T} \vert$ blocks of size $(2p, 2p)$. Hence, we computationally efficiently assemble $\hat{M}^{-1}$ simply by computing the inverse of all separate blocks leading to a sparse block-diagonal matrix. As such, the computational costs of testing whether the condition $\mathbf{\hat{d}} > \mathbf{0}$ is fulfilled reduces to the assembly of $\mathbf{\hat{f}}$ along with one sparse matrix-vector multiplication. Assembly of the constraint gradient of $\mathbf{\hat{d}}(\mathbf{c}_{\mathcal{I}})$, where $\mathbf{c}_{\mathcal{I}}$ is a vector containing the inner control points of $\mathbf{s}$, requires the assembly of $\partial_{\mathbf{c}_{\mathcal{I}}} \mathbf{\hat{f}}$ and a sparse matrix-matrix multiplication. The assembly is hence feasible. However, for large values of $p$ this may lead to an infeasibly large number of constraints. \\
Inspired by \cite{gravesen2012planar}, we formulate an alternative constraint by projecting $\det \partial_{\boldsymbol{\xi}} \mathbf{s}$ onto the coarser THB-space $\mathcal{V}^{2p-1, \alpha-1}(\mathcal{T})$. Similar to (\ref{eq:Adaptive_Bezier_Constraint}), this leads to a constraint of the form
\begin{align}
\label{eq:Adaptive_Gravesen_Constraint}
    \mathbf{d}(\mathbf{s}) = M^{-1} \mathbf{f}(\mathbf{s}) > \mathbf{0}, \quad \text{where} \quad \mathbf{f}_i(\mathbf{s}) = \int \limits_{\hat{\Omega}} \phi_i \det \partial_{\boldsymbol{\xi}} \mathbf{s} \mathrm{d}S, \end{align}
with
\begin{align}
\phi_i \in \left[ \mathcal{V}^{2p-1, \alpha - 1}(\mathcal{T}) \right] \quad \text{and} \quad M_{i,j} = \int \limits_{\hat{\Omega}} \phi_i \phi_j \mathrm{d}S.
\end{align}
Increasing the values of $p$ and $\alpha$, unlike for (\ref{eq:Adaptive_Bezier_Constraint}), the length of $\mathbf{d}$ in (\ref{eq:Adaptive_Gravesen_Constraint}) increases only slowly (thanks to $K$-refinement). On the other hand, the matrix $M$ is not block-diagonal and neither is it separable (unlike in a structured spline setting). As such, the assembly of the constraint gradient is prohibitively expensive. Here, a remedy is to introduce the vector of slack variables $\mathbf{e} > 0$. The constraint from (\ref{eq:Adaptive_Gravesen_Constraint}) can be reformulated as follows:
\begin{align}
    C_{\alpha}(\mathbf{s}, \mathbf{e}) = \mathbf{f}(\mathbf{s}) - M \mathbf{e} = \mathbf{0}, \quad \text{with} \quad \mathbf{e} > \mathbf{0}. 
\end{align}
Hence, we avoid inversion with $M$ at the expense of introducing an additional inequality constraint and changing the existing inequality constraint to an equality constraint. Note that we have:
\begin{align}
    \frac{\partial C_{\alpha}(\mathbf{s}, \mathbf{e})}{\partial (\mathbf{c}_{\mathcal{I}}, \mathbf{e})} = \left[ \frac{\partial \mathbf{f}}{\partial \mathbf{c}_{\mathcal{I}}}, -M \right] \quad \text{and} \quad \frac{\partial \mathbf{e}}{\partial (\mathbf{c}_{\mathcal{I}}, \mathbf{e})} = \left[ 0, I \right],
\end{align}
where $I$ denotes the identity matrix of appropriate dimension. 
\begin{figure}[h!]
\centering
  \subfloat[Depiction of the cones associated with the linear constraint $\mathbf{C}_L(\mathbf{x}_h)$ generated from the control net of a structured spline mapping for a bat-shaped geometry. The constraint is violated (the two associated cones intersect) despite the bijectivity of the mapping, demonstrating the restrictiveness of $\mathbf{C}_L(\mathbf{x}_h)$.]{\includegraphics[align=c, width=.9\linewidth]{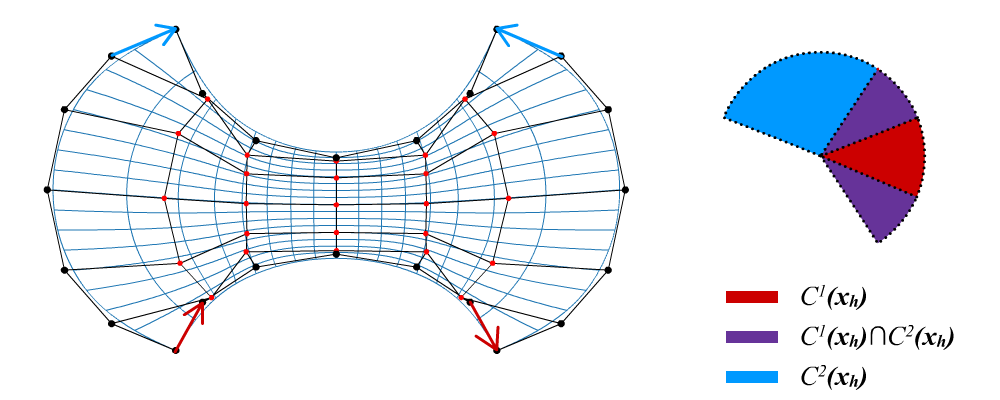}\label{fig:linear_constraint_bat}} \\
  \subfloat[Control net of the Cartesian parametric domain. In this case, both cones collapse into half rays generated by $\mathbb{R}^+(1, 0)^T$ and $\mathbb{R}^+(0, 1)^T$, respectively. The feasible space may then be comprised of all parameterizations of $\hat{\Omega}$ with cones contained within $-\pi/4 < \theta < \pi/4$ and $\pi/4 < \theta < 3\pi/4$, respectively. Hereby, the initial guess is located exactly in the center of the feasible space.]{\includegraphics[align=c, width=.9\linewidth]{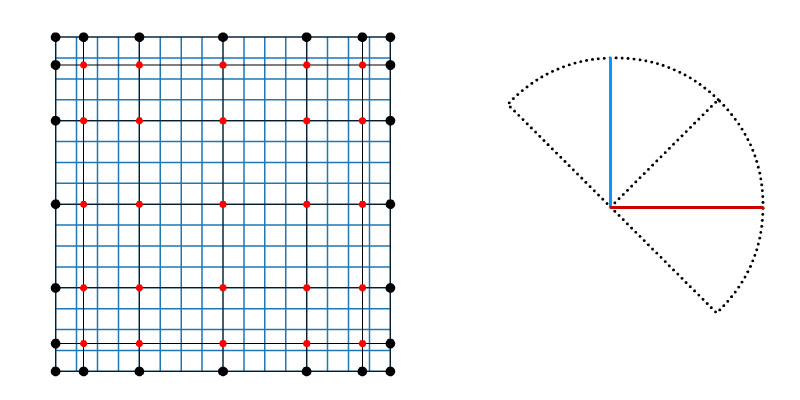}\label{fig:linear_constraint_domain}}
  \caption{Depiction of the bijectivity constraint $\mathbf{C}_L$ for a generic geometry (a), and the initially Cartesian parametric domain (b).}
  \label{fig:linear_constraint}
\end{figure}


\begin{figure}[h!]
\centering
  \subfloat[]{\includegraphics[align=c, width=.45\linewidth]{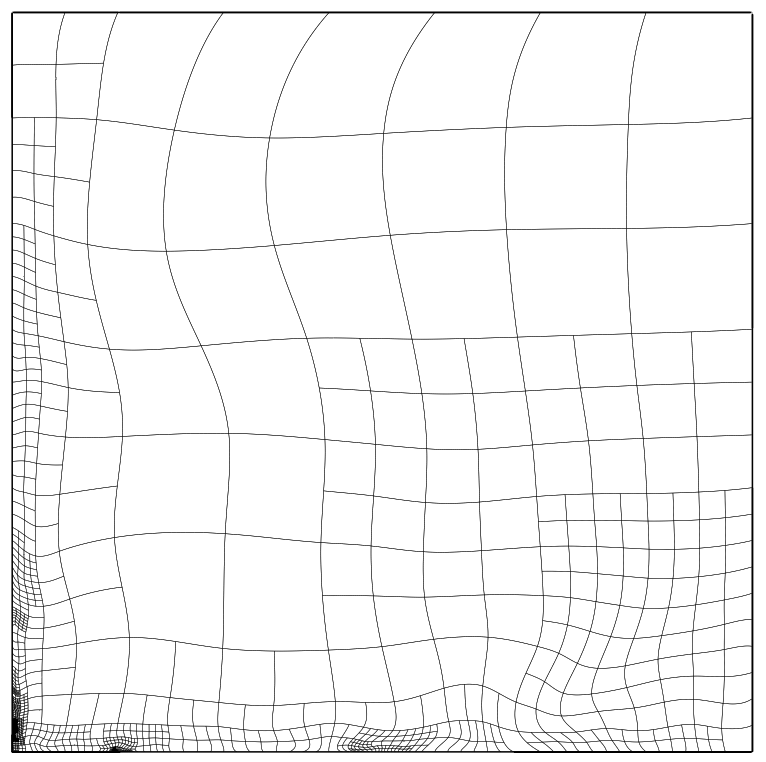}\label{fig:Indiana_AO_domain_reparameterized}} $\qquad$
  \subfloat[]{\includegraphics[align=c, width=.4\linewidth]{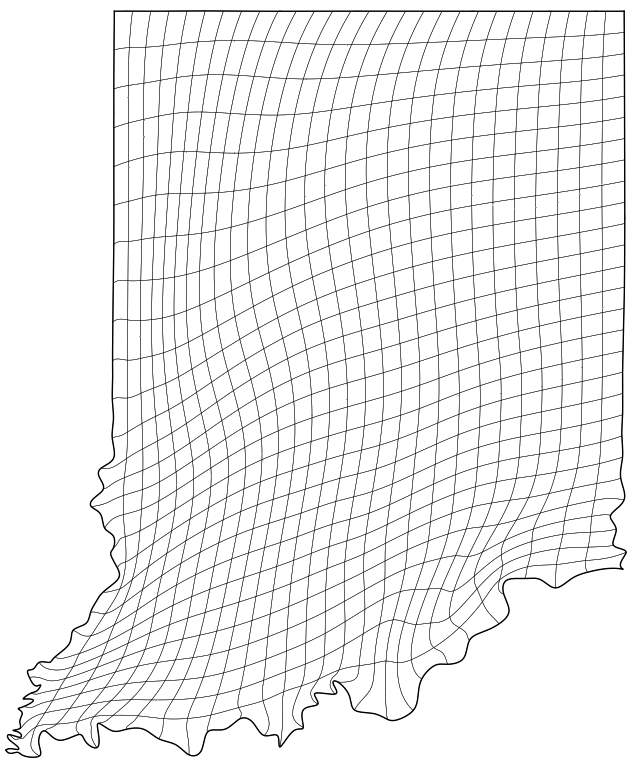}\label{fig:Indiana_AO_geometry}}
  \caption{Result of reparameterizing the domain corresponding to the U.S. state of Indiana with $Q = Q_{\text{AreaOrthogonality}}^{\mathbf{s}}$ (a) and the resulting recomputed geometry parameterization (b).}
  \label{fig:Indiana_AO}
\end{figure}

\begin{figure}[h!]
\centering
  \subfloat[]{\includegraphics[align=c, width=.3\linewidth]{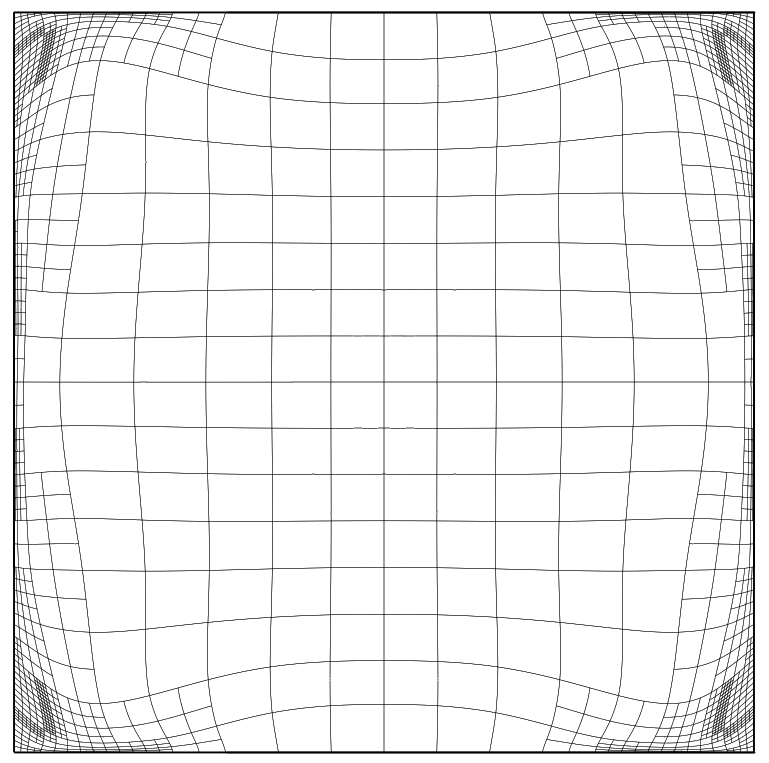}\label{fig:Puzzle_detJ_domain_linear_refined_reparameterized}} $\qquad$
  \subfloat[]{\includegraphics[align=c, width=.5\linewidth]{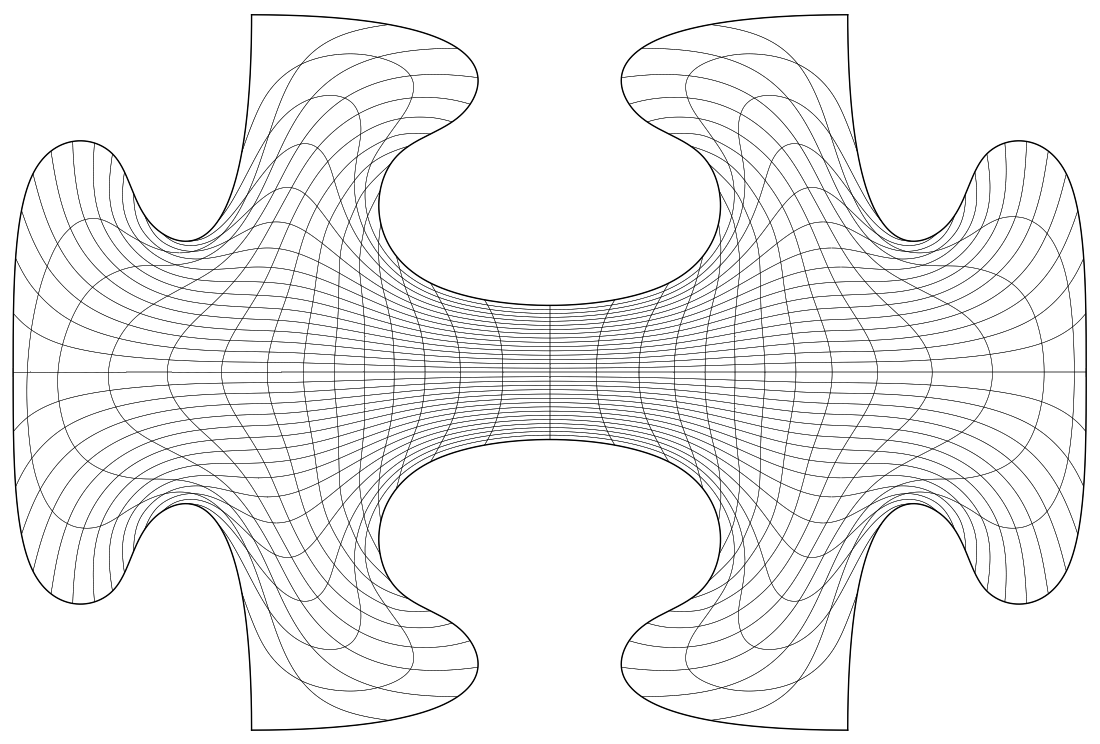}\label{fig:Puzzle_detJ_linear_geometry}}
  \caption{Result of optimizing the puzzle piece domain with $Q = Q_{\text{Area}}^{\mathbf{s}}$ under the constraint $\mathbf{C}(\mathbf{s}) = \mathbf{C}_L(\mathbf{s})$ (a) and the corresponding recomputed mapping (b).}
  \label{fig:Puzzle_detJ_linear}
\end{figure}

\begin{figure}[h!]
\centering
  \subfloat[]{\includegraphics[align=c, width=.3\linewidth]{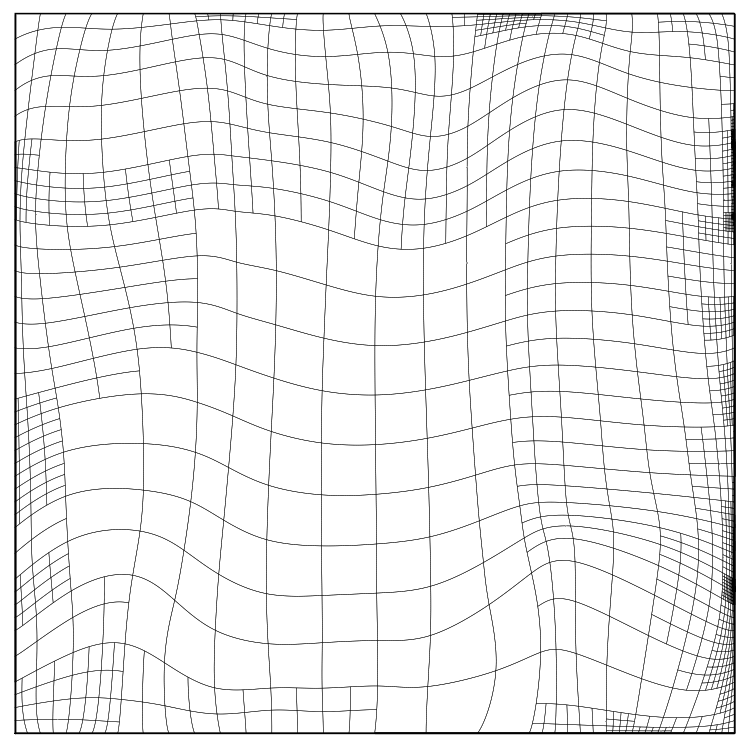}\label{fig:NRW_orth_domain_refined_reparamterized}} $\qquad$
  \subfloat[]{\includegraphics[align=c, width=.5\linewidth]{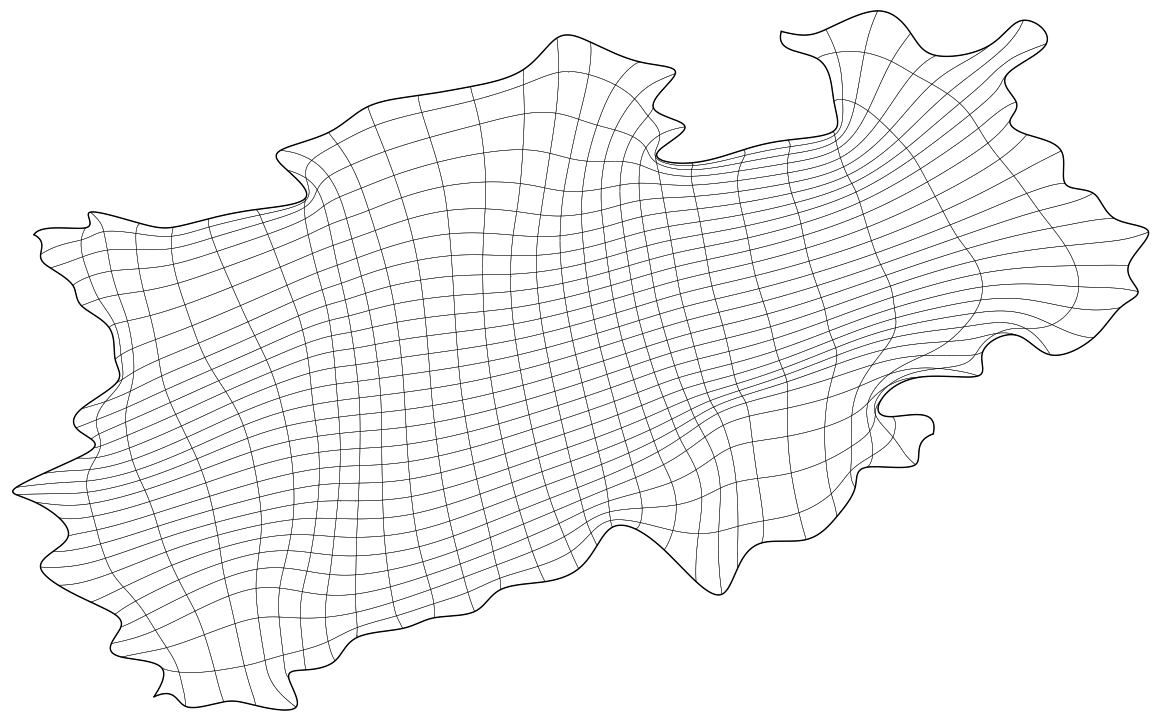}\label{fig:NRW_orth_linear_geometry}}
  \caption{Result of reparameterizing the reference parameterization of the German province of North Rhine-Westphalia (see Figure \ref{fig:NRW}), with $Q(\mathbf{s}) = Q_{\text{Orthogonality}}^{\mathbf{s}}(\mathbf{s})$. The reparameterized domain is shown in (a), while (b) shows the recomputed parameterization.}
  \label{fig:NRW_orth_linear}
\end{figure}

\noindent Given a set of abscissae $\boldsymbol{\Xi} = \{\boldsymbol{\xi}_1^c, \ldots \boldsymbol{\xi}_m^c \} \subset \mathbb{R}^2$, an alternative constraint $\mathbf{C}^{\boldsymbol{\Xi}}(\mathbf{s})$ follows from requiring that
\begin{align}
\label{eq:Adaptive_detJ_discrete_necessary}
    \boldsymbol{\epsilon}_i^L \leq \det \partial_{\boldsymbol{\xi}} \mathbf{s} (\boldsymbol{\xi}_i^c) \leq \boldsymbol{\epsilon}_i^U, \quad \forall i \in \{1, \ldots, m\},
\end{align}
where $\mathbb{R}^m \ni \boldsymbol{\epsilon}^{L, U} \geq 0$ are lower and upper thresholds. Note that (\ref{eq:Adaptive_detJ_discrete_necessary}) is nonlinear and nonconvex but not a sufficient condition for bijectivity of $\mathbf{s}$. However, it makes bijectivity likely for $m$ sufficiently large. \\
Finally, assuming that $\mathbf{s}$ is built from a structured basis $\left[ \mathcal{V}_h \right]$ resulting from a tensor product of the univariate bases
\begin{align*}
    \{ N_1^\square, \ldots, N_n^\square \} \quad \text{and} \quad \{ M_1^\square, \ldots, M_m^\square \},
\end{align*}
we may alternatively utilize the linear constraint proposed in \cite{xu2011parameterization}. Typically, we take $\left[ \mathcal{V}_h \right]$ as the cardinality-wise largest structured basis compatible with $\mathcal{T}$. Given
\begin{align}
    \mathbf{s}(\xi, \eta) = \sum_{i,j} \mathbf{c}_{i, j} N_i^\square(\xi) M_j^\square(\eta),
\end{align}
let the cones $\mathcal{C}^1(\mathbf{s})$ and $\mathcal{C}^2(\mathbf{s})$ be generated by the half rays $\mathbb{R}^+ \Delta_{i,j}^1$ and $\mathbb{R}^+ \Delta_{i, j}^2$ with
\begin{align*}
    \Delta_{i, j}^1 = \mathbf{c}_{i+1, j} - \mathbf{c}_{i, j} \quad \text{and} \quad \Delta_{i, j}^2 = \mathbf{c}_{i, j+1} - \mathbf{c}_{i, j},
\end{align*}
respectively. The constraint is based on the observation that if $\mathcal{C}^1(\mathbf{s})$ and $\mathcal{C}^2(\mathbf{s})$ only intersect in $\boldsymbol{\xi} = \mathbf{0}$, then $\mathbf{s}$ is bijective. In a direct optimization of $\mathbf{x}_h$, above constraint may be impractical since for most $\mathbf{x}_D$, the set 
\begin{align*}
\left \{ \mathbf{x}_h \in \mathcal{U}^{\mathbf{x}_D}_h \enskip \vert \enskip \mathcal{C}^1(\mathbf{x}_h) \cap \mathcal{C}^2(\mathbf{x}_h) = \{ \mathbf{0}\} \right \}
\end{align*}
is empty or the constraint is too restrictive. However, in the case of optimizing $\mathbf{s}$, for $\mathbf{s}^0 = \boldsymbol{\xi}$, the cones $\mathcal{C}^1(\mathbf{s}^0)$ and $\mathcal{C}^2(\mathbf{s}^0)$ are generated by $\mathbb{R}^+(1, 0)^T$ and $\mathbb{R}^+ (0, 1)^T$, respectively. A linear constraint $\mathbf{C}_L(\mathbf{s})$ follows from requiring that $\mathcal{C}^1(\mathbf{s})$ and $\mathcal{C}^2(\mathbf{s})$ be contained in the cones generated by 
\begin{align*}
     \mathbb{R}^+ \times \left \{(1, -1 + \epsilon)^T, (1, 1 + \epsilon)^T \right \} \quad \text{and} \quad  \mathbb{R}^+ \times \left \{(1, 1 + \epsilon)^T, (-1, 1 + \epsilon)^T \right \},
\end{align*}
respectively. Here $\epsilon \ll 1$ is a small positive parameter. Clearly, $\mathbf{s}^0$ is located exactly in the center of the feasible region (see Figure \ref{fig:linear_constraint}), making the constraint much less restrictive at the expense of having to compute $\mathbf{x}_h^*$ first.
\begin{remark}
We can combine the proposed constraints with the principles from Section \ref{sect:Adaptive_Eploiting_Maximum_Principle} to suppress overshoots due to extreme diffusive anisotropy. If $\mathbf{C}(\mathbf{s}) = \mathbf{C}_L(\mathbf{s})$, the problem remains convex.
\end{remark}
\noindent Figure \ref{fig:Indiana_AO_domain_reparameterized} shows the domain corresponding to the U.S. state of Indiana (see Figure \ref{fig:Indiana}) after optimizing with 
\begin{align*}
    Q = Q_{\text{AreaOrthogonality}}^{\mathbf{s}}
\end{align*}
under the constraint $\mathbf{C}(\mathbf{s}) = \hat{\mathbf{d}}(\mathbf{s})$ (see equation (\ref{eq:Adaptive_Bezier_Constraint})). The domain mapping $\mathbf{s}(\boldsymbol{\xi})$ is built from the same THB-basis as $\mathbf{x}^*_h$, comprised of $2338$ DOFs. Since Newton failed to converge, we recomputed $\mathbf{x}_h$ using the Picard approach, which converged after $21$ iterations. The result is depicted in Figure \ref{fig:Indiana_AO}. No a posteriori refinements were required. The reparameterization reduces the value of $L_{\text{AreaOrthogonality}}$ from the initial
\begin{align*}
L_{\text{AreaOrthogonality}}(\mathbf{x}_h^*) = 1.77 \times 10^2, \quad \text{to} \quad L_{\text{AreaOrthogonality}}(\mathbf{x}_h) = 1.36 \times 10^2.
\end{align*}
Next, we optimize the domain corresponding to the puzzle piece geometry (see Figure \ref{fig:solution_strategies_id}) with $\mathbf{C}(\mathbf{s}) = \mathbf{C}_L(\mathbf{s})$ and $Q = Q_{\text{Area}}^{\mathbf{s}}$. Hereby, $\mathbf{s}(\boldsymbol{\xi})$ is built from a structured spline space comprised of $646$ DOFs. The reparameterized domain is depicted in Figure \ref{fig:Puzzle_detJ_domain_linear_refined_reparameterized}.
Bijectivity of $\mathbf{x}_h$ is achieved with $2632$ DOFs and the resulting parameterization is depicted in Figure \ref{fig:Puzzle_detJ_linear_geometry}. With $L_\text{Area}(\mathbf{x}_h) = 142.710$, it is roughly as effective as the reparameterization from Figure \ref{fig:PoissonAreaPuzzle} with $k = 1$. \\
Figure \ref{fig:NRW_orth_linear} shows the German province of North Rhine-Westphalia upon reparameterization with 
\begin{align*}
Q = Q_\text{Orthogonality}^\mathbf{s},
\end{align*}
where $\mathbf{s}(\boldsymbol{\xi})$ is built from a structured spline space comprised of $578$ DOFs, with $\mathbf{C}(\mathbf{s}) = \mathbf{C}_L(\mathbf{s})$. Initially, \begin{align*}
    L_{\text{Orthogonality}}(\mathbf{x}_h^*) = 18.929, \quad \text{while} \quad L_{\text{Orthogonality}}(\mathbf{x}_h) = 5.160
\end{align*}
upon recomputation. Bijectivity is achieved with $4584$ DOFs, which is roughly double the initial $2724$ DOFs.

\begin{figure}[h!]
\centering
  \subfloat[]{\includegraphics[align=c, width=.3\linewidth]{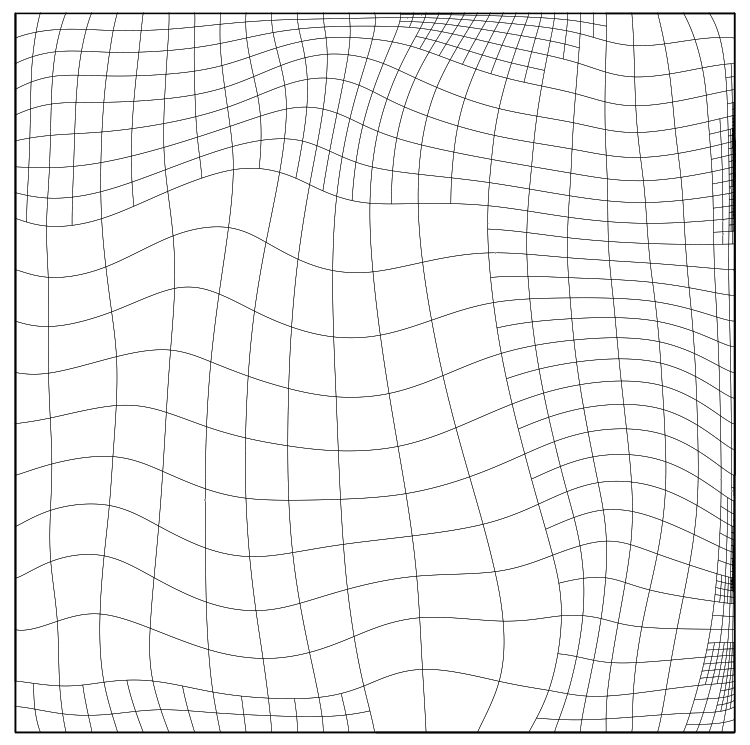}\label{fig:NRW_AO_domain_refined_reparamterized}} $\qquad$
  \subfloat[]{\includegraphics[align=c, width=.5\linewidth]{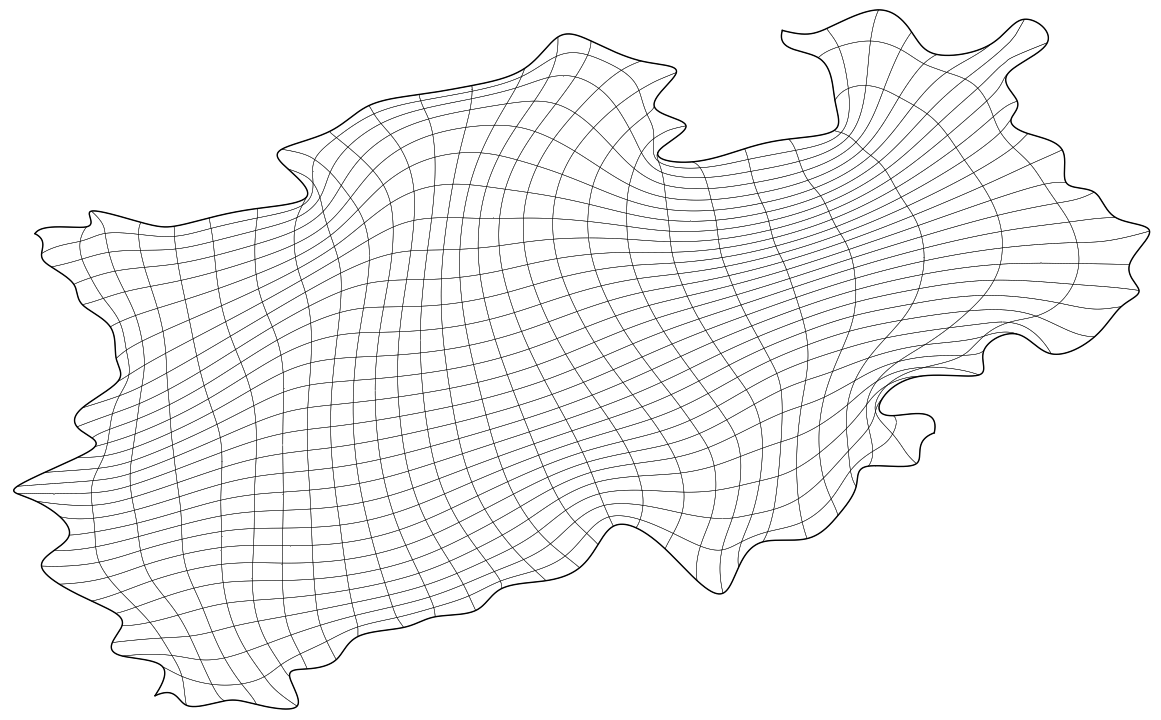}\label{fig:NRW_AO_linear_geometry}}
  \caption{Result of reparameterizing the reference parameterization of the German province of North Rhine-Westphalia (see Figure \ref{fig:NRW}), with $Q(\mathbf{s}) = Q_{\text{AreaOrthogonality}}^{\mathbf{s}}(\mathbf{s})$. The reparameterized domain is shown in (a), while (b) shows the recomputed parameterization.}
  \label{fig:NRW_AO_linear}
\end{figure}
Finally, Figure \ref{fig:NRW_AO_linear} shows the result of reparameterizing the same geometry with 
\begin{align*}
Q = Q^{\mathbf{s}}_\text{AreaOrthogonality}
\end{align*}
and the same constraints. Initially, 
\begin{align*}
    L_{\text{AreaOrthogonality}}(\mathbf{x}_h^*) = 51.244, \quad \text{while} \quad L_{\text{AreaOrthogonality}}(\mathbf{x}_h) = 30.896
\end{align*}
upon recomputation. Bijectivity is achieved with only $2928$ DOFs.

\subsection{Direct Optimization}
As an alternative to operating in the parametric domain, we may choose to directly optimize the geometry parameterization with respect to a quality cost function. As an advantage, we avoid the (possibly expensive) recomputation of $\mathbf{x}_h$. In order to avoid folding, constraints should be employed. As a disadvantage, the linear constraint $\mathbf{C}_L(\mathbf{x}_h)$ cannot be used and the initial guess $\mathbf{x}^{*}_h$ may fail to satisfy the conditions $\mathbf{\hat{d}}(\mathbf{x}_h) > \mathbf{0}$ (cf. (\ref{eq:Adaptive_Bezier_Constraint})) and $\mathbf{d}(\mathbf{x}_h) > \mathbf{0}$ (cf. (\ref{eq:Adaptive_Gravesen_Constraint})) despite being bijective. Heuristically, for complicated geometries, this is usually the case. In such cases, the only viable constraint is $\mathbf{C}^{\Xi}(\mathbf{x}_h)$ (cf. (\ref{eq:Adaptive_detJ_discrete_necessary})). \\
We optimize the puzzle piece geometry with $Q(\mathbf{x}_h) = Q_{\text{Area}}(\mathbf{x}_h)$ under the constraint $\mathbf{\hat{d}}(\mathbf{x}_h) > 0$, where the initial guess $\mathbf{x}_h^{*}$ is the parameterization from Figure \ref{fig:PoissonAreaPuzzle_k_0}.
\begin{figure}[h!]
\centering
  \subfloat[]{\includegraphics[align=c, width=.3\linewidth]{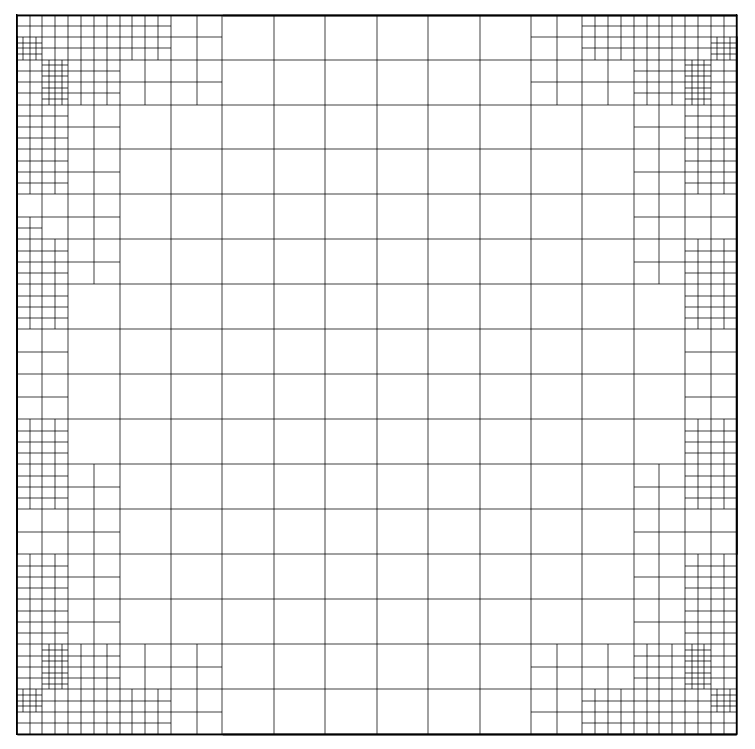}\label{fig:puzzle_piece_Area_direct_domain}} $\qquad$
  \subfloat[]{\includegraphics[align=c, width=.5\linewidth]{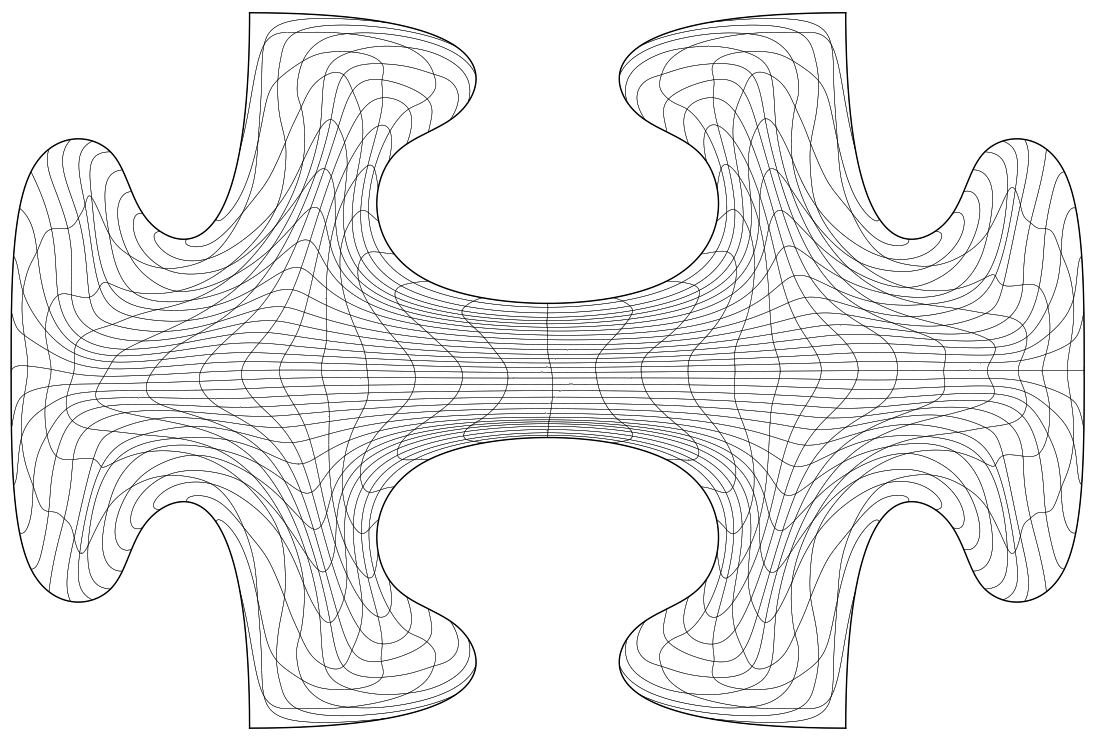}\label{fig:puzzle_piece_Area_direct_geometry}}
  \caption{The puzzle piece geometry after $21$ iterations of minimizing $Q(\mathbf{x}_h) = Q_{\text{Area}}$ under the constraint $\mathbf{\hat{d}}(\mathbf{x}_h) > 0$ (b) and the corresponding domain (a). The minimization was initialized with the parameterization from Figure \ref{fig:PoissonAreaPuzzle_k_0}.}
  \label{fig:puzzle_piece_Area_direct_30}
\end{figure}
Figure \ref{fig:puzzle_piece_Area_direct_30} shows the resulting parameterization. Convergence is achieved after $21$ constrained iterations. The reparameterization reduces $L_{\text{Area}}$ from the initial $L_{\text{Area}}(\mathbf{x}_h^*) = 3.29 \times 10^2$ to $L_{\text{Area}}(\mathbf{x}_h) = 0.96 \times 10^2$, which is slightly more pronounced than the reduction from Figure \ref{fig:PoissonAreaPuzzle} with $k=1.5$. However, the resulting parameterization is less regular compared to Figure \ref{fig:PoissonAreaPuzzle_k_15}, which can be remedied by adding a regularization of the form $Q(\mathbf{x}_h) = Q_{\text{Area}}(\mathbf{x}_h) + \beta Q_{\text{Uniformity}}(\mathbf{x}_h)$.
\begin{figure}[h!]
\centering
  \subfloat[]{\includegraphics[align=c, width=.4\linewidth]{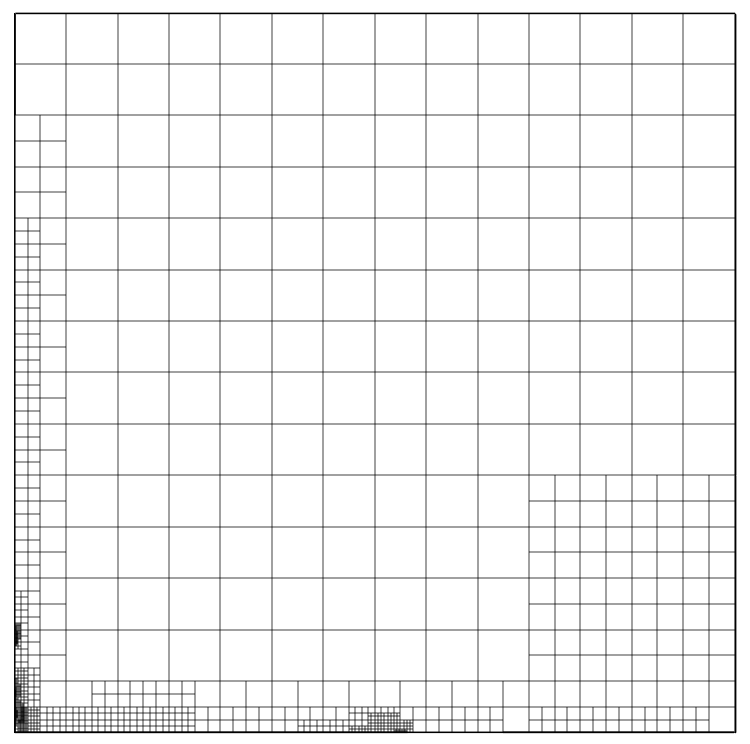}\label{fig:Indiana_Area_direct_domain}} $\qquad$
  \subfloat[]{\includegraphics[align=c, width=.4\linewidth]{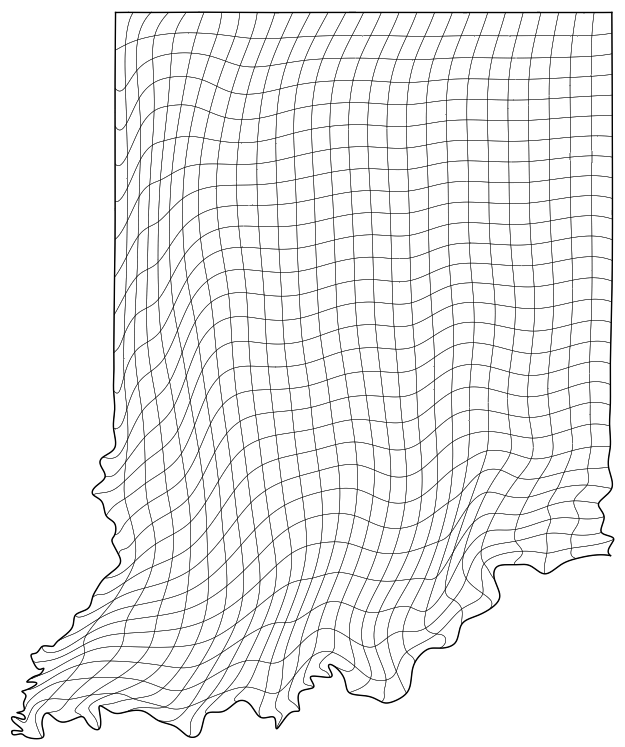}\label{fig:Indiana_Area_direct_geometry}}
  \caption{The parameterization of the U.S. state of Indiana after $30$ iterations of minimizing $Q(\mathbf{x}_h) = Q_{\text{Area}}$ (b) and the corresponding domain (a). The minimization was initialized with the parameterization from Figure \ref{fig:PoissonAreaIndiana_k_0}.}
  \label{fig:Indiana_Area_direct_40}
\end{figure}
Next, we optimize the U.S. state of Indiana with $Q(\mathbf{x}_h) = Q_{\text{Area}}(\mathbf{x}_h)$ under the constraint $\mathbf{C}^{\boldsymbol{\Xi}}(\mathbf{x}_h) \geq 0$ with
\begin{align}
\label{eq:Adaptive_Lower_Upper_puzzle_direct_orth}
    \boldsymbol{\epsilon}_i^L  = \alpha_L \times \det J(\mathbf{x}_h^*)(\boldsymbol{\xi}_i^c) \quad \text{and} \quad \boldsymbol{\epsilon}_i^U  = \alpha_U \times \det J(\mathbf{x}_h^*)(\boldsymbol{\xi}_i^c),
\end{align}
(see equation (\ref{eq:Adaptive_detJ_discrete_necessary})). \\
Figure \ref{fig:Indiana_Area_direct_40} shows the resulting parameterization after $30$ iterations. With 
\begin{align*}
    L_{\text{Area}}(\mathbf{x}_h^*) = 1.049 \times 10^2 \quad \text{and} \quad L_{\text{Area}}(\mathbf{x}_h) = 0.989 \times 10^2,
\end{align*}
the reduction is mild, yet somewhat more pronounced than in Figure \ref{fig:PoissonAreaIndiana}. Here, $\boldsymbol{\Xi}$ results from uniform sampling with $36$ points per element. The choice of the relaxation factors $0 \leq \alpha_L \leq 1$ and $1 \leq \alpha_U$ in (\ref{eq:Adaptive_Lower_Upper_puzzle_direct_orth}) tunes to which degree trading an increase in $L_\text{Area}$ for a decrease in the employed cost function is acceptable. Here, more conservative choices lead to less cost function reduction but to more uniform cell sizes and vice versa. Furthermore, values of $\alpha_L$ closer to $1$ prevent the grid from folding, even if fewer sampling points are used. We used $\alpha_L = 0.05$ and $\alpha_U = 4$.

\subsection{Achieving Boundary Orthogonality}
\begin{figure}[h!]
\centering
  \subfloat[]{\includegraphics[align=c, width=.4\linewidth]{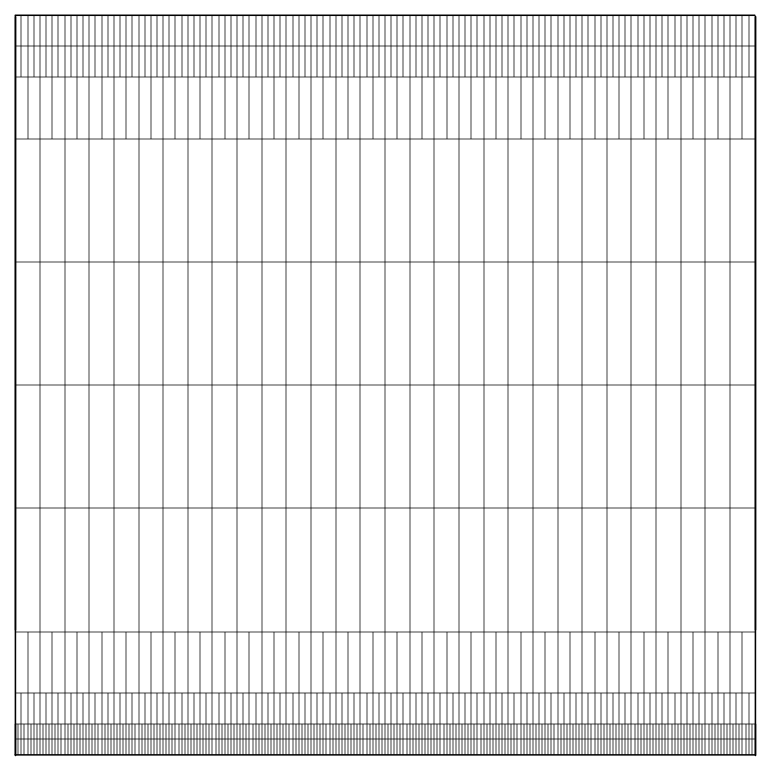}\label{fig:sin_bndry_orth_reference_domain}} $\qquad$
  \subfloat[]{\includegraphics[align=c, width=.5\linewidth]{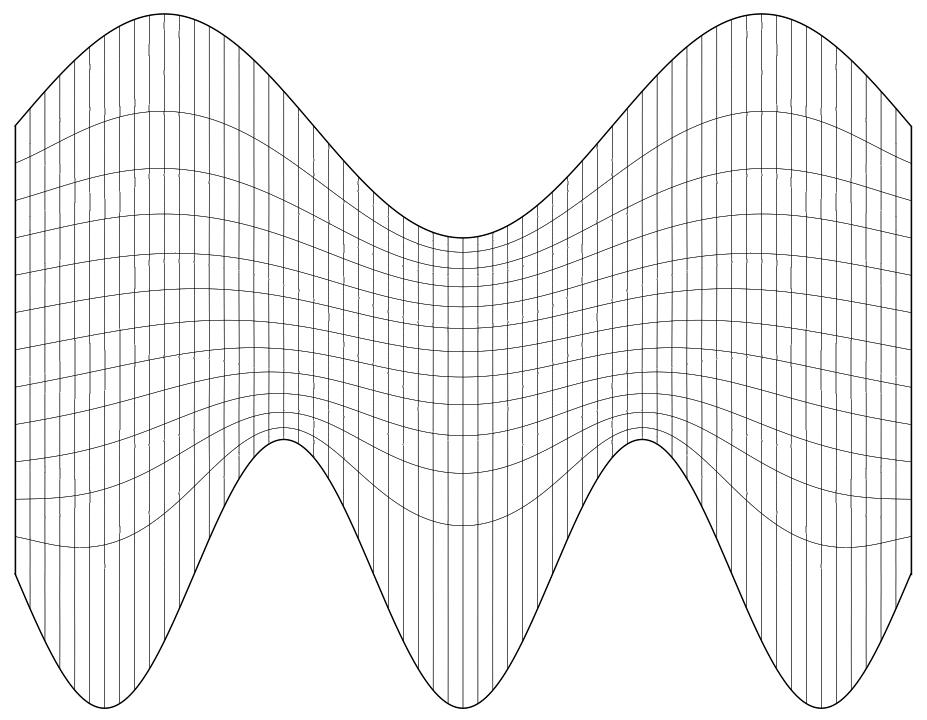}\label{fig:sin_bndry_orth_reference_optimized}}
  \caption{Reference parameterization of a tube-like shaped geometry which is to be orthogonalized by the northern and southern boundaries.}
  \label{fig:sin_bndry_orth_reference}
\end{figure}
\begin{figure}[h!]
\centering
  \subfloat[]{\includegraphics[align=c, width=.4\linewidth]{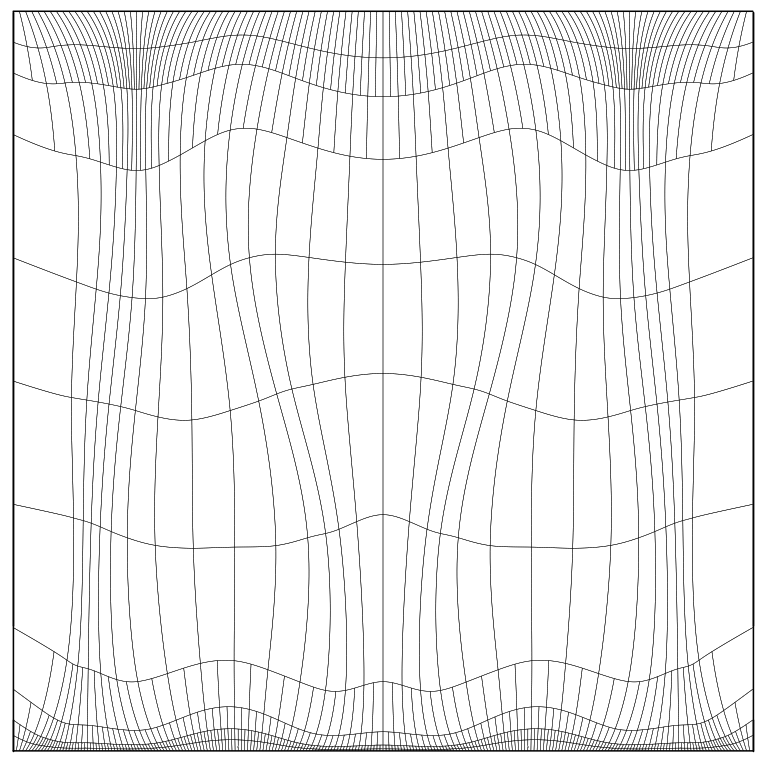}\label{fig:sin_bndry_orth_domain_refined_reparamterized}} $\qquad$
  \subfloat[]{\includegraphics[align=c, width=.5\linewidth]{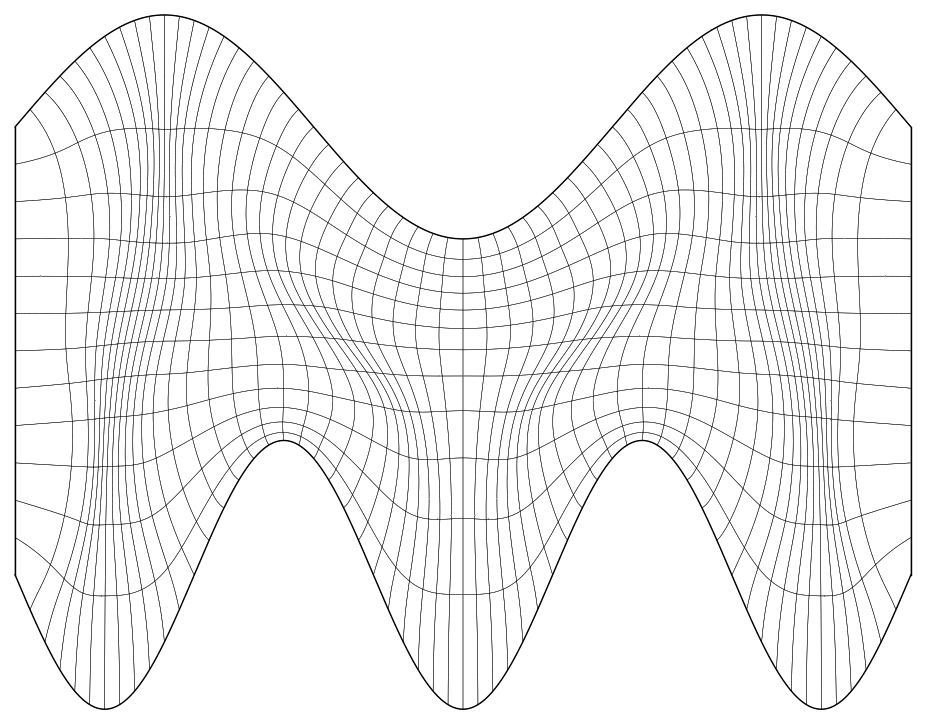}\label{fig:sin_bndry_orth_optimized}}
  \caption{Result of reparameterizing the geometry mapping from Figure \ref{fig:sin_bndry_orth_reference} by weakly enforcing boundary orthogonality through a large penalty term (b) and the corresponding reparameterized domain (a).}
  \label{fig:sin_bndry_orth}
\end{figure}

Many applications favor parameterizations with isolines that are orthogonal to the boundary contours. One way to achieve this is allowing $\lambda_i^{\mathbf{s}} = \lambda_i(\boldsymbol{\xi})^{\mathbf{s}}$ in (\ref{eq:Adaptive_Domain_Q}) and taking $\lambda_\text{Orthogonality}^{\mathbf{s}}$ large close to $\partial \hat{\Omega}$. We are considering the example of achieving orthogonality at the northern and southern boundaries of the geometry depicted in Figure \ref{fig:sin_bndry_orth_reference}. To this end, we minimize the cost function \begin{align*}
    Q(\mathbf{s}) = (1 + \lambda_{\text{O}}(\boldsymbol{\xi})) Q^{\mathbf{s}}_{\text{Orthogonality}},
\end{align*}
where $\lambda_{\text{O}}(\boldsymbol{\xi})$ takes on large values close to the northern and southern boundaries of $\partial \hat{\Omega}$. We employ the constraint $\mathbf{C}(\mathbf{s}) = \mathbf{C}_L(\mathbf{s})$, where $\mathbf{s}(\boldsymbol{\xi})$ is built from a structured spline space comprised of $594$ DOFs. The resulting parameterization is depicted in Figure \ref{fig:sin_bndry_orth}.
The figure indeed shows a large degree of orthogonalization, which is somewhat weaker in the protruded parts of the geometry. This is due to orthogonality only being enforced \textit{weakly} through a penalty term. More pronounced boundary orthogonalization may be achieved by taking $\lambda_\text{O}$ larger close to $\partial \hat{\Omega}$. \\
Let $\gamma_e, \gamma_w, \gamma_s$ and $\gamma_n$ refer to the eastern, western, southern and northern parts of $\partial \hat{\Omega}$, respectively. For a more drastic boundary orthogonalization, we follow the approach from \cite[Chapter ~6]{thompson1998handbook}, which consists of solving the problem
\begin{align}
\label{eq:Adaptive_boundary_orth_laplace}
    \Delta_{\mathbf{x}^*_h} f = 0 \quad \text{s.t.} \quad f = 0 \text{ on } \gamma_e, \enskip f = 1 \text{ on } \gamma_w \quad \text{and} \quad \frac{\partial f}{\partial \mathbf{n}} = 0 \text{ on } \gamma_s \cup \gamma_n
\end{align}
on an initially folding-free geometry parameterization $\mathbf{x}_h^*$. Here, $\mathbf{n}$ denotes the unit outward normal vector on $\partial \Omega$. Upon completion, the control mapping $\mathbf{s} = (\mathbf{s}_1, \mathbf{s}_2)^T \equiv (s, t)^T$ is computed from
\begin{align}
\label{eq:Adaptive_control_mapping_boundary_orth_PDE}
    s(\xi, \eta) = f(\xi, 0) H_0(\eta) + f(\xi, 1) H_1(\eta) \quad \text{and} \quad t(\xi, \eta) = \eta,
\end{align}
where
\begin{align}
    H_0(\eta) = (1 + 2 \eta)(1 - \eta)^2 \quad \text{and} \quad H_1(\eta) = (3 - 2 \eta) \eta^2
\end{align}
are cubic Hermite interpolation functions. It can be shown that with this choice of $s$ and $t$, the solution of (\ref{eq:Adaptive_Strong_Equations_reparameterized}) is orthogonal at $\gamma_s$ and $\gamma_n$.  We approximately solve for $f$ by computing the solution $f_h$ of the discretized counterpart of (\ref{eq:Adaptive_boundary_orth_laplace}) over some structured spline space $\mathcal{V}_h$. Hereby, the Neumann boundary conditions are weakly imposed through partial integration. The control mapping follows from replacing $f \rightarrow f_h$ in (\ref{eq:Adaptive_control_mapping_boundary_orth_PDE}). Should orthogonality at $\gamma_w$ and $\gamma_e$ be desired, we simply exchange the roles of $s \rightarrow t$, $(\gamma_s, \gamma_n) \rightarrow (\gamma_w, \gamma_e)$ and $\xi \rightarrow \eta$.
\begin{remark}
Unlike $f$, $f_h$ may fail to be monotone increasing on $\gamma_s$ or $\gamma_n$, leading to a folded control mapping $\mathbf{s}(\boldsymbol{\xi})$.
\end{remark}
\begin{figure}[h!]
\centering
  \subfloat[]{\includegraphics[align=c, width=.4\linewidth]{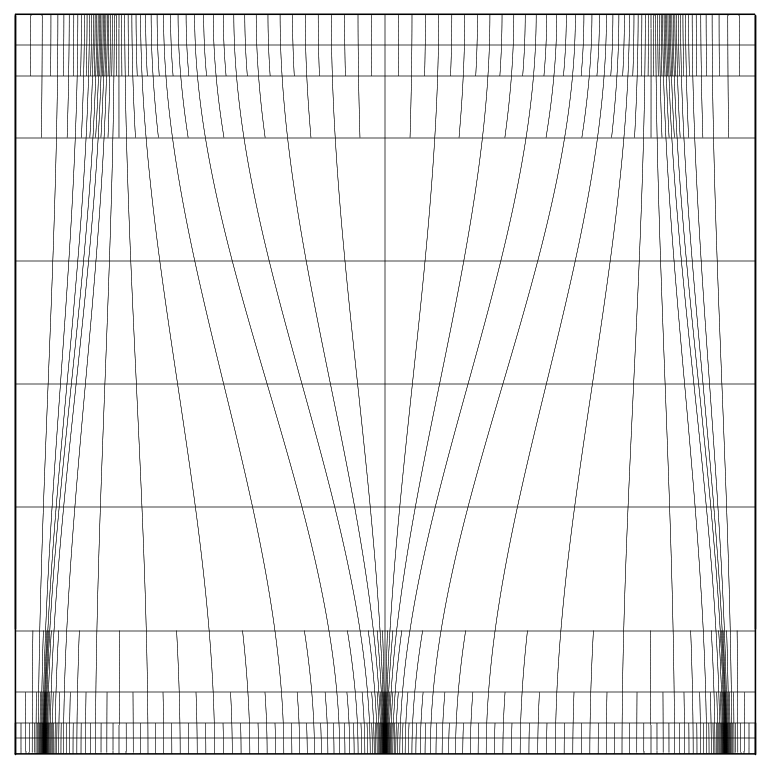}\label{fig:sin_bndry_orth_PDE_domain_refined_reparamterized}} $\qquad$
  \subfloat[]{\includegraphics[align=c, width=.5\linewidth]{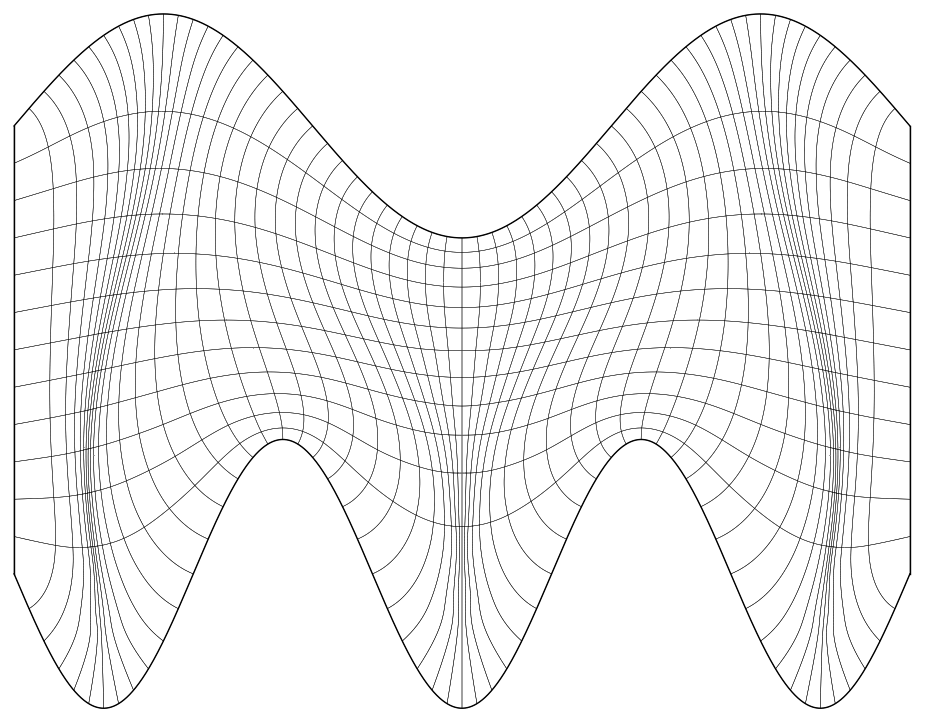}\label{fig:sin_bndry_orth_PDE_optimized}}
  \caption{Result of reprameterizing the geometry mapping from Figure \ref{fig:sin_bndry_orth_reference} using the approach proposed in \cite[Chapter ~6]{thompson1998handbook} (b) and the corresponding reparameterized domain (a).}
  \label{fig:sin_bndry_orth_PDE}
\end{figure}
Figure \ref{fig:sin_bndry_orth_PDE} shows the recomputed parameterization of the same geometry using the preceding methodology, along with the reparameterized parametric domain, which has been computed from the same structured spline basis as in Figure (\ref{fig:sin_bndry_orth}). The figure shows an outstanding boundary orthogonalization, which comes at the expense of larger elements in the protruded parts compared to Figure \ref{fig:sin_bndry_orth}. We introduce another control mapping $\mathbf{s}^\prime(\boldsymbol{\xi})$, which we compute from the solution of
\begin{align}
\label{eq:Adaptive_minimization_s_prime_orth}
    \int_{\hat{\Omega}} (\det \partial_{\mathbf{s}} \mathbf{x}_h)^k (g_{11} + \beta g_{22})_{\mathbf{s} \rightarrow \mathbf{s}^\prime} \det \partial_{\boldsymbol{\xi}} \mathbf{s} \mathrm{d} S \rightarrow \min \limits_{\mathbf{s}^\prime \in \mathcal{V}_h^2}, \quad \text{s.t.} \quad \mathbf{s}^\prime(\boldsymbol{\xi}) = \mathbf{s}(\boldsymbol{\xi}) \text{ on } \partial \hat{\Omega},
\end{align}
where $\mathbf{s} = (s, t)^T$ and $\mathbf{x}_h$ correspond to Figures \ref{fig:sin_bndry_orth_PDE} (a) and (b), respectively. Here, the $g_{ii}$ correspond to diagonal entries of the metric tensor associated with the diffeomorphism between $\mathbf{s} \vert_{\hat{\Omega}}$ and $\mathbf{s}^\prime \vert_{\hat{\Omega}}$. As before, $k > 0$ tunes to which degree the spread in cell size is penalized, while $\beta > 1$ tunes the degree to which $\mathbf{s}^\prime$ is contracted / expanded in the direction of $\partial_{\eta} \mathbf{s}$, in order to compensate for large / small cells in $\mathbf{x}_h$. Taking $\beta$ large essentially freezes $\mathbf{s}^\prime$ in the direction of $\partial_{\xi} \mathbf{s}$, such that boundary orthogonality is preserved. Note that in (\ref{eq:Adaptive_minimization_s_prime_orth}), we are essentially solving the discrete counterpart of
\begin{align}
    \nabla_{\mathbf{s}} \cdot( D \nabla_\mathbf{s} \mathbf{s}^\prime_i ) = 0, \quad i \in \{1, 2\}, \quad \text{s.t.} \quad \mathbf{s}^\prime(\boldsymbol{\xi}) = \mathbf{s}(\boldsymbol{\xi}), \quad \text{with} \quad D = (\det \partial_{\mathbf{s}} \mathbf{x}_h)^k \begin{pmatrix} 1 & 0 \\ 0 & \beta \end{pmatrix}.
\end{align}
\begin{figure}[h!]
\centering
  \subfloat[]{\includegraphics[align=c, width=.4\linewidth]{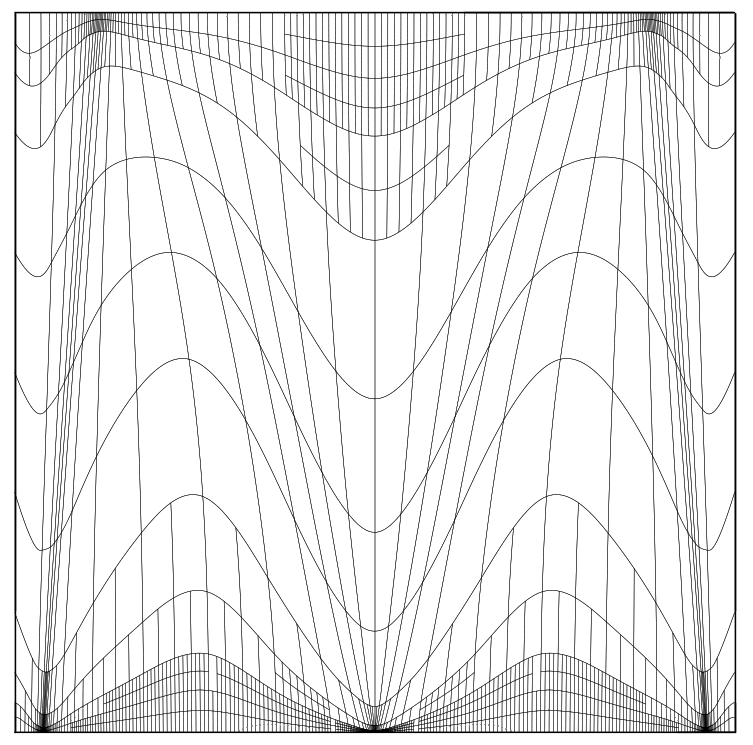}\label{fig:sin_bndry_orth_PDE_rep_domain_refined_reparamterized}} $\qquad$
  \subfloat[]{\includegraphics[align=c, width=.5\linewidth]{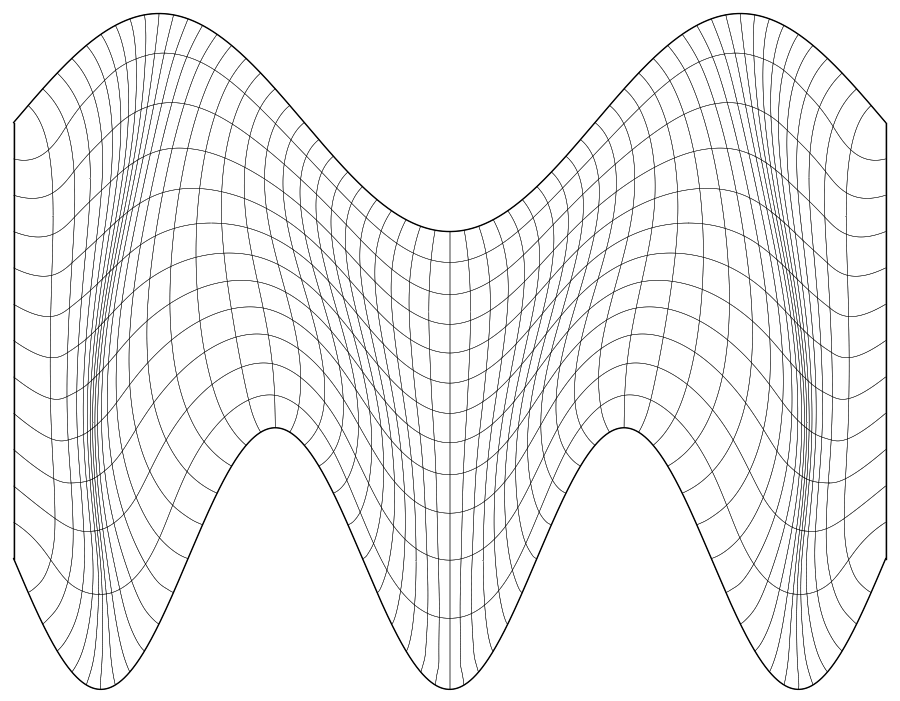}\label{fig:sin_bndry_orth_PDE__rep_optimized}}
  \caption{Result of reprameterizing the geometry mapping from Figure \ref{fig:sin_bndry_orth_PDE} using the principles from Section \ref{sect:Adaptive_Eploiting_Maximum_Principle} (b) and the corresponding reparameterized domain (a).}
  \label{fig:sin_bndry_orth_PDE_rep}
\end{figure}
Figure \ref{fig:sin_bndry_orth_PDE_rep} shows the geometry parameterization along with the reparameterized domain upon recomputation with $k = 0.75$ and $\beta = 300$. Compared to Figure \ref{fig:sin_bndry_orth_PDE}, the figure shows a much better cell size distribution, in particular close to the boundaries. Large cells can be further penalized by increasing the value of $k$.
\section{Conclusion}
In this work, we presented a goal-oriented adaptive THB-spline framework for PDE-based planar parameterization. For this, we adopted the a posteriori refinement technique of dual weighted residual and proposed several goal-oriented refinement cost functions. This resulted in numerical schemes that combine iterative solution techniques with THB-enabled local a posteriori refinement strategies, hence avoiding over-refinement in computing a folding-free geometry parameterization. \\
In order to fine-tune the parametric properties of the resulting mapping, we combined aforementioned schemes with the concept of domain optimization. Hereby, the (convex) parametric domain, which constitutes the target domain of the mapping inverse, is reparameterized in order to alter the parametric properties of the recomputed mapping. For this, we proposed several optimization constraints that avoid the loss of bijectivity.

\section*{Acknowledgements}
The authors gratefully acknowledge the research funding which was partly provided by the MOTOR project that has received funding from the European Unions Horizon 2020 research and innovation program under grant agreement No 678727. \\
Furthermore, the authors are grateful for the coding help they received from the Nutils core development team.

\bibliographystyle{elsarticle-num}
\bibliography{bibliography}


\end{document}